\documentclass[11pt]{article}
\usepackage{amsmath}
\usepackage[colorlinks=true]{hyperref}
 \hypersetup{urlcolor=blue, citecolor=blue}

\usepackage{amssymb}
\usepackage{mathdots}
\makeatletter
\renewcommand{\@begintheorem}[2]{\begin{trivlist}
\item[\hspace{\labelsep}{\bf \mbox{~~~}#1\ #2.}]}
\renewcommand{\@opargbegintheorem}[3]{\begin{trivlist}
\item[\hspace{\labelsep}{\bf \mbox{~~~}#1\ #2 {\rm (#3).}}]}
\renewcommand{\@endtheorem}{\end{trivlist}}
\makeatother

\newtheorem{theorem}{Theorem}
\newtheorem{corollary}{Corollary}
\newtheorem{lemma}{Lemma}
\newtheorem{proposition}{Proposition}
\newtheorem{conjecture}{Conjecture}
\newtheorem{problem}{Problem}
\newtheorem{definition}{Definition}
\newtheorem{remark}{Remark}
\newtheorem{notation}{Notation}
\newtheorem{example}{Example}

\if0
\newtheorem{theorem}{Theorem}[section]
\newtheorem{corollary}[theorem]{Corollary}

\newtheorem{lemma}[theorem]{Lemma}

\newtheorem{definition}[theorem]{Definition}
\newtheorem{remark}[theorem]{Remark}

\fi

\begin{document}

\begin{center}{\Large \noindent \bf Arbitrary matrix coefficient assignment for block matrix linear control systems by static output feedback}
\end{center}

\bigskip

\centerline{\scshape Vasilii Zaitsev}
\medskip
{\footnotesize
 \centerline{Udmurt State University}
      \centerline{426034 Izhevsk, Russia, e-mail: verba@udm.ru}
}

\medskip

\centerline{\scshape Inna Kim}
\medskip
{\footnotesize
 \centerline{Udmurt State University}
      \centerline{426034 Izhevsk, Russia, e-mail: kimingeral@gmail.com}
}

\bigskip

\begin{abstract}
This work has introduced a generalized formulation of the problem of eigenvalue spectrum assignment for block matrix systems. In this problem, it is required to construct a feedback that provides that the matrix of the closed-loop system is similar to a block companion matrix with arbitrary predetermined block matrix coefficients. Sufficient conditions for the resolvability of this problem by linear static output feedback are obtained when the coefficients of the system have a special form, namely, the state matrix is a lower block Frobenius matrix or a lower block Hessenberg matrix, and the input and output block matrix coefficients contain some zero blocks. These conditions are controllability-like rank conditions. Sufficient conditions are constructive. It is proved that, in particular cases, when the system has block scalar matrix coefficients, these conditions can be weakened. The results generalize the previous results obtained for the case of one-dimensional blocks and for the case of systems given by a linear differential equation of higher orders with a multidimensional state. Based on the main results, algorithms are developed that ensure the construction of a gain matrix. The algorithms are implemented on a modeling  examples.
\end{abstract}


\medskip


\section{Introduction \label{sect1}}

Consider a linear control system
\begin{equation}\label{026}
\dot z=\mathcal{F}z+\mathcal{G}v.
\end{equation}
Here $z\in\mathbb{K}^r$ is a state, $v\in\mathbb{K}^q$ is a control, $\mathbb{K}=\mathbb{R}$ or $\mathbb{K}=\mathbb{C}$, $\mathcal{F}\in M_r(\mathbb{K})$, $\mathcal{G}\in M_{r,q}(\mathbb{K})$, $M_{r,q}(\mathbb{K})$ is a space of $r\times q$-matrices with elements of $\mathbb{K}$, $M_r(\mathbb{K}):=M_{r,r}(\mathbb{K})$ (we will denote $M_{r,q}:=M_{r,q}(\mathbb{K})$, $M_r:=M_r(\mathbb{K})$, if the set $\mathbb{K}$ is predefined).

Suppose that the controller is constructed as linear static state feedback (LSSF)
\begin{equation}\label{027}
v=\mathcal{L}z
\end{equation}
with $\mathcal{L}\in M_{q,r}(\mathbb{K})$. The closed-loop system has the form
\begin{equation}\label{027-1}
\dot z=(\mathcal{F}+\mathcal{G}\mathcal{L})z.
\end{equation}
Denote by
\begin{equation}\label{028}
\chi(\mathcal{F}+\mathcal{G}\mathcal{L},\lambda)= \lambda^r+\delta_{1}\lambda^{r-1}+ \ldots  + \delta_{r}
\end{equation}
the characteristic polynomial of the matrix of the closed-loop system \eqref{027-1}. Let $\lambda_j\in\mathbb{C}$, $j=\overline{1,n}$, be the roots of the polynomial \eqref{028}. They form the spectrum of the system \eqref{027-1}.

In the theory of modal control \cite{SIMON1968316}, one can distinguish the problem of assigning coefficients  $(\delta_1,\ldots,\delta_n)$ of the characteristic polynomial \eqref{028}  and the problem of spectrum assignment (pole assignment) \cite{Wonham_1967}. In the spectrum assignment problem, for an arbitrary (admissible) set $\{\lambda_1,\ldots,\lambda_n\}$, it is required to construct a control \eqref{027} such that the spectrum of the closed-loop system \eqref{027-1} coincides with the given set.

If $\mathbb{K}=\mathbb{C}$ then there is a bijection $(\delta_1,\ldots,\delta_n)\in\mathbb{C}^n$ $\longleftrightarrow$ $\sigma=\{\lambda_1,\ldots,\lambda_n\}$ $(\lambda_j\in\mathbb{C})$ between the (ordered) set of coefficients of the characteristic polynomial \eqref{028} and the set of roots $\sigma=\{\lambda_1,\ldots,\lambda_n\}$ (that is, the spectrum) of this polynomial. 
If $\mathbb{K}=\mathbb{R}$, then such a bijection $(\delta_1,\ldots,\delta_n)\in\mathbb{R}^n$ $\longleftrightarrow$ $\sigma=\{\lambda_1,\ldots,\lambda_n\}$ $(\lambda_j\in\mathbb{C})$ exists if $\sigma$ is a set of real type, i.e., $\overline{\sigma}=\sigma$ (that is, the set $\sigma$ is invariant under the complex conjugation operation). Thus, for the system \eqref{027-1} the problem of assigning coefficients is equivalent to the spectrum assignment problem.

It is said that the spectrum assignment problem is resolvable for system \eqref{026} by LSSF if for  any $\delta_i\in\mathbb{K}$, $i=\overline{1,r}$, there exists a controller \eqref{027}
such that the characteristic polynomial of the matrix of the closed-loop system \eqref{027-1}
satisfies the condition \eqref{028}.

The condition  of complete controllability of system \eqref{026}, i.e., 
${\rm rank}\,[\mathcal{G},\mathcal{F}\mathcal{G},\ldots,\mathcal{F}^{r-1}\mathcal{G}]=r$,
is necessary and sufficient for the resolvability of the spectrum assignment problem for system  \eqref{026} by LSSF (see \cite{Popov_1964} for $\mathbb{K}=\mathbb{C}$ and \cite{Wonham_1967} for $\mathbb{K}=\mathbb{R}$).

Now let us consider the  multidimensional differential equation. Let  $s\in\mathbb{N}$ be given. Consider a linear control system
\begin{equation}\label{040}
x^{(n)}+A_{1}x^{(n-1)}+ \ldots  + A_{n}x  = B_1 u.
\end{equation}
Here $x\in\mathbb{K}^s$ is a state vector, $u\in\mathbb{K}^s$ is a control vector, $A_i\in M_s(\mathbb{K})$, $i=\overline{1,n}$, $B_1 \in M_s(\mathbb{K})$ are constant matrices.

Suppose that $\det B_1\ne 0$. Then, for any $\Gamma_i\in M_s(\mathbb{K})$, $i=\overline{1,n}$, one can choose  a linear state feedback control 
\begin{equation}\label{050}
u=K_1 x^{(n-1)}+\ldots+K_n x,
\end{equation}
where $K_i\in M_s(\mathbb{K})$, such that the closed-loop system has the form
\begin{equation}\label{055}
x^{(n)}+\Gamma_{1}x^{(n-1)}+ \ldots  + \Gamma_{n}x  = 0.
\end{equation}
For this purpose, one should take  
\begin{equation}\label{056}
K_i=B_1^{-1}(A_i-\Gamma_i), \quad i=\overline{1,n}.
\end{equation}
This means that one can assign arbitrary matrix coefficients to the differential equation \eqref{055}. 

In the paper \cite{Zaitsev_2021}, this property is named as {\it arbitrary matrix eigenvalue spectrum assignability (by linear static state feedback)} \cite[Definition~1]{Zaitsev_2021}, and the condition $\det B_1\ne 0$ is necessary and sufficient for this property to hold \cite[Proposition~1]{Zaitsev_2021}.
In fact, the name of this property is not quite correct. It is more correct to say that we can {\it assign arbitrary matrix coefficients to the  characteristic matrix polynomial}.

Consider the characteristic matrix polynomial of the system \eqref{055}:
\begin{equation}\label{057}
\Psi(\lambda)=I \lambda^n + \Gamma_{1} \lambda^{n-1} + \ldots + \Gamma_{n-1} \lambda + \Gamma_{n}, \quad I\in M_s.
\end{equation}

In contrast to the scalar case, the problem of assigning arbitrary matrix coefficients to the  characteristic matrix polynomial \eqref{057} and the matrix spectrum assignment problem for system \eqref{055} are not equivalent due to a number of reasons. The question of the roots of a matrix polynomial is not simple and is not unambiguous \cite{Lancaster-Tismenetsky,Dennis_1976,Tsai_1988}. Firstly, for matrix polynomials $\Psi(\lambda)$ (unlike scalar polynomials), the fundamental  theorem of algebra that every polynomial has at least one root does not hold. Example: $n=2$, $s=2$, $X,J\in M_2(\mathbb{C})$; $J=\left[\begin{matrix} 0 & 1 \\ 0 & 0 \end{matrix}\right]$; $P(X)=X^2 -J$. There is no matrix $X$ (neither real nor complex) such that $X^2=J$. Next, a matrix equation can have an infinite continuum set of solutions. 
Example: $n=2$, $s=2$, $X,I\in M_2(\mathbb{C})$, $I=\left[\begin{matrix} 1 & 0 \\ 0 & 1 \end{matrix}\right]$. The equation $X^2 -I=0$ has partial solutions $X(\tau)=\left[\begin{matrix} \cos\tau & \sin\tau \\ \sin\tau & -\cos\tau \end{matrix}\right]$, $\tau\in\mathbb{R}$. Next, under certain conditions, there are complete sets of distinct left solvents $(L_j)$ and right solvents $(R_j)$ satisfying the following matrix
equations, respectively:
\begin{equation*}
L_j^n+L_j^{n-1} \Gamma_1+ L_j^{n-2} \Gamma_2+\ldots + L_j \Gamma_{n-1} + \Gamma_{n}=0, \quad j=1,\ldots,n,
\end{equation*}
and
\begin{equation*}
R_j^n+ \Gamma_1 R_j^{n-1}+ \Gamma_2 R_j^{n-2}+\ldots +  \Gamma_{n-1} R_j + \Gamma_{n}=0, \quad j=1,\ldots,n;
\end{equation*}
the spectral factorization of the matrix polynomial $\Psi(\lambda)$ can be described by
\begin{equation*}
\Psi(\lambda)=(\lambda I -S_1)(\lambda I -S_2)\cdots(\lambda I -S_n), \quad I,S_j\in M_s(\mathbb{C});
\end{equation*}
the relationships among the left solvents $(L_j)$, the right solvents $(R_j)$ and the spectral factors $(S_j)$ of $\Psi(\lambda)$ can be found in \cite{Shieh_1981}; and in general, $S_i\ne L_j$, $S_i\ne R_j$, $L_i \ne R_j$. Thus, the matrix spectrum assignment problem for system \eqref{055}  can be posed only under certain restrictions, while the problem of assigning arbitrary matrix coefficients to the characteristic matrix polynomial \eqref{057}, for equations like \eqref{040}, can be posed without conditions, and such a problem, in fact, was solved in the paper \cite{Zaitsev_2021}, for  multidimensional differential equation by means of static output feedback.

Standard change of variables $z_1=x$, $z_2=x'$, \ldots, $z_n=x^{(n-1)}$ reduces the control system \eqref{040}, \eqref{050} to the form \eqref{026}, \eqref{027} where $z={\rm col}[z_1,\ldots,z_n]\in\mathbb{K}^{ns}$, $v=u\in\mathbb{K}^{s}$,
\begin{equation}\label{062}
\begin{gathered}
\mathcal{F}=\left[\begin{matrix}
0 & I & 0 & \ldots & 0 \\
0 & 0 & I & \ldots & 0 \\
\vdots & \vdots & \vdots & \ddots & \vdots  \\
0 & 0 & 0 & \ldots & I \\
-A_n & -A_{n-1} & -A_{n-2} & \ldots & -A_1
\end{matrix} \right], \qquad
\mathcal{G}=\left[\begin{matrix} 0 \\ 0 \\ \vdots \\ 0 \\ B_1 \end{matrix} \right],\\
\mathcal{L}=[K_n,K_{n-1},\ldots,K_1], \qquad \qquad r=ns,  \quad q=s,
\end{gathered}
\end{equation}
$0\in M_s$, $I\in M_s$. The state feedback controller \eqref{027} with coefficients \eqref{056} of the matrix $\mathcal{L}$ leads system \eqref{026}, \eqref{062} to the closed-loop system \eqref{027-1} with the matrix $\mathcal{F}+\mathcal{G}\mathcal{L}=:\Phi$. The matrix $\Phi$ has the block form
\begin{equation}\label{064}
\Phi=\left[\begin{matrix}
0 & I & 0 & \ldots & 0 \\
0 & 0 & I & \ldots & 0 \\
\vdots & \vdots & \vdots & \ddots & \vdots  \\
0 & 0 & 0 & \ldots & I \\
-\Gamma_n & -\Gamma_{n-1} & -\Gamma_{n-2} & \ldots & -\Gamma_1
\end{matrix} \right].
\end{equation}
This matrix is a {\it block companion matrix} associated to the matrix polynomial $\Psi(\Lambda)$, $\Lambda\in M_s$ \cite{Martins_2014}.
For arbitrary block matrices, the formal determinant and the characteristic (matrix) polynomial are defined only under the condition that the blocks of the matrix (at least those located in different block rows and block columns) are commuting in pairs \cite{Vitoria_1982,Martins_2014}. However, for a block matrix of the form \eqref{064}, this condition is satisfied and the formal characteristic (matrix) polynomial $\Delta(\Lambda)\!:=\det \big((I\otimes \Lambda) - \Phi\big)$ coincides with $\Psi(\Lambda)$. Thus, a matrix of the form \eqref{064} is, in a certain sense, the canonical form for block matrices.

Based on the above, let us consider the corresponding problem for a block matrix system not of a special form \eqref{062}, but of a general form. Consider the following block matrix system with $s\times s$-blocks:
\begin{gather} \label{066}
\dot{x}=F x + G u ,\quad x\in\mathbb{K}^{ns}, \quad u\in\mathbb{K}^{ms},\\
\label{067}
{F}=\begin{bmatrix}
{F}_{11} &  \ldots & {F}_{1n} \\
\vdots &  & \vdots  \\
{F}_{n1} &  \ldots & {F}_{nn} 
\end{bmatrix}, \;
{G}=\left[\begin{matrix} 
{G}_{11} &  \ldots & {G}_{1m} \\
\vdots &  & \vdots  \\
{G}_{n1} &  \ldots & {G}_{nm} 
\end{matrix} \right], \; {F}_{ij},{G}_{j \alpha}\in M_s(\mathbb{K}), \; i,j=\overline{1,n}, \; \alpha=\overline{1,m}. 
\end{gather} 
Here $x\in\mathbb{K}^{ns}$ is a state vector, $u\in\mathbb{K}^{ms}$ is a control vector.
Let the control in system~\eqref{066}, \eqref{067} have the form of linear static state feedback:
\begin{equation} \label{068}
{u}={L}{x}.
\end{equation}
The closed-loop system has the form
\begin{equation} \label{069}
\dot{x}=(F  + G L )x.
\end{equation}
We give the following definition.
\begin{definition}\label{def-01} 
We say that, for system \eqref{066}, \eqref{067}, the problem of {\it arbitrary matrix coefficient assignment (AMCA) for the characteristic matrix polynomial (CMP) by linear static state feedback (LSSF)  is resolvable} if for any $\Gamma_i\in M_s(\mathbb{K})$, $i=\overline{1,n}$, there exists a gain matrix $L\in M_{ms,ns}(\mathbb{K})$  such that the closed-loop system \eqref{069} is reducible by some non-degenerate change of variables $z=Sx$ to the system
\begin{equation}\label{070}
\dot z=\Phi z, \quad z\in\mathbb{K}^{ns},
\end{equation}
with the matrix $\Phi$ of \eqref{064}, that is the matrix $F+GL$ of the system \eqref{069} is similar to the matrix \eqref{064}: $S(F+GL)S^{-1}=\Phi$.
\end{definition}

This definition is an extension of the (corrected) Definition 1 of \cite{Zaitsev_2021} from systems \eqref{040}, \eqref{050} to block matrix systems \eqref{066}, \eqref{067}, \eqref{068}. 

\begin{remark}
\label{rem-dop-1} Systems of the form \eqref{040} of the second order are found in spectral problems for vibrating systems, see,  e.g., \cite{Lancaster_2005,Lancaster_2005-2,Mottershead_2006,Lancaster_2007,Ram_2013,Wei_2016}; see also \cite{Datta_2006} for systems of the third order. For such systems, the problem of spectrum assignment or partial spectrum assignment is solved. In \cite{Yu_2009,Cai_2012}, systems of higher orders are considered. In eigenstructure assignment problems (see, e.g., \cite{Yu_2009}; see also \cite{Duan_2015}  and references therein), in addition to assigning eigenvalues, it is necessary to assign eigenvectors. For example, Lancaster \cite{Lancaster_2005} considers the following problem: Given complete spectral data for a system (i.e., complete information on eigenvalues and eigenvectors), define a corresponding system. The problem of AMCA for CMP allows solving the task of assigning not only eigenvalues but also eigenvectors (see \cite{Yaici_2014}; see also \cite[Theorem 11]{Zaitsev_2021} and \cite[Theorems 8 and 9]{Zaitsev_2024-2}). In other works, block poles placement was used in problems of control for two-wheeled robot 
\cite{Nail_2019}, in applications to gas turbine power plant \cite{Nail_2018}, to gas compressor system \cite{Nail_2018-2}, in tank processes \cite{Bekhiti_2017}, and
in other problems \cite{Bekhiti_2015,Bekhiti_2016,Nail_2019-2}. 
Thus, the problem of AMCA is relevant for applications.
\end{remark}

\begin{remark} \label{rem-dop-2} 
Let $s=1$. Then, the problem of AMCA turns into the following problem. For a Frobenius matrix with arbitrary given coefficients, it is required to construct a control that ensures similarity between the matrix of the closed-loop system and the given Frobenius matrix. Note that this problem is stronger than simply the problem of assigning coefficients of a characteristic polynomial: from the similarity of a certain matrix to a Frobenius matrix with given coefficients, it follows that the characteristic polynomials coincide, but the converse is not true. Nevertheless, for the solvability of both of these problems, the condition of complete controllability of an open-loop system is necessary and sufficient.

These tasks are classical problems of mathematical control theory. Here, we consider tasks that are  extensions of these problems to the case when $s>1$, so, these problems have important theoretical significance. Let us describe the previously obtained results. Let $s\ge 1$ and let $m=1$.
In \cite{Leang_Shieh_1982,Leang_Shieh_1983} the following statement was proved: (St1) {\it Condition 
\begin{equation}\label{072}
{\rm rank}\,[{G},{F}{G},\ldots,{F}^{n-1}{G}]=ns
\end{equation}
is necessary and sufficient for the following property to be satisfied: there exists a transformation $z=Sx$ that reduces the system  \eqref{066}, \eqref{067} to the 
 system \eqref{026} with matrices \eqref{062} with $B_1=I$.} It is clear that then, by choosing the appropriate gain coefficients, one can assign arbitrary block coefficients in the last block row of the matrix  $\mathcal{F}+\mathcal{G}\mathcal{L}$. From this, it follows that the following statement is true: (St2) {\it Let $m=1$. Condition \eqref{072} is sufficient for the solvability of the problem of AMCA for CMP by LSSF.}
(Note that, in statement (St2), the necessity of condition \eqref{072} does not follow from the above and proving the necessity requires additional research).

In addition, in \cite{Leang_Shieh_1983}, the block matrix spectrum assignment problem was solved for system \eqref{066}, \eqref{067} for the case $m=1$, in a certain sense (see Sect.~\ref{sect-assignment}). Results of \cite{Leang_Shieh_1983} was used in problems 
of block partial fraction expansion
\cite{Leang_Shieh_1986,Leyva_Ramos_1991},
cascade decomposition \cite{Leang_Shieh_1985},
eigenvalue assignment \cite{Shafai_1987,Shafai_1988},
eigenstructure assignment \cite{Yaici_2014},
and in some applied problems \cite{Hinkkanen_2015,Nail_2018}.

Now consider the case when $m>1$. For the case $s=1$, the theory is well known and is classical. Let $s>1$. In this case,  the questions of the analogues of  statements (St1) (and (St2)) are open and unproven to this day. We have a hypothesis that the analogues of these statements are true for arbitrary $m>1$ and $s>1$ but these statements have not yet been proven.

Further, in the case when $s>1$, new features appear that are absent in the case when $s=1$. Namely,
in Definition~\ref{def-01}, the matrices $\Gamma_i$ are chosen arbitrarily from the space $M_s(\mathbb{K})$. Let us assume that they are not chosen arbitrarily, but, for example, are upper triangular, or lower triangular, or diagonal. If $s=1$, then these definitions will obviously coincide. If $s>1$, then these definitions are generally different.
What will be the relationship between these definitions?
The study of this question was started in \cite{Zaitsev_2024-2}.
In particular, it was proved (for the case $m=2$ and $s=2$ and $\mathbb{K}=\mathbb{C}$) that if we take upper triangular matrices instead of arbitrary ones in Definition \ref{def-01}, then the resulting definition will be equivalent. For arbitrary $m>2$ and $s>2$ (for $\mathbb{K}=\mathbb{C}$), these statements have not yet been proven.
All of the above questions are important from a theoretical point of view and are not trivial.
\end{remark}

\begin{remark} \label{rem-dop-3} 
Block matrix systems have been studied in other works. In \cite{Pereira-2007}, the block versions of the classical Hermite, Routh, Hurwitz and Schwarz matrices are presented for the free differential equation \eqref{040} (i.e., with $u=0$). In \cite{Vitoria_1982}, a block Cayley--Hamilton theorem is proved. In \cite{Victoria_1982-2}, some results about the block eigenvalues of block compound matrices and additive block compound matrices are obtained. The results on the stability of block matrix systems are collected and presented in \cite{Martins_2014}. 
We note the following difference with our study. In works \cite{Vitoria_1982,Victoria_1982-2,Martins_2014}, when considering block matrix systems, it is assumed that the matrix blocks commute pairwise. This condition is quite restrictive. In the present work, we do not impose such a condition.
\end{remark}

The main aim of the present work is to pose a problem of AMCA by {\it linear static output feedback} and to solve it (in some partial cases). Let $s\in \mathbb{N}$. Consider an input-output linear control system with block matrix coefficients:
\begin{gather}\label{080}
\dot x= Fx+Gu, \quad x\in \mathbb{K}^{ns}, \quad u \in \mathbb{K}^{ms}, \\
\label{082}
y=Hx, \quad y\in \mathbb{K}^{ks}, \\
\label{084}
{F}=\begin{bmatrix}
{F}_{11} &  \ldots & {F}_{1n} \\
\vdots &  & \vdots  \\
{F}_{n1} &  \ldots & {F}_{nn} 
\end{bmatrix}, \quad
{G}=\begin{bmatrix}
{G}_{11} &  \ldots & {G}_{1m} \\
\vdots &  & \vdots  \\
{G}_{n1} &  \ldots & {G}_{nm} 
\end{bmatrix}, \quad
{H}=\begin{bmatrix}
{H}_{11} &  \ldots & {H}_{1n} \\
\vdots &  & \vdots  \\
{H}_{k1} &  \ldots & {H}_{kn} 
\end{bmatrix}.
\end{gather}
Here $x\in \mathbb{K}^{ns}$ is a state vector, $u \in \mathbb{K}^{ms}$ is a control vector, 
$y\in \mathbb{K}^{ks}$ is an output vector; $F_{ij}, G_{j\alpha}, H_{\beta i} \in M_{s}(\mathbb{K})$, $i,j=\overline{1,n}$, $\alpha=\overline{1,m}$, $\beta=\overline{1,k}$.
Suppose that the control in system \eqref{080}, \eqref{082}, \eqref{084} has the form of linear static output feedback (LSOF):
\begin{equation}\label{086}
u=Qy.
\end{equation}
Here $Q=\{Q_{\alpha\beta}\}\in M_{ms,ks}(\mathbb{K})$, $Q_{\alpha\beta}\in M_s(\mathbb{K})$,
$\alpha=\overline{1,m}$, $\beta=\overline{1,k}$. The closed-loop system has the form
\begin{equation}\label{088}
\dot x= (F+GQH)x, \quad x\in \mathbb{K}^{ns}.
\end{equation}
\begin{definition}\label{def-02} 
We say that, for system \eqref{080}, \eqref{082}, \eqref{084}, the problem of {\it arbitrary matrix coefficient assignment (AMCA) for the characteristic matrix polynomial (CMP) by linear static output feedback (LSOF) is resolvable} if for any $\Gamma_i\in M_s(\mathbb{K})$, $i=\overline{1,n}$, there exists a gain matrix $Q\in M_{ms,ks}(\mathbb{K})$  such that the closed-loop system \eqref{088} is reducible by some change of variables $z=Sx$ to the system \eqref{070} with the matrix \eqref{064}, that is the matrix $F+GQH$ of the system \eqref{088} is similar to the matrix \eqref{064}.
\end{definition}

If $s=1$, then the problem of AMCA for CMP by LSOF formulated in Definition \ref{def-02}
coincides with the classical problem of eigenvalue assignment (or assignment of coefficients of the characteristic polynomial) by linear static output feedback. So, we are considering a generalization of the classical problem of eigenvalue assignment by linear static output feedback to systems with block matrix coefficients. 

The static output feedback problem of eigenvalue assignment is one of the most important open questions in control theory \cite{Brockett_1999,Rosenthal_1999}, see also reviews \cite{Syrmos_1997,Sadabadi_2016,Shumafov_2019}. This problem has been studied by many authors for a long time, including recent time (see some references in \cite{Zaitsev_2021}). However, this problem does not yet have a complete constructive solution in the general case. 

This problem has been solved in \cite{Zaitsev_2003} in the following partial case.
Consider the following control system:
\begin{equation}
\label{210}
\begin{gathered}
x^{(n)}+A_{1}x^{(n-1)}+
\ldots  {}+ A_{n}x  = \\
\begin{aligned}
&=B_{p1}u_1^{(n-p)}+B_{p+1,1}u_1^{(n-p-1)}+\ldots+B_{n1}u_1+\ldots \\
&+B_{pm}u_m^{(n-p)}+B_{p+1,m}u_m^{(n-p-1)}+
\ldots+B_{nm}u_m,
\end{aligned}
\end{gathered}
\end{equation}
\begin{equation}
\label{220}
\begin{aligned}
y_1&={C}_{11}x+{C}_{21}x'+\ldots+{C}_{p1}x^{(p-1)}, \quad \ldots, \\
y_k&={C}_{1k}x+{C}_{2k}x'+\ldots+{C}_{pk}x^{(p-1)},
\end{aligned}
\end{equation}
with linear static output feedback
\begin{gather}\label{230}
u=Qy.
\end{gather}
Here $x\in\mathbb{K}^s$ is a state variable, $u_{\alpha}\in\mathbb{K}^s$  are control variables, $y_{\beta}\in\mathbb{K}^s$ are
output variables,
$A_i$, $B_{l\alpha}$, $C_{\nu\beta}\in M_{s}(\mathbb{K})$,
$i=\overline{1,n}$, $l=\overline{p,n}$, $\nu=\overline{1,p}$, $\alpha=\overline{1,m}$,
$\beta=\overline{1,k}$; $p\in\{1,\ldots,n\}$ is a pregiven number; $u={\rm
col}(u_1,\ldots,u_m)\in\mathbb{K}^{ms}$, $y={\rm col}(y_1,\ldots,y_k)\in\mathbb{K}^{ks}$. 
Let us construct the matrices  $B=\{B_{l\alpha}\}$, $l=\overline{1,n}$, $\alpha=\overline{1,m}$, and $C=\{c_{\nu \beta}\}$, $\nu=\overline{1,n}$, $\beta=\overline{1,k}$, where $B_{l\alpha}\!:=0\in M_{s}(\mathbb{K})$ for
$l<p$ and $C_{\nu \beta}\!:=0\in M_{s}(\mathbb{K})$  for $\nu>p$. Let $J\!:=\{\vartheta_{ij}\}\in M_n(\mathbb{R})$ where $\vartheta_{ij}=1$ for $j=i+1$ and $\vartheta_{ij}=0$ for $j\ne i+1$. Let $T$ denote the transposition of a matrix. The following theorem has been proved in \cite{Zaitsev_2003} for the case $s=1$ (that is for the case when matrix coefficients are scalar coefficients).
\begin{theorem}\label{theo-0-1}
\it
Let $s=1$. System \eqref{210}, \eqref{220} is arbitrary scalar coefficient assignable by LSOF \eqref{230} iff the matrices
\begin{equation*}
C^T B,\quad C^T JB, \quad  \ldots, \quad  C^T J^{n-1}B
\end{equation*}
are linearly independent.
\end{theorem}
Next, in \cite{Zaitsev_2009,Zaitsev_2010} this result was generalized to systems \eqref{080}, \eqref{082}, \eqref{084} for the case if the coefficients of the system have the following special form: 
\begin{equation}\label{232}
\begin{gathered}
s=1; \quad F_{i,i+1}\ne 0, \; i=\overline{1,n-1}; \quad F_{ij}=0,\; j>i+1; \\
G_{j\alpha}=0, \; j=\overline{1,p-1}, \; \alpha=\overline{1,m}; \quad
H_{\beta i}=0, \; \beta=\overline{1,k}, \; i=\overline{p+1,n}.
\end{gathered}
\end{equation}
\begin{theorem}[{see \cite[Theorem 1]{Zaitsev_2010}}]\label{theo-0-2}
\it
Suppose that conditions \eqref{232} hold. System \eqref{080}, \eqref{082}, \eqref{084}  is arbitrary scalar coefficient assignable by LSOF \eqref{086} iff the matrices
\begin{equation*}
H G,\quad H FG, \quad  \ldots, \quad  H F^{n-1}G
\end{equation*}
are linearly independent.
\end{theorem}

In the paper \cite{Zaitsev_2021}, for system \eqref{210}, \eqref{220}, \eqref{230}, a generalization of Theorem \ref{theo-0-1} has been obtained from the case $s=1$ to the case of arbitrary $s\ge 1$. The present work is a continuation of \cite{Zaitsev_2021}. Here we study the problem of AMCA for CMP by LSOF, for system \eqref{080}, \eqref{082}, \eqref{084}. 
We generalize Theorem \ref{theo-0-2} from the case $s=1$ to the case of arbitrary $s\ge 1$.
The results of the present work generalize the results of \cite{Zaitsev_2021}.

\section{Notations, definitions, and auxiliary statements\label{sect-definitions}}

We will use some notations, definitions, and statements from \cite{Zaitsev_2021}.
Here and throughout, we suppose that the numbers $s,n,m,k\in\mathbb{N}$, and $p\in\{1,\ldots,n\}$ are fixed. For any matrix $A\in M_{\omega}$, we suppose, by definition, $A^0=I\in M_{\omega}$, where $I$ is the identity matrix; $[e_1,\ldots,e_{\omega}]\!:=I\in M_{\omega}$.
Denote by $\otimes$ the right Kronecker product of matrices  $A=\{a_{ij}\}\in M_{\omega,\rho}$, $i=\overline{1,\omega}$, $j=\overline{1,\rho}$, and $B\in M_{\sigma,\tau}$ \cite[Ch. 12]{Lancaster-Tismenetsky} defined by the formula $$A\otimes B:=\left[\begin{matrix}
a_{11}B & a_{12}B & \dots & a_{1\rho}B \\
a_{21}B & a_{22}B & \dots & a_{2\rho}B \\
\vdots & \vdots & & \vdots \\
a_{\omega 1}B & a_{\omega 2}B & \dots & a_{\omega \rho}B
\end{matrix}
\right]
\in M_{\omega\sigma,\rho\tau}.
$$
Denote $\mathcal{J}\!:=J\otimes I\in M_{ns}$ where $I\in M_s$ and
$J:=\{\vartheta_{ij}\}\in M_n$, $\vartheta_{ij}=1$ for $j=i+1$ and $\vartheta_{ij}=0$ for $j\ne i+1$. 
We will use the mappings ${\rm vecc},\,{\rm vecr}$ that unroll a matrix $A=\{a_{ij}\}\in M_{\omega,\rho}(\mathbb{K})$, $i=\overline{1,\omega}$, $j=\overline{1,\rho}$, column-by-column and row-by-row respectively into the column vector and the row vector respectively:
\begin{gather*}
{\rm vecc}\,A={\rm col}(a_{11},\ldots,a_{\omega 1},\ldots,a_{1\rho},\ldots,a_{\omega\rho})\in M_{\omega\rho,1}(\mathbb{K}), \\
{\rm vecr}\,A=[a_{11},\ldots,a_{1\rho},\ldots,a_{\omega 1},\ldots,a_{\omega\rho}]\in M_{1,\omega\rho}(\mathbb{K}).
\end{gather*}

\if0

\begin{lemma}[{see \cite[Lemma 1]{Zaitsev_2021}}]\label{lem-1}
If $X\in M_{\omega,\rho}$, $Y\in M_{\rho,\sigma}$, $Z\in M_{\sigma,\tau}$, then
\begin{equation*}
{\rm vecc}\,(XYZ)=(Z^T\otimes X)\,{\rm vecc}\, Y.
\end{equation*}
\end{lemma}

\fi

\begin{definition}[{see \cite[Definition 4]{Zaitsev_2021}}]\label{def-SP}
For the fixed $s\in\mathbb{N}$, let us introduce the operation of the block trace ${\rm SP}_s \colon
M_{qs}\rightarrow M_{s}$ by the following rule: if $A=\{A_{ij}\}\in M_{qs}$, $A_{ij}\in M_s$, $i,j=\overline{1,q}$, then ${\rm SP}_s A=\sum\limits_{i=1}^q A_{ii}$.
\end{definition}

\if0
\begin{lemma}[{see \cite[Lemma 2]{Zaitsev_2021}}] \label{lem-Sp-scalar}
Suppose that $X$ and $Y$ are block matrices with $s\times s$-blocks, there exist $XY$ and $YX$, and the blocks of the matrix $Y$ are scalar matrices, i.e.,
\begin{gather*}
X=\{X_{ij}\}\in M_{qs,rs}, \quad X_{ij}\in M_s, \quad i=\overline{1,q}, \quad j=\overline{1,r}; \\
Y=\{Y_{ji}\}\in M_{rs,qs}, \quad Y_{j i}=y_{ji}I, \quad y_{ji}\in\mathbb{K}, \quad I\in M_s, \quad  j=\overline{1,r}, \quad i=\overline{1,q}.
\end{gather*}
Then ${\rm SP}_s(XY)={\rm SP}_s(YX)$.
\end{lemma}

\begin{lemma}[{see \cite[Lemma 3]{Zaitsev_2021}}] \label{lem-SP-J}
Let $D=\{D_{ \omega \rho}\}\in M_{ns}$, $D_{ \omega \rho}\in M_s$, $\omega, \rho=\overline{1,n}$. Then
\begin{equation*} 
{\rm
SP}_s(\mathcal{J}^{i-1}D)=\sum\limits_{\eta=1}^{n-i+1}D_{\eta+i-1,\eta},\quad
i=\overline{1,n}.
\end{equation*}
\end{lemma}

\fi

In what follows, we will need the properties of ${\rm SP}_s$ from Lemmas 2 and 3 of \cite{Zaitsev_2021}.

\begin{definition}[{see \cite[Definition 5]{Zaitsev_2021}}] \label{def-star}
Suppose that $X$ and $Y$ are block matrices with $s\times s$-blocks such that
the number of the (block) columns of $X$ is equal to the number of the (block) rows of $Y$:
\begin{gather*}
X=\{X_{ij}\}\in M_{qs,rs}, \quad X_{ij}\in M_s, \quad i=\overline{1,q}, \quad j=\overline{1,r}; \\
Y=\{Y_{j\nu}\}\in M_{rs,ts}, \quad Y_{j \nu}\in M_s, \quad  j=\overline{1,r}, \quad \nu=\overline{1,t}.
\end{gather*}
For the matrices $X$ and $Y$, let us introduce the operation of the block multiplication by the following rule:
\begin{equation*}
Z=X \star Y:=\{Z_{i\nu}\}, \quad
Z_{i\nu}:=\sum_{j=1}^r X_{ij}\otimes Y_{j \nu},
\quad i=\overline{1,q}, \quad \nu=\overline{1,t}.
\end{equation*}
We have $Z_{i\nu}\in M_{s^2}$ for all  $i=\overline{1,q}$, $\nu=\overline{1,t}$, therefore, $Z:=X \star Y \in M_{qs^2, ts^2}$.
\end{definition}
For convenience, so as not to write brackets, we assume 
\begin{equation*}
P \star RS:=P\star(R\cdot S), \qquad PR\star S:=(P\cdot R) \star S
\end{equation*}
 where  matrices $P,R,S$ have the corresponding dimensions.

\begin{lemma}\label{lem-F-star}
\it
Let
\begin{gather*}
X=\{X_{i\rho}\}\in M_{qs,ns}, \quad X_{i\rho}\in M_s, \quad i=\overline{1,q}, \quad \rho=\overline{1,n}; \\
R=Y\otimes I\in M_{ns}, \quad Y=\{y_{\rho\omega}\}\in M_n, \quad y_{\rho\omega}\in \mathbb{K}, \quad \rho,\omega=\overline{1,n}, \quad I\in M_s, \\
Z=\{Z_{\omega\nu}\}\in M_{ns,rs}, \quad Z_{\omega \nu}\in M_s, \quad  \omega=\overline{1,n}, \quad \nu=\overline{1,r}.
\end{gather*}
Then, 
\begin{equation*}%\label{270}
X R \star Z= X\star RZ.
\end{equation*}
\end{lemma}

The proof of Lemma \ref{lem-F-star} is carried out by direct calculation of blocks of matrices $X R \star Z$ and $X\star RZ$.

\if0

\emph{Proof.}
Denote
\begin{gather*}
X R=:V=\{V_{i\omega}\}\in M_{qs,ns}, \quad 
V_{i\omega}\in M_s,  \quad i=\overline{1,q}, \quad \omega=\overline{1,n};\\
RZ=:W=\{W_{\rho\nu}\}\in M_{ns,rs}, \quad 
W_{\rho\nu}\in M_s,   \quad \rho=\overline{1,n}, \quad \nu=\overline{1,r}.
\end{gather*}
Then, 
\begin{gather*}
V_{i\omega}=\sum_{\rho=1}^n X_{i\rho} \cdot(y_{\rho\omega} I),\quad
W_{\rho\nu}=\sum_{\omega=1}^n (y_{\rho\omega} I)\cdot Z_{\omega \nu}.
\end{gather*}
Thus, by Definition \ref{def-star}, we have
\begin{gather}\label{290}
(XR \star Z)_{i\nu}=(V \star Z)_{i\nu}=\sum_{\omega=1}^{n}V_{i\omega}\otimes Z_{\omega \nu}
=\sum_{\omega=1}^n \sum_{\rho=1}^{n}(X_{i\rho} \cdot y_{\rho\omega}) \otimes Z_{\omega \nu}
= \sum_{\rho=1}^{n} \sum_{\omega=1}^n y_{\rho\omega} \cdot (X_{i\rho} \otimes Z_{\omega \nu}),\\
\label{300}
(X\star RZ)_{i\nu}=(X\star  W)_{ i \nu }=\sum_{\rho=1}^{n} X_{i\rho} \otimes W_{\rho\nu} =\sum_{\rho=1}^{n} X_{i\rho} \otimes (\sum_{\omega=1}^{n} y_{\rho\omega} Z_{\omega \nu})
= \sum_{\rho=1}^{n} \sum_{\omega=1}^n y_{\rho\omega} \cdot (X_{i\rho} \otimes Z_{\omega \nu}).
\end{gather}
By \eqref{290} and \eqref{300}, equality \eqref{270} is true. Q.E.D.
\hfill $\square$
\fi

\begin{remark}
\textrm{Lemma \ref{lem-F-star} has been proven in \cite[Lemma 4]{Zaitsev_2021} for the case  $R=J\otimes I$.}
\end{remark}

\begin{definition}[{see \cite[Definition 6]{Zaitsev_2021}}]  \label{def-block-transpose}
For the fixed $s\in\mathbb{N}$, let us introduce the operation of the block transposition  $\mathcal{T}$ by the following rule: if  $A=\{A_{ij}\}\in M_{qs,rs}$, $A_{ij}\in M_s$, $i=\overline{1,q}$, $j=\overline{1,r}$, then
\begin{equation*}
A^{\mathcal{T}}:=B=\{B_{ji}\}\in M_{rs,qs},\quad B_{ji}:=A_{ij}, \quad  j=\overline{1,r}, \quad i=\overline{1,q}.
\end{equation*}
\end{definition}

\if0
\begin{lemma}[{see \cite[Lemma 5]{Zaitsev_2021}}]\label{lem-transpose-property}
The following properties hold.

1. $(A^{\mathcal{T}})^{\mathcal{T}}=A$.

2. If $X$ and $Y$ are block matrices with $s\times s$-blocks, there exists $XY$, and the blocks of the matrix $Y$ are scalar matrices, i.e.,
\begin{gather*}
X=\{X_{ij}\}\in M_{qs,rs}, \quad X_{ij}\in M_s, \quad i=\overline{1,q}, \quad j=\overline{1,r}; \\
Y=\{Y_{j\sigma}\}\in M_{rs,ts}, \quad Y_{j \sigma}=y_{j\sigma}I, \quad y_{j\sigma}\in\mathbb{K}, \quad I\in M_s, \quad  j=\overline{1,r}, \quad \sigma=\overline{1,t},
\end{gather*}
then 
\begin{equation}\label{336-1}
(XY)^{\mathcal{T}}=Y^{\mathcal{T}} X^{\mathcal{T}}.
\end{equation}
\end{lemma} 

\fi

In what follows, we will need the properties of $\mathcal{T}$  from \cite[Lemma 5]{Zaitsev_2021}.

\begin{definition} [{see \cite[Definition 7]{Zaitsev_2021}}] \label{def-VECCC}
Let $X$ be a block matrix with $s\times s$-blocks:
\begin{gather*}
X=\{X_{ij}\}\in M_{qs,rs}, \quad X_{ij}\in M_s, \quad i=\overline{1,q}, \quad j=\overline{1,r}.
\end{gather*}
Let us construct the mappings ${\rm VecCR}_s,\,{\rm VecRR}_s:M_{qs,rs}\to M_{s,qrs}$ that unroll the matrix $X=\{X_{ij}\}\in M_{qs,rs}$ by block columns and by block rows  respectively into the block row with $s\times s$-blocks:
\begin{gather*}
{\rm VecCR}_s\,X=[X_{11},\ldots,X_{q1},\ldots,X_{1r},\ldots,X_{qr}], \\
{\rm VecRR}_s\,X=[X_{11},\ldots,X_{1r},\ldots,X_{q1},\ldots,X_{qr}],
\end{gather*}
and the mappings ${\rm VecRC}_s,\,{\rm VecCC}_s:M_{qs,rs}\to M_{qrs,s}$ that unroll the matrix $X=\{X_{ij}\}\in M_{qs,rs}$ by block rows and by block columns respectively into the block column with $s\times s$-blocks:
\begin{gather*}
{\rm VecRC}_s\,X=\left[\begin{matrix}X_{11}\\ \vdots \\ X_{1r} \\ \vdots \\ X_{q1} \\ \vdots \\ X_{qr} \end{matrix}\right], \quad
{\rm VecCC}_s\,X=\left[\begin{matrix}X_{11}\\ \vdots \\  X_{q1} \\ \vdots \\ X_{1r} \\ \vdots \\ X_{qr} \end{matrix}\right].
\end{gather*}
\end{definition}

The following equalities are clear:
\begin{gather}\label{336-3}
{\rm VecCR}_s\,X={\rm VecRR}_s\,(X^{\mathcal{T}}),\\
\label{336-5}
({\rm VecRC}_s\,X)^T={\rm VecCR}_s\,(X^{T}).
\end{gather}

\if0

\begin{lemma} [{see \cite[Lemma 6]{Zaitsev_2021}}] \label{lem-SP-vec}
If
\begin{gather*}
X=\{X_{ij}\}\in M_{qs,rs}, \quad X_{ij}\in M_s, \quad i=\overline{1,q}, \quad j=\overline{1,r},\\
Y=\{Y_{ji}\}\in M_{rs,qs}, \quad Y_{ji}\in M_s, \quad j=\overline{1,r}, \quad i=\overline{1,q},
\end{gather*}
then
\begin{align}\label{337}
{\rm SP}_s (XY)& ={\rm VecRR}_s X \cdot {\rm VecCC}_s Y = \\
\label{338}
& = {\rm VecCR}_s X \cdot {\rm VecRC}_s Y.
\end{align}
\end{lemma}

\fi 

In what follows, we will need the assertions of \cite[Lemma 6]{Zaitsev_2021}.

\begin{lemma} [{see \cite[Lemma 7]{Zaitsev_2021}}]\label{lem-unroll}
\it
Let
\begin{gather*}
X=\{X_{l\alpha}\}\in M_{ns,ms}, \quad X_{l\alpha}\in M_s, \quad l=\overline{1,n}, \quad \alpha=\overline{1,m}; \\
Q=\{Q_{\alpha \beta}\}\in M_{ms,ks}, \quad Q_{\alpha \beta}\in M_s, \quad  \alpha=\overline{1,m}, \quad \beta=\overline{1,k};\\
Y=\{Y_{\beta l}\}\in M_{ks,ns}, \quad Y_{\beta l}\in M_s, \quad \beta=\overline{1,k}, \quad l=\overline{1,n}.
\end{gather*}
Suppose that $R={\rm SP}_s (XQY)$.
Then,
\begin{equation*}%\label{340}
{\rm vecc}\,R = {\rm VecRR}_{s^2}\big((Y^{\mathcal{T}})^T\star X) \cdot{\rm vecc}\,({\rm VecCR}_{s}\,Q).
\end{equation*}
\end{lemma}

Consider a block matrix 
\begin{equation}\label{350}
Z=\{Z_{ij}\}\in M_{ns}, \quad Z_{ij}\in M_s, \quad i,j=\overline{1,n}.
\end{equation}
It is said that the matrix $Z$ is a {\it lower block Hessenberg matrix} if $Z_{ij}=0\in M_s$, $j>i+1$. If, in addition, $\det Z_{i,i+1}\ne 0$, then this lower block Hessenberg matrix is called {\it unreduced}.
We will consider only those lower block Hessenberg matrices that are unreduced, therefore, for brevity, we will omit the word ``unreduced''.  If, for the block matrix \eqref{350}, we have $Z_{i,i+1}=I\in M_s$, $i=\overline{1,n-1}$, and $Z_{ij}=0\in M_s$, $i=\overline{1,n-1}$, $j=\overline{1,n}$, $i+1\ne j$ (that is $Z$ has the form \eqref{064}), then  it is said that $Z$  is a {\it lower block Frobenius matrix}.

\begin{lemma} \label{lemma-00}
\it
Suppose that a block matrix \eqref{350} is a lower block Hessenberg matrix. Then, there exists a non-degenerate lower block triangular matrix $S$ such that the matrix $SZS^{-1}$ is a lower block Frobenius matrix.
\end{lemma}

\emph{Proof.}
Let us remove the last block row in the matrix $Z$ and denote the resulting matrix by $Q\in M_{(n-1)s,ns}$. Take $e_1\in\mathbb{R}^n$, $I\in M_s$. Then, $e_1^T\otimes I \in M_{s,ns}$. Construct $S_1\!:=\left[ \begin{matrix} e_1^T\otimes I \\ Q \end{matrix} \right]$. Then, the matrix $S_1$ has the form
\begin{equation}\label{355}
S_1=\left[\begin{matrix}
I & 0 & \ldots & 0 \\
Z_{11} & Z_{12}  & \ldots & 0 \\
\vdots & \vdots & \ddots & \vdots  \\
Z_{n-1,1} & Z_{n-1,2}  & \ldots & Z_{n-1,n}
\end{matrix} \right], \qquad 0,I\in M_s.
\end{equation}
Hence, $S_1\in M_{ns}$ is a lower  block triangular matrix and $\det S_1\ne 0$. For each $l=2,\ldots,n-1$, on the basis of the matrix $S_{l-1}=\{S_{ij}^{(l-1)}\}\in M_{ns}$, $S_{ij}^{(l-1)}\in M_{s}$, $i,j=\overline{1,n}$, we construct the matrix $S_{l}=\{S_{ij}^{(l)}\}\in M_{ns}$, $S_{ij}^{(l)}\in M_{s}$, $i,j=\overline{1,n}$, as follows: $S_{11}^{(l)}\!:=I\in M_s$, $S_{1j}^{(l)}\!:=S_{j1}^{(l)}\!:=0\in M_s$, $j=\overline{2,n}$; $S_{ij}^{(l)}\!:=S_{i-1,j-1}^{(l-1)}$, $i,j=\overline{2,n}$. Then, the matrices $S_l$ are lower block triangular non-degenerate matrices for all $l=1,\ldots,n-1$. Let $S=S_{n-1}\cdot\ldots \cdot S_1$. Then, $S$ is a lower block triangular non-degenerate matrix as well. Let us show that the matrix $S$ is required.

Due to \eqref{355}, we have: $QS_1^{-1}=[0\ \ I]\in M_{(n-1)s,ns}$, $0\in M_{(n-1)s,s}$, $I\in M_{(n-1)s}$. Therefore, 
\begin{equation*}%\label{358}
ZS_1^{-1}=\left[\begin{matrix} Q\\ *\end{matrix} \right]\cdot S_1^{-1} = \left[\begin{matrix} 0 & I\\
* & * \end{matrix} \right], \qquad
0\in M_{(n-1)s,s}, \quad I\in M_{(n-1)s}.
\end{equation*}
Next, we have, $S_1=\left[\begin{matrix} Q' & 0\\
* & * \end{matrix} \right]$, where $0\in M_{(n-1)s,s}$ and $Q'\in M_{(n-1)s}$ is lower block triangular. Hence,
\begin{equation}\label{360}
S_1ZS_1^{-1}=\left[\begin{matrix} Q' & 0\\
* & * \end{matrix} \right] \cdot  \left[\begin{matrix} 0 & I\\
* & * \end{matrix} \right]= \left[\begin{matrix} 0 & Q'\\
* & * \end{matrix} \right], \qquad 0\in M_{(n-1)s,s}.
\end{equation}
Denote $X_1\!:=S_1ZS_1^{-1}$. Hence, from \eqref{360}, we have
\begin{equation*}%\label{362}
X_1=\left[\begin{matrix}
0 & I & 0 & \ldots & 0 \\
0 & Z_{11} & Z_{12}  & \ldots & 0 \\
\vdots & \vdots & \vdots & \ddots  & \vdots \\
0 & Z_{n-2,1} & Z_{n-1,2}  & \ldots & Z_{n-2,n-1} \\
* & * & * & \ldots & * 
\end{matrix} \right].
\end{equation*}

Let us remove the last block row in the matrix $X_1$ and add from above the block row $e_1^T\otimes I = [ I \ \ 0 \ \ \dots \ \ 0 ] \in M_{s,ns}$. Then, we obtain the matrix $S_2$. Construct $X_2\!:=S_2X_1S_2^{-1}$. Then, $X_2=S_2S_1ZS_1^{-1}S_2^{-1}$. By using the similar arguments as above, we obtain that
\begin{equation*}%\label{364}
X_2=\left[\begin{matrix}
0 & I & 0 & 0 & \ldots & 0 \\
0 & 0 & I & 0 &\ldots & 0 \\
0 & 0 & Z_{11} & Z_{12}  & \ldots & 0 \\
\vdots & \vdots & \vdots & \vdots & \ddots  & \vdots \\
0 & 0 & Z_{n-3,1} & Z_{n-3,2}  & \ldots & Z_{n-3,n-2} \\
* & * & * & * &\ldots & * 
\end{matrix} \right].
\end{equation*}
Applying this operation $n-1$ times, we get the matrix $X_{n-1}=S_{n-1}\cdot\ldots \cdot S_1 Z S_1^{-1} \cdot \ldots \cdot S_{n-1}^{-1} = SZS^{-1}$. By construction, we obtain that $X_{n-1}$ is a lower block Frobenius matrix. Q.E.D.
\hfill $\square$

\section{Sufficient conditions to solving the problem of AMCA for CMP by LSOF for systems with a lower block Frobenius matrix \label{sect-main}}

Consider system \eqref{080}, \eqref{082}, \eqref{084}. Suppose that the coefficients of this system have the following special form: for some $p\in\{1,\ldots,n\}$, the first $p-1$ block rows of the matrix $G$ are zero, the last $n-p$ block columns of the matrix $H$ are zero, the matrix $F$ is a lower block Frobenius matrix, i.e.,
\begin{gather}\label{400}
F=\left[\begin{matrix}
0 & I & 0 & \ldots & 0 \\
0 & 0 & I & \ldots & 0 \\
\vdots & \vdots & \vdots & \ddots & \vdots  \\
0 & 0 & 0 & \ldots & I \\
-A_n & -A_{n-1} & -A_{n-2} & \ldots & -A_1
\end{matrix} \right], \qquad 0,I,A_i\in M_s, \quad i=\overline{1,n}, \\
\label{410}
G=\begin{bmatrix} 0 & \ldots & 0\\
\vdots & & \vdots\\
0 & \ldots & 0\\
G_{p1}& \ldots & G_{pm}\\
\vdots & & \vdots\\
G_{n1} & \ldots & G_{nm}
\end{bmatrix}, \qquad 0,G_{j\alpha}\in M_s, \quad j=\overline{p,n}, \quad \alpha=\overline{1,m}, \\
\label{420}
H=\begin{bmatrix} H_{11} & \ldots & H_{1p} & 0 & \ldots & 0\\
\vdots & & \vdots & \vdots & & \vdots\\
H_{k1} & \ldots & H_{kp} & 0 & \ldots & 0
\end{bmatrix}, \qquad 0,H_{\beta i}\in M_s, \quad \beta=\overline{1,k}, \quad i=\overline{1,p}.
\end{gather}

Consider the matrices
\begin{equation*}
(H^{\mathcal{T}})^T\star G, \quad (H^{\mathcal{T}})^T\star FG, \quad \ldots, \quad 
(H^{\mathcal{T}})^T\star F^{n-1}G.
\end{equation*}
We have  $(H^{\mathcal{T}})^T\in M_{ks,ns}$, $F^{i-1}G\in M_{ns,ms}$, hence, $(H^{\mathcal{T}})^T\star F^{i-1}G \in M_{ks^2,ms^2}$ for all $i=\overline{1,n}$. Let us construct the matrices ${\rm VecRR}_{s^2}\big((H^{\mathcal{T}})^T\star F^{i-1}G\big) \in M_{s^2,kms^2}$, $i=\overline{1,n}$, and the matrix
\begin{equation}\label{425}
\Theta=\left[\begin{matrix}
{\rm VecRR}_{s^2}\big((H^{\mathcal{T}})^T\star G\big)\\
{\rm VecRR}_{s^2}\big((H^{\mathcal{T}})^T\star FG\big)\\
\ldots\ldots\ldots\ldots\ldots\ldots\ldots \\
{\rm VecRR}_{s^2}\big((H^{\mathcal{T}})^T\star F^{n-1}G\big)
\end{matrix}\right]\in M_{ns^2,kms^2}.
\end{equation}

\begin{theorem}\label{teo-AMCA-LSOF-1}
\it
For system \eqref{080}, \eqref{082},  with coefficients \eqref{400}, \eqref{410}, \eqref{420}, the problem of AMCA for CMP by LSOF is resolvable, if the following condition holds:
\begin{equation}
\label{428}
{\rm rank}\,\Theta = ns^2.
\end{equation}
\end{theorem}

\section{Proof of Theorem \ref{teo-AMCA-LSOF-1} \label{section-proof-1} }

Let the matrix $F$ have the form \eqref{400}. Here and everywhere below we assume that $A_0:=I\in M_s$. From the matrix $F$, construct the following matrix:
\begin{equation}\label{430}
{P}:=\begin{bmatrix}
A_0&0&\ldots&0\\
A_1&A_0&\ldots&0\\
\vdots&\ddots&\ddots&\vdots\\
A_{n-1}&\ldots&A_1&A_0
\end{bmatrix}\in M_{ns}.
\end{equation}

Let a block matrix $D\in M_{ns}$ with $s\times s$-blocks have the following form: the first $p-1$ block rows and the last $n-p$ block columns of the matrix $D$ are zero, i.e.,
\begin{equation}\label{435}
D=\begin{bmatrix} 
0 & \ldots & 0 & 0 & \ldots & 0 \\
\vdots & & \vdots & \vdots & & \vdots \\
0 & \ldots & 0 & 0 & \ldots & 0 \\
D_{p1} & \ldots & D_{pp} & 0 & \ldots & 0\\
\vdots & & \vdots & \vdots & & \vdots\\
D_{n1} & \ldots & D_{np} & 0 & \ldots & 0
\end{bmatrix}, \qquad 0,D_{\tau \sigma}\in M_s, \quad \tau=\overline{p,n}, \quad \sigma=\overline{1,p}.
\end{equation}
Construct the matrix $Z:=F+D$. This matrix is a lower block Hessenberg matrix. By Lemma \ref{lemma-00}, this matrix is reducible to a lower block Frobenius matrix $\Phi$ of \eqref{064} by some lower block triangular matrix $S$. 

Let us prove preliminarily an auxiliary assertion.

\begin{lemma} \label{lemma-02}
\it
Let two block matrices $C_1,C_2\in M_{ns}$ be given: $C_1=\left[\begin{matrix} 0 & 0 \\ C & 0 \\ \Psi & 0 \end{matrix}\right]$ and $C_2=\left[\begin{matrix} 0 & 0 & 0 \\ O_2 & C & 0 \end{matrix}\right]$, where $O_2=0\in M_{(n-p)s,s}$,
$
C=\left[\begin{matrix}
C_{p1} & \ldots & C_{pp} \\
\vdots & & \vdots \\
C_{n-1,1} & \ldots & C_{n-1,p}
\end{matrix}\right]
\in M_{(n-p)s,ps}$,
$\Psi =\left[\begin{matrix}
C_{n1} & \ldots & C_{np}
\end{matrix}\right]
\in M_{s,ps}$, $C_{ij}\in M_s$, $i=\overline{p,n}$, $j=\overline{1,p}$.
Then,
\begin{align}\label{445}
&{\rm SP}_s (\mathcal{J}^{r} PC_1)={\rm SP}_s (\mathcal{J}^{r} P C_2), \quad \text{ for all } \quad r=0,\ldots,n-p-1, \\
\label{447}
&{\rm SP}_s (\mathcal{J}^{r} PC_1)={\rm SP}_s (\mathcal{J}^{r} P C_2)+\sum_{i=p}^n A_{n-i} C_{i,n-r}, \quad \text{ for all } \quad r=n-p,\ldots,n-1.
\end{align}
\end{lemma}

\emph{Proof.}
Let us represent the matrix $C_1$ as the sum $C_1=C_3+C_4$, where $C_3=\left[\begin{matrix} 0 & 0 \\ C & 0 \\ 0 & 0 \end{matrix}\right]$ and $C_4=\left[\begin{matrix} 0 & 0 \\ 0 & 0 \\ \Psi & 0 \end{matrix}\right]$. Then,
\begin{equation}\label{451}
{\rm SP}_s (\mathcal{J}^{r} PC_1)={\rm SP}_s (\mathcal{J}^{r} PC_3)+{\rm SP}_s (\mathcal{J}^{r} PC_4). 
\end{equation}
First, we find ${\rm SP}_s (\mathcal{J}^{r} PC_4)$. By construction, we have: $PC_4=C_4$. Hence, 
\begin{equation}\label{452}
{\rm SP}_s (\mathcal{J}^{r} PC_4)={\rm SP}_s (\mathcal{J}^{r} C_4).
\end{equation}
By \cite[Lemma 3]{Zaitsev_2021}, we have, for $i=\overline{1,n}$,
\begin{equation}\label{454}
{\rm SP}_s (\mathcal{J}^{i-1} C_4)=\left\{
\begin{aligned}
& C_{n,n-i+1}, & & \text{ if } n-i+1\le p, \\
& 0, & &  \text{ if } n-i+1 > p.
\end{aligned} \right.
\end{equation}
Hence, by \eqref{452} and \eqref{454}, we get:
\begin{equation}\label{456}
{\rm SP}_s (\mathcal{J}^{r} P C_4)=\left\{
\begin{aligned}
& C_{n,n-r}, & & \text{ for } r=n-p,\ldots,n-1, \\
& 0, & &  \text{ for } r=0,\ldots,n-p-1.
\end{aligned} \right.
\end{equation}

Now, let us consider $C_2$ and $C_3$. Construct $B_2:=PC_2$, $B_3:=PC_3$. Then, the matrices $B_2$ and $B_3$ have the following form: $B_3=\left[\begin{matrix} 0 & 0 \\ B & 0 \\ \widetilde \Psi & 0 \end{matrix}\right]$, $B_2=\left[\begin{matrix} 0 & 0 & 0 \\ 0 & \widehat{B} & 0 \\ \end{matrix}\right]$, $B,\widehat{B}\in  M_{(n-p)s,ps}$, $\widetilde\Psi \in M_{s,ps}$, and it is easy to see that $B=\widehat{B}$. 

Represent the matrix $B_3$ as the sum $B_3=B_5+B_6$, where $B_5=\left[\begin{matrix} 0 & 0 \\ B & 0 \\ 0 & 0 \end{matrix}\right]$ and $B_6=\left[\begin{matrix} 0 & 0 \\ 0 & 0 \\ \widetilde \Psi & 0 \end{matrix}\right]$. Then,
\begin{equation}\label{460}
{\rm SP}_s (\mathcal{J}^{r} B_3)={\rm SP}_s (\mathcal{J}^{r} B_5)+{\rm SP}_s (\mathcal{J}^{r} B_6). 
\end{equation}
Due to the equality $B=\widehat{B}$, from \cite[Lemma 3]{Zaitsev_2021}, it follows that, for any $r=0,\ldots,n-1$,
\begin{equation}\label{462}
{\rm SP}_s (\mathcal{J}^{r} B_5)={\rm SP}_s (\mathcal{J}^{r} B_2). 
\end{equation}

Next, by construction, we have $\widetilde\Psi= \left[\begin{matrix}
\Psi_1 & \ldots & \Psi_p \end{matrix}\right]$, $\Psi_j\in M_s$, where
\begin{equation}\label{464}
\Psi_j=A_{n-p}C_{pj}+\ldots +A_1 C_{n-1,j} = \sum_{i=p}^{n-1} A_{n-i} C_{i,j}.
\end{equation}
By \cite[Lemma 3]{Zaitsev_2021}, we have, for $i=\overline{1,n}$,
\begin{equation}\label{466}
{\rm SP}_s (\mathcal{J}^{i-1} B_6)=\left\{
\begin{aligned}
& \Psi_{n-i+1}, & & \text{ if } n-i+1\le p, \\
& 0, & &  \text{ if } n-i+1 > p.
\end{aligned} \right.
\end{equation}
Hence, by \eqref{464} and \eqref{466}, we get:
\begin{equation}\label{468}
{\rm SP}_s (\mathcal{J}^{r} B_6)=\left\{
\begin{aligned}
& \Psi_{n-r}, & &  \text{ for } r=\overline{n-p,n-1}, \\
& 0, & &  \text{ for } r=\overline{0,n-p-1},
\end{aligned} \right. 
=
\left\{
\begin{aligned}
& \sum_{i=p}^{n-1} A_{n-i} C_{i,n-r}, & &  \text{ for } r=\overline{n-p,n-1}, \\
& 0, & &  \text{ for } r=\overline{0,n-p-1}.
\end{aligned} \right. 
\end{equation}

Taking into account \eqref{451}, \eqref{456}, \eqref{460}, \eqref{462} and \eqref{468}, we obtain: $(a)$ for $r=0,\ldots,n-p-1$,
\begin{gather*}
{\rm SP}_s (\mathcal{J}^{r} PC_1)={\rm SP}_s (\mathcal{J}^{r} PC_3)+{\rm SP}_s (\mathcal{J}^{r} PC_4)= {\rm SP}_s (\mathcal{J}^{r} PC_3) = {\rm SP}_s (\mathcal{J}^{r} B_3) \\
= {\rm SP}_s (\mathcal{J}^{r} B_5)+{\rm SP}_s (\mathcal{J}^{r} B_6)= {\rm SP}_s (\mathcal{J}^{r} B_2)+{\rm SP}_s (\mathcal{J}^{r} B_6)= {\rm SP}_s (\mathcal{J}^{r} B_2) = {\rm SP}_s (\mathcal{J}^{r} PC_2);
\end{gather*}

$(b)$ for  $r=n-p,\ldots,n-1$,
\begin{gather*}
{\rm SP}_s (\mathcal{J}^{r} PC_1)={\rm SP}_s (\mathcal{J}^{r} PC_3)+{\rm SP}_s (\mathcal{J}^{r} PC_4)= {\rm SP}_s (\mathcal{J}^{r} PC_3) + C_{n,n-r} = {\rm SP}_s (\mathcal{J}^{r} B_3) + C_{n,n-r} \\
= {\rm SP}_s (\mathcal{J}^{r} B_5)+{\rm SP}_s (\mathcal{J}^{r} B_6) + C_{n,n-r} = {\rm SP}_s (\mathcal{J}^{r} B_2)+
\sum_{i=p}^{n-1} A_{n-i} C_{i,n-r} + A_0 C_{n,n-r}\\ = {\rm SP}_s (\mathcal{J}^{r} PC_2) + \sum_{i=p}^{n} A_{n-i} C_{i,n-r}.
\end{gather*}
Thus, equalities \eqref{445} and \eqref{447} are fulfilled. Q.E.D.
\hfill $\square$

\begin{remark}
For the case $s=1$, Lemma \ref{lemma-02} has been proved in \cite[Lemma 5]{Zaitsev_2009}.
\end{remark}

The following lemma is true.

\begin{lemma} \label{lemma-01}
\it
There exists a lower block triangular matrix $S$ such that the following holds:
the matrix $S(F+D)S^{-1}$ is equal to the matrix $\Phi$ of \eqref{064}, the block matrix coefficients of which are related to the coefficients  of the matrices $F$ and $D$ as follows:
\begin{gather}\label{440}
\Gamma_i=A_i-{\rm SP}_s(\mathcal{J}^{i-1} P D), \quad i=\overline{1,n}.
\end{gather}
\end{lemma}

\emph{Proof.}
Consider the index $p$ of the block matrix $D$ in \eqref{435}. Denote the $p$ by the variable $k\in\{1,\ldots,n\}$. Let us make the proof of Lemma \ref{lemma-01} by  induction on $k$ changing from $n$ to $1$. The basis of induction is the following: let $k=n$. Then, $D=\left[\begin{matrix} 0 \\ \Psi \end{matrix} \right]$, $\Psi = \left[\begin{matrix}  D_{n1} & \ldots & D_{nn} \end{matrix} \right] \in M_{s,ns}$, $D_{nj}\in M_s$, $j=\overline{1,n}$. 

The matrix $F+D$ is itself a lower block Frobenius matrix. Let us take $S=I\in M_{ns}$. Then $\Phi = F+D$, and, hence, $\Gamma_i=A_i-D_{n,n-i+1}$, $i=\overline{1,n}$. On the other hand, we have: $PD=D$. From this, by applying \cite[Lemma 3]{Zaitsev_2021}, we obtain ${\rm SP}_s (\mathcal{J}^{i-1} PD)={\rm SP}_s (\mathcal{J}^{i-1} D)=D_{n,n-i+1}$, for all $i=1,\ldots,n$. Thus, equalities \eqref{440} are fulfilled. The proof of the basis is complete.

Let us make the inductive assumption: let the assertion of the lemma hold for any $k\in\{p+1,\ldots,n\}$. Let us show that it holds for $k=p$ as well. Consider the matrix $D$ in \eqref{435}. By the inductive assumption, $p < n$. Denote
$L\!:=\left[\begin{matrix}
D_{p1} & \ldots & D_{pp} \\
\vdots & & \vdots \\
D_{n-1,1} & \ldots & D_{n-1,p}
\end{matrix}\right]
\in M_{(n-p)s,ps}$,
$\Sigma\!:=\left[\begin{matrix}
D_{n1} & \ldots & D_{np}
\end{matrix}\right]
\in M_{s,ps}$, $\Xi\!:=\left[\begin{matrix}
D_{n1} & \ldots & D_{np} & 0 & \ldots & 0
\end{matrix}\right]
\in M_{s,ns}$, $0,D_{ij}\in M_s$, $i=\overline{p,n}$, $j=\overline{1,p}$. Then,  $D=\left[\begin{matrix} 0 & 0 \\ L & 0 \\ \Sigma & 0 \end{matrix}\right]$.
Let us represent the matrix $D$ in the form $D=D_1+D_2$, where $D_1=\left[\begin{matrix} 0 & 0 \\ L & 0 \\ 0 & 0 \end{matrix}\right]$ and $D_2=\left[\begin{matrix} 0 & 0 \\ 0 & 0 \\ \Sigma & 0 \end{matrix}\right]$. Set $D_3:=\mathcal{J}^T D$. Then, $D_3=\left[\begin{matrix} 0 & 0 \\ L & 0 \end{matrix}\right]\in M_{ns}$. We have $\mathcal{J} D_3=D_1$. Set $D_4:=D_3 \mathcal{J}$. Then, $D_4=\left[\begin{matrix} 0 & 0 & 0 \\ O_4 & L & 0 \end{matrix}\right]\in M_{ns}$, $O_4=0\in M_{(n-p)s,s}$. Set $K:=\left[\begin{matrix} 0 \\ \Omega \end{matrix} \right]\in M_{ns}$, $\Omega = \left[\begin{matrix}  -A_n & \ldots & -A_1 \end{matrix} \right] \in M_{s,ns}$. Then, $F= \mathcal{J}  + K$. Set $T:=I+D_3\in M_{ns}$. Then, $T=\left[\begin{matrix} L_1 & 0 \\ L & L_2  \end{matrix}\right]$, $L_1=I\in M_{ps}$, $L_2=I\in M_{(n-p)s}$. Hence, $T^{-1}=I-D_3\in M_{ns}$. We have
\begin{equation*}\label{470}
(F+D)T^{-1}=(F+D)(I-D_3)=F+D-FD_3-DD_3.
\end{equation*}
The first $p$ block rows of $D_3$ are zero, and the last $n-p$ block columns of $D$ are zero; hence, $DD_3=0\in M_{ns}$. Further,
\begin{equation*}\label{472}
FD_3=(\mathcal{J}+K)D_3=\mathcal{J}D_3+KD_3=D_1+KD_3.
\end{equation*}
Therefore, $D-FD_3=D-D_1-KD_3=D_2-KD_3$.  Thus, $(F+D)T^{-1}=F+D_2-KD_3=F+K_1$ where $K_1=\left[\begin{matrix} 0 \\ \Delta \end{matrix} \right]\in M_{ns}$, $\Delta = \Xi-\Omega D_3 \in M_{s,ns}$. Further, we have
\begin{equation*}\label{474}
T(F+D)T^{-1}=(I+D_3)(F+K_1)= F+K_1 +D_3 F +D_3 K_1.
\end{equation*}
The first $(n-1)$ block rows of $K_1$ are zero, and the last block column of $D_3$ is zero because $p<n$; hence, $D_3K_1=0\in M_{ns}$. Further, $D_3F=D_3(\mathcal{J}+K)=D_3\mathcal{J}+D_3 K$. We have $D_3\mathcal{J}=D_4$ and $D_3 K=0 \in M_{ns}$. Consequently,
\begin{equation*}
T(F+D)T^{-1}=F+K_1 +D_4.
\end{equation*}
Denote $F_1:=F+K_1$. Then,
\begin{equation}\label{480}
T(F+D)T^{-1}=F_1 +D_4.
\end{equation}
The matrix $F_1$ has the same form as the matrix $F$ and differs from it only by the last block row. Let  $\Omega_1:= \left[\begin{matrix}  -B_n & \ldots & -B_1 \end{matrix} \right] \in M_{s,ns}$ be the last block row of the matrix $F_1$. Then $\Omega_1=\Omega+ \Xi -\Omega D_3$, or, equivalently, $-\Omega_1=-\Omega- \Xi +\Omega D_3$. We rewrite the last equality of the block rows in terms of ``block coordinates'' of these rows starting from the last coordinate. We obtain
\begin{equation} \label{482}
\begin{aligned}
&B_1=A_1, \quad \ldots , \quad B_{n-p}=A_{n-p}, \\
&B_{n-p+1}=A_{n-p+1} - D_{np} - A_1 D_{n-1,p} - \ldots - A_{n-p} D_{pp} , \quad \ldots , \\
&B_{n}=A_{n} - D_{n1} - A_1 D_{n-1,1} - \ldots - A_{n-p} D_{p1}.
\end{aligned}
\end{equation}

By the inductive assumption, the assertion of Lemma \ref{lemma-01} holds for the matrix $F_1+D_4$, because, for the matrix 
$F_1+D_4$, the induction parameter $k$ is equal to $p+1$. Hence, there exists  a lower block triangular matrix $S$ such that
\begin{equation}\label{486}
S(F_1+D_4)S^{-1}=\Phi
\end{equation}
and the coefficients $\Gamma_i$ of $\Phi$ satisfy
\begin{equation}\label{488}
\Gamma_1= B_1- {\rm SP}_s (\mathcal{J}^{0} P D_4), \quad \ldots, \quad
\Gamma_n= B_n- {\rm SP}_s (\mathcal{J}^{n-1} P D_4).
\end{equation}
Construct $S_1:=ST$. Then, $S_1$ is lower block triangular non-degenerate and, by \eqref{480} and \eqref{486}, $S_1(F+D)S_1^{-1}=\Phi$.

The matrices $D$ and $D_4$ have the form of the matrices $C_1$ and $C_2$ in Lemma \ref{lemma-02}, respectively. Therefore, by Lemma \ref{lemma-02}, the following equalities hold:
\begin{equation} \label{490}
\begin{aligned}
&{\rm SP}_s (\mathcal{J}^{0} P D)={\rm SP}_s (\mathcal{J}^{0} P D_4),  \quad \ldots, \quad
{\rm SP}_s (\mathcal{J}^{n-p-1} P D)={\rm SP}_s (\mathcal{J}^{n-p-1} P D_4),  \\
&
{\rm SP}_s (\mathcal{J}^{n-p} P D)={\rm SP}_s (\mathcal{J}^{n-p} P D_4) + \sum_{i=p}^n A_{n-i}D_{ip}, \quad \ldots, \\
&{\rm SP}_s (\mathcal{J}^{n-1} P D)={\rm SP}_s (\mathcal{J}^{n-1} P D_4) + \sum_{i=p}^n A_{n-i}D_{i1}.
\end{aligned}
\end{equation}
By substituting \eqref{482} and \eqref{490} into \eqref{488}, we obtain
\begin{equation*} 
\begin{aligned}
&\Gamma_1=A_1-{\rm SP}_s (\mathcal{J}^{0} P D),  \quad \ldots, \quad
\Gamma_{n-p}=A_{n-p}-{\rm SP}_s (\mathcal{J}^{n-p-1} P D), \\
&
\Gamma_{n-p+1}=A_{n-p+1}-{\rm SP}_s (\mathcal{J}^{n-p} P D),  \quad \ldots, \quad
\Gamma_{n}=A_{n}-{\rm SP}_s (\mathcal{J}^{n-1} P D).
\end{aligned}
\end{equation*}
The proof of Lemma \ref{lemma-01} is complete.
\hfill $\square$

\begin{remark}
For the case $s=1$, Lemma \ref{lemma-01} has been proved in \cite[Lemma 4]{Zaitsev_2009}.
\end{remark}

From the matrix $F$, construct the following matrices: 
$$N_0:=I\in M_{ns}, \quad N_{\nu}:=N_{\nu-1} \cdot F + (L\otimes A_{\nu}) \in M_{ns}, \quad \nu=\overline{1,n-1},$$ where $L=I\in M_n$.
Then,
\begin{equation}\label{495}
\begin{aligned}
&N_0=I,\\
&N_1=F+ (I\otimes A_{1}),\\
&N_2=F^2+ (I\otimes A_{1})\cdot F + (I\otimes A_{2}),\\
&\ldots, \\
&N_{\nu}=F^{\nu}+ (I\otimes A_{1})\cdot F^{\nu-1} + (I\otimes A_{2})\cdot F^{\nu-2} + \ldots + 
(I\otimes A_{\nu-1})\cdot F + (I\otimes A_{\nu}),\\
&\ldots, \\
&N_{n-1}=F^{n-1}+ (I\otimes A_{1})\cdot F^{n-2} + (I\otimes A_{2})\cdot F^{n-3} + \ldots + 
(I\otimes A_{n-2})\cdot F + (I\otimes A_{n-1}).
\end{aligned}
\end{equation}

\begin{lemma}\label{lemma-03}
\it
Let a block matrix $D\in M_{ns}$ have the form \eqref{435}. Then,
\begin{equation}\label{500}
{\rm SP}_s (\mathcal{J}^{\nu} P D)={\rm SP}_s (N_{\nu}  D),
\end{equation}
for all $\nu=0,\ldots,n-1$.
\end{lemma}

\emph{Proof.}
From the form of the matrix $D$, it follows that ${\rm SP}_s (C D)=0\in M_{ns}$ for every strictly lower block triangular matrix $C=\{C_{ij}\}\in M_{ns}$, $C_{ij}\in M_s$, $i,j=\overline{1,n}$ (that is, a matrix such that $C_{ij}=0\in M_s$ for $i\le j$). 
Therefore, it suffices to show that the upper block triangular parts (that is, the diagonal and superdiagonal blocks) of the matrices $\mathcal{J}^{\nu} P$ and $N_{\nu}$ coincide for each $\nu=0,\ldots,n-1$. Consider the matrices $\mathcal{J}^{\nu} P$. We have
\begin{equation}\label{510}
\begin{gathered}
\mathcal{J}^{\nu} P=
\left[~\begin{matrix} P_{\nu} \\ \hline O_1 \end{matrix}~\right],\quad
P_{\nu}=\left[\begin{matrix}
A_{\nu} & \ldots & A_0 & 0 & \ldots & 0 \\
A_{{\nu}+1}&  \hdots & \hdots  & A_0 & \ldots & 0 \\
\hdots & \hdots & \hdots & \hdots & \hdots & \hdots \\
A_{n-1} & \hdots & \hdots & \hdots & \hdots  & A_0 
\end{matrix} \right] \in M_{(n-{\nu})s,ns},\\ O_1=0\in M_{{\nu s},ns}, \quad 0\in M_s.
\end{gathered}
\end{equation}
Consider the matrices $N_{\nu}$. Since $N_0=I$, it follows that \eqref{500} holds for $\nu=0$. 

The matrix $F$ has the form $F=\mathcal{J}+(e_n\otimes I)\Omega=(J\otimes I)+(e_n\otimes I)\Omega$ where $e_n\in\mathbb{R}^n$, $I\in M_s$, $e_n\otimes I\in M_{ns,s}$, $\Omega = \left[\begin{matrix}  -A_n & \ldots & -A_1 \end{matrix} \right] \in M_{s,ns}$. So, 
\begin{gather*}
N_1=F+ (I\otimes A_{1})=(J\otimes I)+(e_n\otimes I)\Omega + (I\otimes A_{1})=
(J^0\otimes A_1)+(J^1\otimes A_{0})+(e_n\otimes I)\Omega \\
=
\left[\begin{matrix}
A_1 & A_0 & 0 & \ldots & 0 & 0 \\
0 & A_1 & A_0 & \ldots & 0 & 0 \\
\hdots & \hdots & \hdots & \hdots & \hdots & \hdots \\
0 & 0 & 0 & \ldots &A_1 & A_0 \\
-A_{n} & -A_{n-1} & -A_{n-2} & \ldots & -A_{2} & 0 
\end{matrix}
\right].
\end{gather*}
Next,
\begin{gather*}
N_2=N_1 F + (I\otimes A_2) = N_1 \big((J\otimes I)+(e_n\otimes I)\Omega\big)  + (I\otimes A_2)  \\= N_1 (J\otimes I) + N_1 (e_n\otimes I)\Omega + (I\otimes A_2).
\end{gather*}
We have: $N_1 (e_n\otimes I)=(e_{n-1}\otimes I)$; $\Omega (J\otimes I) = \left[\begin{matrix}  0 & -A_n & \ldots & -A_2 \end{matrix} \right] \in M_{s,ns}$. Therefore,
$$
N_2=\Big((J^0\otimes A_1)+(J^1\otimes A_{0})+(e_n\otimes I)\Omega\Big)(J\otimes I)  + (e_{n-1}\otimes I)\Omega + (I\otimes A_2) $$ 
$$=
\Big((J^0\otimes A_2)+(J^1\otimes A_{1})+ (J^2\otimes A_{0})\Big) + \Big((e_{n-1}\otimes I)\Omega + (e_{n}\otimes I)\Omega (J\otimes I)\Big)  $$
$$=\left[~\begin{matrix}
A_2 & A_1 & A_0 & 0 & \ldots & 0 & 0 \\
0 & A_2 & A_1 & A_0 & \ldots & 0 & 0\\
\hdots & \hdots & \hdots & \hdots & \hdots & \hdots & \hdots \\
0 & 0 & 0 & 0 & \ldots  &A_2 & A_1 \\
0 & 0 & 0 & 0 & \ldots  & 0 & A_2 
\end{matrix}~\right] + \ 
\left[~\begin{matrix}
0 & 0 & \ldots & 0 & 0 \\
%\vdots & & & &\vdots\\
\hdots & \hdots & \hdots & \hdots & \hdots \\
0 & 0 & \ldots & 0 & 0 \\
-A_{n} & -A_{n-1} & \ldots & -A_{2} & -A_{1} \\
0 & -A_{n} & \ldots & -A_{3} & -A_{2} 
\end{matrix}~\right]$$
$$=\left[~\begin{matrix}
A_2 & A_1 & A_0 & 0 & \ldots & 0 & 0 & 0 \\
0 & A_2 & A_1 & A_0 & \ldots & 0 & 0 & 0\\
\hdots & \hdots & \hdots & \hdots & \hdots & \hdots & \hdots & \hdots \\
0 & 0 & 0 & 0 & \ldots &A_2 &A_1 & A_0  \\
\hline -A_{n} & -A_{n-1} & \hdots & \hdots & \hdots &  -A_{3} & 0 & 0 \\
0 &-A_{n} & \hdots & \hdots & \hdots & -A_{4} & -A_{3} & 0 
\end{matrix}~\right].
$$

Next, let us show by induction that
\begin{equation}\label{520}
\begin{gathered}
N_{\nu}=\Big((J^0\otimes A_{\nu})+(J^1\otimes A_{\nu-1})+ \ldots +(J^{\nu}\otimes A_{0})\Big) \\ + \Big((e_{n-\nu+1}\otimes I)\Omega + (e_{n-\nu+2}\otimes I)\Omega (J\otimes I)+ \ldots + (e_{n}\otimes I)\Omega (J^{\nu-1}\otimes I)\Big),
\end{gathered}
\end{equation}
for all ${\nu}=\overline{1,n-1}$. This formula is valid for ${\nu}=1,2$ as shown above.
Suppose that \eqref{520} holds for $\nu=k<n-1$. Let us show that it is true for $\nu=k+1$.
We have 
\begin{equation}\label{530}
\begin{gathered}
N_{k+1}=N_k \cdot F +(I\otimes A_{k+1})= N_k\big((J\otimes I)+(e_n\otimes I)\Omega\big) +
(I\otimes A_{k+1}) \\ = N_k (J\otimes I)+N_k (e_n\otimes I)\Omega + (I\otimes A_{k+1}) \\ =
\Big( \sum_{i=0}^k (J^i \otimes A_{k-i}) + \sum_{j=1}^k (e_{n-k+j}\otimes I) \Omega (J^{j-1} \otimes I) \Big) (J \otimes I) + N_k (e_n\otimes I)\Omega + (I\otimes A_{k+1}).
\end{gathered}
\end{equation}
By induction assumption, from the form of the matrix $N_k$, it follows that the last block column of $N_k$ is $(e_{n-k}\otimes I)$. Hence, 
\begin{equation}\label{532}
N_k (e_{n}\otimes I)= (e_{n-k}\otimes I).
\end{equation}
Expanding the brackets in \eqref{530}, we get
\begin{gather}\label{534}
(J^i \otimes A_{k-i}) \cdot (J \otimes I) = (J^{i+1} \otimes A_{k-i}),\\
\label{536}
(J^{j-1} \otimes I) \cdot (J \otimes I) = (J^{j} \otimes I).
\end{gather}
Thus, from \eqref{530}, \eqref{532}, \eqref{534}, and \eqref{536},  it follows that
\begin{equation*}%\label{538}
\begin{gathered}
N_{k+1}=\Big((J^0\otimes A_{k+1})+(J^1\otimes A_{k})+ \ldots +(J^{k+1}\otimes A_{0})\Big) \\ + \Big((e_{n-k}\otimes I)\Omega + (e_{n-k+1}\otimes I)\Omega (J\otimes I)+ \ldots + (e_{n}\otimes I)\Omega (J^{k}\otimes I)\Big).
\end{gathered}
\end{equation*}
Hence, the formula \eqref{520} is true for all ${\nu}=\overline{1,n-1}$. By using this formula, we construct the matrix
\begin{equation}\label{540}
N_{\nu}=\left[~\begin{matrix}
A_{\nu} & A_{{\nu}-1} & \hdots & \hdots  & A_0 & 0 & \ldots& 0\\
0 & A_{{\nu}} & \hdots & \hdots  & A_1 & A_0 & \ldots & 0\\
\hdots &  \hdots &  \hdots &  \hdots &  \hdots &  \hdots &  \hdots &  \hdots \\
0 & \hdots & 0  & A_{\nu} & \hdots & \hdots & A_1 & A_0\\
\hline
-A_n & -A_{n-1} & \ldots & -A_{{\nu}+1} & 0 & \hdots & \hdots  & 0\\
0 & -A_n &  \ldots & -A_{{\nu}+2} & -A_{{\nu}+1} & 0 & \ldots & 0\\
\hdots &  \hdots &  \hdots &  \hdots &  \hdots &  \hdots &  \hdots &  \hdots \\
0 & \ldots & -A_n & \hdots & \hdots & \hdots  & -A_{{\nu}+1} & 0 
\end{matrix}~\right].
\end{equation}
It follows from \eqref{510} and \eqref{540} that the upper block triangular parts of  $\mathcal{J}^{\nu} P$ and $N_{\nu}$ coincide. The proof of the lemma is complete.
\hfill $\square$

\begin{remark}
For the case $s=1$, Lemma \ref{lemma-03} has been proved in \cite[Lemma 2]{Zaitsev_2010}.
\end{remark}

\begin{lemma}\label{lemma-04}
\it
Let matrices $F$ and $D$ have the form \eqref{400} and \eqref{435}, respectively.
Then, there exists a lower block triangular matrix $S$ such that the following holds:
the matrix $S(F+D)S^{-1}$ is equal to the matrix $\Phi$ of \eqref{064}, the block matrix coefficients of which are related to the coefficients  of the matrices $F$ and $D$ as follows:
\begin{gather}\label{640}
\Gamma_i=A_i-{\rm SP}_s(N_{i-1} D), \quad i=\overline{1,n}.
\end{gather}
\end{lemma}

Lemma \ref{lemma-04} follows from Lemmas \ref{lemma-01} and \ref{lemma-03}.

Denote
\begin{equation}\label{645}
T_1:={\rm SP}_s(D), \quad T_2:={\rm SP}_s( F D), \quad \ldots, \quad T_n:={\rm SP}_s( F^{n-1} D), 
\end{equation}
\begin{equation}\label{650}
\widehat{\Gamma}:=\left[\begin{matrix} \Gamma_1 \\ \Gamma_2 \\ \vdots \\ \Gamma_n  \end{matrix}\right] \in M_{ns,s}, \quad
\widehat{A}:=\left[\begin{matrix} A_1 \\ A_2 \\ \vdots \\ A_n  \end{matrix}\right]\in M_{ns,s}, \quad
\widehat{T}:=\left[\begin{matrix} T_1 \\ T_2 \\ \vdots \\ T_n  \end{matrix}\right] \in M_{ns,s}.
\end{equation}

\begin{lemma} \label{lemma-05}
\it
Equalities \eqref{640} are equivalent to the equality
\begin{gather}\label{670}
\widehat\Gamma=\widehat A- P \widehat T.
\end{gather}
\end{lemma}

\emph{Proof.}
Equality \eqref{670} has the form
\begin{equation} \label{673}
\left[\begin{matrix} \Gamma_1 \\ \Gamma_2 \\ \vdots \\ \Gamma_n  \end{matrix}\right]=
\left[\begin{matrix} A_1 \\ A_2 \\ \vdots \\ A_n  \end{matrix}\right]-
\begin{bmatrix}
A_0&0&\ldots&0\\
A_1&A_0&\ldots&0\\
\vdots&\ddots&\ddots&\vdots\\
A_{n-1}&\ldots&A_1&A_0
\end{bmatrix} \cdot
\left[\begin{matrix} {\rm SP}_s(D) \\ {\rm SP}_s( F D) \\ \vdots \\ {\rm SP}_s( F^{n-1} D)  \end{matrix}\right].
\end{equation}
It is easy to see that the following equality is true: ${\rm SP}_s\big((I\otimes B) \cdot C\big)=B \cdot  {\rm SP}_s(C)$  where $I\in M_n$, $B\in M_s$, $C\in M_{ns}$. Hence, for any $i=\overline{0,n-1}$, $\nu=\overline{0,n-1}$, 
\begin{equation} \label{675}
{\rm SP}_s\big((I\otimes A_i) \cdot F^{\nu}D\big)=A_i \cdot  {\rm SP}_s(F^{\nu} D).
\end{equation}
Due to \eqref{675} and \eqref{495}, we obtain that \eqref{673} is equivalent to \eqref{640}. Q.E.D.
\hfill $\square$

\begin{remark}
Since $P$ is non-degenerate, then, for any $\widehat\Gamma$, one can express $\widehat{T}$ from \eqref{670}:
\begin{equation} \label{680}
\widehat{T}=P^{-1}(\widehat{A}-\widehat{\Gamma}).
\end{equation}
\end{remark}

Consider the system \eqref{080}, \eqref{082} with coefficients \eqref{400}, \eqref{410}, \eqref{420}. From \eqref{410} and \eqref{420}, it follows that the matrix $GQH$ has the form \eqref{435} of the matrix $D$. Let us replace the matrix $D$ in the equalities \eqref{645} with $GQH$.
Then, equalities \eqref{645} take the form
\begin{equation}\label{690}
\begin{aligned}
&T_1={\rm SP}_s(GQH),\\
&T_2={\rm SP}_s(FGQH),\\
& \ldots \ldots \ldots \ldots \ldots \ldots , \\
&T_n={\rm SP}_s(F^{n-1}GQH).
\end{aligned}
\end{equation}

Denote 
\begin{gather}
\label{700}
v:={\rm vecc}\,({\rm VecCR}_s Q)\in\mathbb{K}^{kms^2}, \\
\label{710}
w:={\rm
col}\,\big({\rm vecc}\,(T_1),\ldots,{\rm
vecc}\,(T_n)\big)\in \mathbb{K}^{ns^2}.
\end{gather} 

\begin{lemma} \label{lemma-07}
\it
Equalities \eqref{690} are equivalent to the equality
\begin{gather}\label{720}
\Theta v= w.
\end{gather}
\end{lemma}

\emph{Proof.}
Let us apply the mapping $\rm vecc$ to equalities \eqref{690}
and apply Lemma \ref{lem-unroll} to the matrices $X=F^{i-1}G$, $Y=H$. Then, equalities \eqref{690} take the form
\begin{equation}\label{725}
{\rm vecc}\, (T_i)={\rm VecRR}_{s^2}\big( (H^{\mathcal{T}})^T\star
F^{i-1}G\big)\cdot {\rm vecc}\,({\rm VecCR}_s Q),\quad
i=\overline{1,n}.
\end{equation}
In the notations \eqref{425}, \eqref{700}, \eqref{710}, system \eqref{725} takes the form \eqref{720}. Q.E.D.
\hfill $\square$

\begin{lemma} \label{lemma-10}
\it
System \eqref{720} is resolvable with respect to $v$ for arbitrary $w$ iff condition \eqref{428} holds.
\end{lemma}
The proof of Lemma \ref{lemma-10} is clear. If condition \eqref{428} holds, then system \eqref{720} has the particular solution
\begin{equation}\label{730}
v=\Theta^T(\Theta\Theta^T)^{-1}w. 
\end{equation}

Let us carry out the final arguments to prove Theorem \ref{teo-AMCA-LSOF-1}. Let condition \eqref{428} hold. We show that the problem of AMCA for CMP by LSOF is resovable. Let arbitrary matrices $\Gamma_i\in M_s(\mathbb{K})$ be given. One needs to construct a matrix $Q\in M_{ms,ks}(\mathbb{K})$ and a matrix $S\in M_{ns}(\mathbb{K})$ such that $S(F+GQH)S^{-1}=\Phi$ and $\Phi$ has the given block coefficients $-\Gamma_{n+1-i}$ in the last block row. By Lemma \ref{lemma-04}, for this it is sufficient to construct the matrix $Q$, which ensures the fulfillment of the equalities
\begin{gather}\label{735}
\Gamma_i=A_i-{\rm SP}_s(N_{i-1} GQH), \quad i=\overline{1,n}.
\end{gather}

Construct the matrices $\widehat A$ and $\widehat \Gamma$ by using formula \eqref{650}. By using \eqref{680}, let us construct the matrix $\widehat T$. Construct the vector $w$ by using formula \eqref{710}. Let us resolve the system \eqref{720} by formula \eqref{730}. From \eqref{730}, we find the matrix $Q$ using the inverse formula to \eqref{700}:
\begin{equation}\label{740}
Q={\rm VecCR}_s^{-1}({\rm vecc}^{-1}v).
\end{equation}
Then, taking into account Lemma \ref{lemma-05},  equalities \eqref{735} are satisfied. The proof of Theorem \ref{teo-AMCA-LSOF-1}
is complete.

\begin{remark}
Let us suppose that the blocks of the matrix $F$ are scalar matrices, i.e., $A_i=a_i I$, $a_i \in\mathbb{K}$, $i=\overline{1,n}$, $I\in M_s$. Hence, the matrix $F$ has the form $F=F_0\otimes I$, where $F_0=J+e_n\psi \in M_n$, $\psi=[-a_n,\ldots,-a_1]\in M_{1,n}$, $I\in M_s$. In this case, the matrices $ (H^{\mathcal{T}})^T\star F^{i-1}G$ in \eqref{425} can be replaced by $ (H^{\mathcal{T}})^T F^{i-1} \star G$, $i=\overline{1,n}$. This assertion follows from Lemma \ref{lem-F-star}.
\end{remark}

\begin{remark}\label{remark-main}
Theorem \ref{teo-AMCA-LSOF-1}  has been proven previously in special cases: $(a)$ for the case $s=1$ (see \cite[Theorem 1]{Zaitsev_2010}); $(b)$ for the case when the system \eqref{080}, \eqref{082}  has the form \eqref{210}, \eqref{220} (see \cite[Theorem 3]{Zaitsev_2021}). In these cases, it was proven that condition \eqref{428} is also necessary for the problem of AMCA for CMP by LSOF to be resolvable. Here, in Theorem \ref{teo-AMCA-LSOF-1}, we cannot yet assert that condition \eqref{428} is necessary for the solvability of this problem. The reason is that, for the matrix $F+GQH$, the matrix $\Phi$, which is similar to it and is a lower block Frobenius matrix, is not uniquely determined (in contrast to the case when $s = 1$). It may turn out that, for the problem to be solvable, a set of matrices $\{(\Gamma_1,\ldots,\Gamma_n)\}$ does not necessarily need to range the entire set $M_s\times \ldots \times M_s$. For example, we have a hypothesis (which has not yet been proven): 
{\it for system \eqref{080}, \eqref{082}, \eqref{084}, the problem of AMCA for CMP by LSOF is resolvable iff for any $\Gamma_i\in T_s(\mathbb{K})$, $i=\overline{1,n}$, there exists a gain matrix $Q\in M_{ms,ks}(\mathbb{K})$  such that the the matrix $F+GQH$ of the system \eqref{088} is similar to the matrix \eqref{064}}; here $T_s(\mathbb{K})\subset M_s(\mathbb{K})$ is a subset of  (upper) triangular matrices. The question of the necessity of condition \eqref{428} in Theorem~\ref{teo-AMCA-LSOF-1} requires additional study.
\end{remark}

\begin{remark}\label{remark-5}
\textrm{Note that the condition $mk\ge n$ is necessary for \eqref{428}.}
\end{remark}

Based on the proof of Theorem \ref{teo-AMCA-LSOF-1}, we present an algorithm for solving the problem of AMCA for CMP by LSOF.

{\bf Algorithm 1.} Let the system
\eqref{080}, \eqref{082} with coefficients \eqref{400}, \eqref{410}, \eqref{420} be given.

1. Construct the matrix $\Theta$ of \eqref{425}.

2. Check the condition \eqref{428}. If this condition is satisfied, then the problem is solvable.

3. Let arbitrary matrices $\Gamma_1,\ldots,\Gamma_n\in M_s$ be given.

4. Construct the matrices $P$ of \eqref{430} and the matrices $\widehat{\Gamma}$ and $\widehat{A}$ of \eqref{650}.

5. Calculate the matrix $\widehat{T}$ using the equality \eqref{680}.

6. The matrix $\widehat{T}$ has the form \eqref{650}. From $\widehat{T}$,  find the matrices $T_1, \ldots, T_n \in M_s$.

7. From the matrices $T_1, \ldots, T_n \in M_s$, construct the vector $w$ using the formula \eqref{710}.

8. Solve the system \eqref{720} with respect to the vector $v$; for example, using the formula \eqref{730}.

9. Applying the formula \eqref{740} to the vector $v$,  find the gain matrix $Q$.

10. Denote the matrix $F+GQH$ of the closed-loop system by $Z$. This matrix is a lower block Hessenberg matrix. By Lemma \ref{lemma-00}, from the matrix $Z$, construct the matrices $S_1,\ldots,S_n$ and the matrix $S$. Then, this matrix $S$ reduces the matrix $Z=F+GQH$ of the closed-loop system to the matrix $\Phi$ of \eqref{064}, i.e., $S(F+GQH)S^{-1}=\Phi$.

Let us demonstrate the algorithm with an example.

\begin{example}
\label{ex0}
Consider system \eqref{080}, \eqref{082} with $n=3$, $s=2$, $m=k=p=2$ with the following matrices:
\begin{equation}\label{ex1-1720}
F=
\left[\begin{matrix}
0 & 0 & 1 & 0 & 0 & 0\\
0 & 0 & 0 & 1 & 0 & 0\\
0 & 0 & 0 & 0 & 1 & 0\\
0 & 0 & 0 & 0 & 0 & 1\\
2 & 0 & -1 & 0 & 0 & 0\\
0 & 0 & 0 & -1 & 0 & -1
\end{matrix} \right], \quad
{G}=\begin{bmatrix} 0 & 0 & 0 & 0\\
0 & 0 & 0 & 0\\
0 & 1 & 1 & 0 \\
0 & -1 & 0 & 1 \\
-1 & 3 & 1 & -1 \\
2 & 1 & 0 & -1
\end{bmatrix},  \quad
{H}=\begin{bmatrix} 
1 & 0 & 0 & 0 & 0 & 0 \\
0 & 0 & 0 & 1 & 0 & 0 \\
1 & 0 & 0 & 0 & 0 & 0 \\
0 & -1 & 1 & 0 & 0 & 0
\end{bmatrix}.
\end{equation}
The matrices \eqref{ex1-1720} have the form \eqref{400}, \eqref{410}, \eqref{420}.
From \eqref{ex1-1720}, we have
\begin{equation}\label{ex1-1730}
 F {G} =
\left[\begin{matrix}
0 & 1 & 1 & 0 \\
0 & -1 & 0 & 1 \\
-1 & 3 & 1 & -1 \\
2 & 1 & 0 & -1 \\
0 & -1 & -1 & 0 \\
-2 & 0 & 0 & 0
\end{matrix} \right], \qquad
F^2 {G}=\begin{bmatrix} 
-1 & 3 & 1 & -1 \\
2 & 1 & 0 & -1 \\
0 & -1 & -1 & 0 \\
-2 & 0 & 0 & 0 \\
1 & -1 & 1 & 1 \\
0 & -1 & 0 & 1
\end{bmatrix}.
\end{equation}
From \eqref{ex1-1720} and \eqref{ex1-1730}, we obtain:
\begin{gather*}
(H^{\mathcal{T}})^T\star {G}=\left[\begin{matrix}X_{11}&X_{12}\\X_{21}&X_{22}\end{matrix}\right],\\
X_{11}=\left[\begin{matrix}0&0&0&0\\0&0&0&0\\0&0&0&1\\0&0&0&-1\end{matrix}\right],\quad 
X_{12}=\left[\begin{matrix}0&0&0&0\\0&0&0&0\\0&0&1&0\\0&0&0&1\end{matrix}\right],\\
X_{21}=\left[\begin{matrix}0&0&0&1\\0&0&0&-1\\0&0&0&0\\0&0&0&0\end{matrix}\right],\quad 
X_{22}=\left[\begin{matrix}0&0&1&0\\0&0&0&1\\0&0&0&0\\0&0&0&0\end{matrix}\right],
\end{gather*}
\begin{gather*}
(H^{\mathcal{T}})^T\star {F}{G}=\left[\begin{matrix}Y_{11}&Y_{12}\\Y_{21}&Y_{22}\end{matrix}\right],\\
Y_{11}=\left[\begin{matrix}0&1&0&0\\0&-1&0&0\\0&0&-1&3\\0&0&2&1\end{matrix}\right],\quad 
Y_{12}=\left[\begin{matrix}1&0&0&0\\0&1&0&0\\0&0&1&-1\\0&0&0&-1\end{matrix}\right],\\
Y_{21}=\left[\begin{matrix}0&1&-1&3\\0&-1&2&1\\0&0&0&-1\\0&0&0&1\end{matrix}\right],\quad 
Y_{22}=\left[\begin{matrix}1&0&1&-1\\0&1&0&-1\\0&0&-1&0\\0&0&0&-1\end{matrix}\right],
\end{gather*}
\begin{gather*}
(H^{\mathcal{T}})^T\star {F}^2 {G}=\left[\begin{matrix}Z_{11}&Z_{12}\\Z_{21}&Z_{22}\end{matrix}\right], \\
Z_{11}=\left[\begin{matrix}-1&3&0&0\\2&1&0&0\\0&0&0&-1\\0&0&-2&0\end{matrix}\right],\quad 
Z_{12}=\left[\begin{matrix}1&-1&0&0\\0&-1&0&0\\0&0&-1&0\\0&0&0&0\end{matrix}\right],\\
Z_{21}=\left[\begin{matrix}-1&3&0&-1\\2&1&-2&0\\0&0&1&-3\\0&0&-2&-1\end{matrix}\right],\quad 
Z_{22}=\left[\begin{matrix}1&-1&-1&0\\0&-1&0&0\\0&0&-1&1\\0&0&0&1\end{matrix}\right].
\end{gather*}
Constructing matrix \eqref{425}, we get
\begin{equation}\label{ex1-1740}
\Theta=\left[\begin{matrix}
X_{11}&X_{12}&X_{21}&X_{22}\\Y_{11}&Y_{12}&Y_{21}&Y_{22}\\Z_{11}&Z_{12}&Z_{21}&Z_{22}\end{matrix}\right].
\end{equation}
Calculating the rank of the matrix \eqref{ex1-1740}, we obtain that ${\rm rank}\,\Theta=12$. So, conditions of Theorem \ref{teo-AMCA-LSOF-1} are fulfilled. Hence, by Theorem \ref{teo-AMCA-LSOF-1}, for system \eqref{080}, \eqref{082} with coefficients \eqref{ex1-1720}, the problem of AMCA for CMP by LSOF is resolvable. Let us construct this feedback control. Suppose, for example, that
\begin{equation*}%\label{ex-1750}
\Gamma_1 =
\left[\begin{matrix}
6 & 0 \\
0 & 6
\end{matrix} \right], \qquad
\Gamma_2 =
\left[\begin{matrix}
11 & 0 \\
0 & 11
\end{matrix} \right], \qquad
\Gamma_3 =
\left[\begin{matrix}
6 & 0 \\
0 & 6
\end{matrix} \right].
\end{equation*}
We have 
\begin{equation*}%\label{ex1-1780}
F=\left[\begin{matrix}0&I&0\\0&0&I\\-A_3&-A_2&-A_1\end{matrix}\right],\quad 
A_1=\left[\begin{matrix}0&0\\0&1\end{matrix}\right],\quad A_2=\left[\begin{matrix}1&0\\0&1\end{matrix}\right],\quad A_3=\left[\begin{matrix}-2&0\\0&0\end{matrix}\right].
\end{equation*}
Construct the matrix $P=\left[\begin{matrix}I&0&0\\A_1&I&0\\A_2&A_1&I\end{matrix}\right]$ of \eqref{430} and $\widehat{\Gamma}=\left[\begin{matrix}\Gamma_1 \\ \Gamma_2 \\ \Gamma_3\end{matrix}\right]$ and 
$\widehat{A}=\left[\begin{matrix}A_1 \\ A_2 \\ A_3\end{matrix}\right]$ of  \eqref{650}. Calculating \eqref{680}, we obtain that
\begin{equation*}
\widehat{T}=\left[\begin{matrix} T_1 \\ T_2 \\ T_3 \end{matrix}\right], \quad
T_1=\left[\begin{matrix}-6&0\\0&-5\end{matrix}\right],\quad T_2=\left[\begin{matrix}-10&0\\0&-5\end{matrix}\right],\quad T_3=\left[\begin{matrix}-2&0\\0&4\end{matrix}\right].
\end{equation*}
Construct $w$ by formula \eqref{710}. Then
$w={\rm col}\,(-6,0,0,-5,-10,0,0,-5,-2,0,0,4)$. Now, resolving the system $\Theta v =w$ of \eqref{720} by the formula \eqref{730}, we obtain that 
$$
v={\rm col}\,(0,-2,-5,-16,-3,-5,16,-21,0,-2,3,9,-3,-5,-15,9).
$$ From \eqref{740}, we get
\begin{equation}\label{ex1-1820}
Q=\left[\begin{matrix}0 & -5 & 0 &3 \\ -2 & -16 & -2 & 9 \\ -3 & 16 & -3 & -15 \\ -5 & -21 & -5 & 9
\end{matrix}\right].
\end{equation}
The controller $u=Qy$ with the matrix \eqref{ex1-1820} leads the system \eqref{080}, \eqref{082} with coefficients \eqref{ex1-1720} to the closed-loop system with the matrix
\begin{equation*}%\label{ex1-1830}
F+GQH=\left[\begin{matrix} 0 & 0 & 1 & 0 & 0 & 0 \\  0 & 0 & 0 & 1 & 0 & 0 \\ -10 & 6 & - 6 & 0 & 1 & 0 \\ -6 & 0 & 0 & -5 & 0 & 1 \\ -6 & 0 & -1 &-6 & 0 & 0 \\  6 & -6 & 6 & -6 & 0 & -1
\end{matrix}\right].
\end{equation*}
Applying step 10 of Algorithm 1, we obtain that 
\begin{equation*}%\label{ex1-1840}
S=\left[\begin{matrix} I  & 0  &  0\\  0 & I & 0  \\ S_{31} & S_{32} & I
\end{matrix}\right], \quad
S_{31}=\left[\begin{matrix} -10  & 6  \\  -6 &  0  
\end{matrix}\right], \quad
S_{32}=\left[\begin{matrix} -6  & 0  \\  0 &  -5  
\end{matrix}\right].
\end{equation*}
One can check that the following equality holds: $\Phi ={{S}}\big(F + G  Q H){{S}}^{-1}$.
\end{example}

\begin{remark}
The property of AMCA allows simultaneously assign eigenvalues and eigenvectors, i.e., eigenstructure (see \cite[Theorem 11]{Zaitsev_2021} for the case when all eigenvalues are different and \cite[Theorem 8]{Zaitsev_2024-2} for the case when not all eigenvalues are necessarily different). System \eqref{070} with the matrix \eqref{064} is equivalent to the differential eqiation \eqref{055} in the space $\mathbb{K}^n$. Consider the following task. It is required to construct a matrix $Q$ such that the closed-loop system \eqref{088} is equivalent to the differential equation \eqref{055} having a given basis of solutions. What could this basis be? In \cite{Zaitsev_2021,Zaitsev_2024-2}, it is shown that it can be quite arbitrary, namely. 
Let an arbitrary set of linearly independent vectors  $h_1,\ldots,h_s\in\mathbb{K}^s$ be given and an arbitrary  list $\Omega=(\lambda_1,\lambda_2,\ldots,\lambda_{ns})$ of $ns$ (not necessarily different) numbers $\lambda_{\xi}\in \mathbb{K}$ be given such that following vector functions are linearly independent:
\begin{equation}\label{7-200}
\begin{aligned}
\psi_{1,1}(t) & =h_1 e^{\lambda_1 t}, &  \psi_{1,2}(t)&=h_2 e^{\lambda_2 t}, & \ldots & \ldots & \psi_{1,s}(t)&=h_s e^{\lambda_s t}, \\
\psi_{2,1}(t) & =h_1 e^{\lambda_{s+1} t}, &  \psi_{2,2}(t)&=h_2 e^{\lambda_{s+2} t}, & \ldots & \ldots & \psi_{2,s}(t)&=h_s e^{\lambda_{2s} t}, \\
\ldots \ldots & \ldots \ldots \ldots \ldots  & 
\ldots \ldots & \ldots \ldots \ldots  \ldots  & \ldots & \ldots & 
\ldots \ldots & \ldots \ldots \ldots  , \\
\psi_{n,1}(t) & =h_1 e^{\lambda_{(n-1)s+1} t}, &  \psi_{n,2}(t)&=h_2 e^{\lambda_{(n-1)s+2} t}, & \ldots & \ldots & \psi_{n,s}(t)&=h_s e^{\lambda_{ns} t}.
\end{aligned}
\end{equation}
Then, the set \eqref{7-200} may serve as such a basis.

For example, let $n=3$, $s=2$, $h_1=\left[\begin{matrix} 1 \\ 0 \end{matrix}\right]$, $h_2=\left[\begin{matrix} 0 \\ 1 \end{matrix}\right]$, $\Omega=(-1,-1,-2,-2,-3,-3)$. Then, the set \eqref{7-200} is
\begin{equation}\label{7-205}
\begin{aligned}
\psi_{1,1}(t) &=\left[\begin{matrix} 1 \\ 0 \end{matrix}\right] e^{- t}, & 
\psi_{1,2}(t) &=\left[\begin{matrix} 0 \\ 1 \end{matrix}\right] e^{- t},\\
\psi_{2,1}(t) &=\left[\begin{matrix} 1 \\ 0 \end{matrix}\right] e^{- 2t}, & 
\psi_{2,2}(t) &=\left[\begin{matrix} 0 \\ 1 \end{matrix}\right] e^{- 2t},\\
\psi_{3,1}(t) &=\left[\begin{matrix} 1 \\ 0 \end{matrix}\right] e^{- 3t}, & 
\psi_{3,2}(t) &=\left[\begin{matrix} 0 \\ 1 \end{matrix}\right] e^{- 3t}.
\end{aligned}
\end{equation}
The vector functions \eqref{7-205} are linearly independent. 
Using the proof of \cite[Theorem 11]{Zaitsev_2021} (and  \cite[Theorem 8]{Zaitsev_2024-2}), construct a differential equation \eqref{055} that has the set \eqref{7-205} as its basis of solutions: then, we obtain that
\begin{equation}\label{7-210}
\Gamma_1 =
\left[\begin{matrix}
6 & 0 \\
0 & 6
\end{matrix} \right], \qquad
\Gamma_2 =
\left[\begin{matrix}
11 & 0 \\
0 & 11
\end{matrix} \right], \qquad
\Gamma_3 =
\left[\begin{matrix}
6 & 0 \\
0 & 6
\end{matrix} \right].
\end{equation}
In Example \ref{ex0}, we constructed a matrix $Q$ that provides the similarity of the matrix of system \eqref{088} to matrix \eqref{064} with coefficients \eqref{7-210}. Thus, in Example \ref{ex0}, for system \eqref{080}, \eqref{082} with matrices \eqref{ex1-1720}, LSOF is constructed, which ensures the following property: The closed-loop system \eqref{088} is equivalent to the differential equation \eqref{055} having the basis of solutions \eqref{7-205}. In particular, the system \eqref{088} (and the equation \eqref{055}) is exponentially stable.
\end{remark}

\section{Theorem 3 generalizes sufficient conditions in Theorems 2 and 3 of \cite{Zaitsev_2021}
\label{generalization-LAIA} }

Consider system \eqref{210}, \eqref{220}. Suppose that the control in system \eqref{210}, \eqref{220} has the form of linear static output feedback \eqref{230} where $Q=\{Q_{\alpha \beta}\}\in M_{ms,ks}(\mathbb{K})$, $Q_{\alpha \beta}\in M_s(\mathbb{K})$, $\alpha=\overline{1,m}$, $\beta=\overline{1,k}$.
The closed-loop system~\eqref{210}, \eqref{220}, \eqref{230}  takes the form
\begin{gather}\label{750}
x^{(n)}+\sum\limits_{i=1}^n
A_{i}x^{(n-i)}-\sum\limits_{\alpha=1}^m\sum\limits_{l=p}^n
B_{l\alpha}\left(\sum\limits_{\beta=1}^k
Q_{\alpha\beta}\Big(\sum\limits_{\nu=1}^p
{C}_{\nu\beta}x^{(n-l+\nu-1)}\Big)\right)=0.
\end{gather}

\begin{definition}\label{def-07} We say that system
\eqref{210}, \eqref{220}  is {\it arbitrary matrix coefficient assignable by LSOF} \eqref{230} if for any $\Gamma_i\in M_s(\mathbb{K})$, $i=\overline{1,n}$, there exists a gain matrix $Q\in M_{ms,ks}(\mathbb{K})$  such that the closed-loop system \eqref{750}  has the form \eqref{055}.
\end{definition}

\begin{remark}
\textrm{Definition \ref{def-07} was given in \cite[Definition 2]{Zaitsev_2021} but the property of arbitrary matrix coefficient assignability (AMCA) was named as arbitrary matrix eigenvalue spectrum assignability (AMESA). As noted in Introduction, the term AMESA is not accurate and it is more correct to use the term AMCA.}
\end{remark}
 
From system \eqref{210}, \eqref{220}, construct the following block matrices:
\begin{equation}\label{760}
\mathcal{A}=
\left[\begin{matrix}
0 & I & 0 & \ldots & 0 \\
0 & 0 & I & \ldots & 0 \\
\vdots & \vdots & \vdots & \ddots & \vdots  \\
0 & 0 & 0 & \ldots & I \\
-A_n & -A_{n-1} & -A_{n-2} & \ldots & -A_1
\end{matrix} \right].
\end{equation}
\begin{equation} \label{770}
\mathcal{B}=\left[\begin{matrix}
O_1 \\ \widehat{\mathcal{B}}
\end{matrix}
\right], \quad
\widehat{\mathcal{B}}
=\begin{bmatrix}
{B}_{p1} &  \ldots & {B}_{pm} \\
\vdots &  & \vdots  \\
{B}_{n1} &  \ldots & {B}_{nm} 
\end{bmatrix}, \quad
\mathcal{C}=\left[\begin{matrix}
\widehat{\mathcal{C}} \\ O_2
\end{matrix}
\right], \quad
\widehat{\mathcal{C}}
=\begin{bmatrix}
{C}_{11} &  \ldots & {C}_{1k} \\
\vdots &  & \vdots  \\
{C}_{p1} &  \ldots & {C}_{pk} 
\end{bmatrix}.
\end{equation}
Here $O_1=0\in M_{(p-1)s,ms}$, $O_2=0\in M_{(n-p)s,ks}$, $\mathcal{A}\in M_{ns}$, $\mathcal{B}\in M_{ns,ms}$, $\mathcal{C}\in M_{ns,ks}$.

Consider the matrices
\begin{equation*}
\mathcal{C}^T\star \mathcal{B}, \quad \mathcal{C}^T\star \mathcal{J}\mathcal{B}, \quad \ldots, \quad \mathcal{C}^T\star \mathcal{J}^{n-1}\mathcal{B}.
\end{equation*}
We have  $\mathcal{C}^T\in M_{ks,ns}$, $\mathcal{B}\in M_{ns,ms}$, hence, $\mathcal{C}^T\star \mathcal{J}^{i-1}\mathcal{B} \in M_{ks^2,ms^2}$ for all $i=\overline{1,n}$. Let us construct the matrices ${\rm VecRR}_{s^2}(\mathcal{C}^T\star \mathcal{J}^{i-1}\mathcal{B}) \in M_{s^2,kms^2}$, $i=\overline{1,n}$, and the matrix
\begin{equation*}%\label{780}
\Theta_0=\left[\begin{matrix}
{\rm VecRR}_{s^2}(\mathcal{C}^T\star \mathcal{B})\\
{\rm VecRR}_{s^2}(\mathcal{C}^T\star \mathcal{J}\mathcal{B})\\
\ldots\ldots\ldots\ldots\ldots\ldots\ldots \\
{\rm VecRR}_{s^2}(\mathcal{C}^T\star \mathcal{J}^{n-1}\mathcal{B})
\end{matrix}\right]\in M_{ns^2,kms^2}.
\end{equation*}

The following theorems have been proved in \cite{Zaitsev_2021}.

\begin{theorem}[{see \cite[Theorem 2]{Zaitsev_2021}}] \label{teo-SP-2021}
\it
System \eqref{210}, \eqref{220}
 is AMCA by LSOF \eqref{230} if and only if  for any $\Gamma_i\in M_s(\mathbb{K})$, $i=\overline{1,n}$, there exists a  matrix $Q\in M_{ms,ks}(\mathbb{K})$ such that the following equalities hold:
\begin{gather}\label{790}
\Gamma_i=A_i-{\rm SP}_s\,(\mathcal{J}^{i-1}\mathcal{B}Q\mathcal{C}^{\mathcal{T}}), \quad i=\overline{1,n}.
\end{gather}
\end{theorem}

\begin{theorem}[{see \cite[Theorem 3]{Zaitsev_2021}}]\label{teo-AMESA-LSOF-1-2021}
\it
System \eqref{210}, \eqref{220}
 is AMCA by LSOF \eqref{230} if and only if
\begin{equation}
\label{800}
{\rm rank}\,\Theta_0 = ns^2.
\end{equation}
\end{theorem}

We will show that the sufficient conditions in Theorems \ref{teo-SP-2021} and \ref{teo-AMESA-LSOF-1-2021} are special cases of the results proved in Section~\ref{section-proof-1}.

From \eqref{760} and \eqref{770}, it follows that the   matrix $\mathcal{A}$ has the form \eqref{400}, the matrix $\mathcal{B}$ has the form \eqref{410}, the matrix $\mathcal{C}^{\mathcal{T}}$ has the form \eqref{420}. From the matrix $\mathcal{A}$, construct the matrix $P$  by the formula \eqref{430}. The matrix $P$ is lower block triangle, hence, $P^{-1}$ is also lower block triangle.
Construct the matrices 
\begin{equation}
\label{810}
\Sigma:=\mathcal{A}, \quad 
\Psi:=P^{-1}\mathcal{B}, \quad
\Omega:=\mathcal{C}^{\mathcal{T}}.
\end{equation}
and consider the system
\begin{gather}\label{820}
\dot{z}=\Sigma z+ \Psi u,\quad z\in \mathbb{K}^{ns}, \quad u\in\mathbb{K}^{ms},\\
\label{830}
y=\Omega z, \quad y\in \mathbb{K}^{ks}.
\end{gather}
It is clear that the  matrices \eqref{810} have the form \eqref{400}, \eqref{410}, \eqref{420} because $P^{-1}$ is  lower block triangle. 
Thus, system \eqref{820}, \eqref{830} with matrices \eqref{810} is a special case of system \eqref{080}, \eqref{082} with matrices \eqref{400}, \eqref{410}, \eqref{420}.

\begin{theorem}\label{teo-7-1}
\it
System  \eqref{210}, \eqref{220} is equivalent to system \eqref{820},  \eqref{830} with matrices \eqref{810}.
\end{theorem}

\emph{Proof.}
The proof method is similar to the proof of \cite[Theorem 4]{Zaitsev_2024}.
Denote $\Psi=:\!\{\Psi_{i\alpha}\}$, $\Psi_{i \alpha}\in M_s$, $i=\overline{1,n}$, $\alpha=\overline{1,m}$. 
Then, by construction, 
\begin{equation}\label{835}
\Psi_{i\alpha}=0 \in M_s, \quad i=\overline{1,p-1}, \quad \alpha=\overline{1,m}.
\end{equation}
Set 
\begin{equation}
\label{840}
x\!:=z_1.
\end{equation}
We will prove that this equality establishes a correspondence between system \eqref{820},  \eqref{830} and system \eqref{210}, \eqref{220}. By \eqref{840}, from system \eqref{820}, taking into account \eqref{835}, we obtain
\begin{equation}
\label{850}
\begin{aligned}
\dot x & = \dot{z}_1=z_2,\\
\ddot x & = \dot{z}_2=z_3,\\
\ldots & \ldots \ldots \ldots ,\\
 x^{(p-1)} & = \dot{z}_{p-1}=z_p.
\end{aligned}
\end{equation}
Next,
\begin{align}
\notag  & x^{(p)}  = \dot{z}_{p}=z_{p+1}+\sum_{\alpha=1}^m \Psi_{p\alpha}u_{\alpha},
\\
\notag  & x^{(p+1)} = \big(x^{(p)} \big)' = \dot{z}_{p+1}+ \sum_{\alpha=1}^m \Psi_{p\alpha}\dot u_{\alpha}
= {z}_{p+2}+ \sum_{\alpha=1}^m \Psi_{p+1,\alpha}{u}_{\alpha}+
\sum_{\alpha=1}^m \Psi_{p\alpha}\dot{u}_{\alpha}, \\
\label{860} & \ldots \ldots\ldots \ldots \ldots \ldots \ldots \ldots \ldots \ldots \ldots \ldots \ldots \ldots \ldots, \\
\notag  & x^{(\varkappa)} = z_{\varkappa+1}+ 
\sum_{i=p}^{\varkappa} \sum_{\alpha=1}^m  \Psi_{i\alpha}{u}_{\alpha}^{(\varkappa-i)},\\
\notag  & \ldots \ldots\ldots \ldots \ldots \ldots \ldots \ldots \ldots \ldots \ldots \ldots \ldots \ldots \ldots, \\
\notag  & x^{(n-1)} = z_{n}+  
\sum_{i=p}^{n-1} \sum_{\alpha=1}^m \Psi_{i\alpha}{u}_{\alpha}^{(n-1-i)}.
\end{align}
Finally,
\begin{equation}
\label{870}
x^{(n)} = -A_{1} z_{n} -A_{2}z_{n-1}- \ldots - A_{n} z_1+
\sum_{i=p}^{n} \sum_{\alpha=1}^m  \Psi_{i\alpha}{u}_{\alpha}^{(n-i)}.
\end{equation}
From \eqref{840} and \eqref{850}, we have
\begin{equation}
\label{880}
z_{\nu}=x^{(\nu-1)}, \quad \nu=\overline{1,p}.
\end{equation}
From \eqref{860}, we have
\begin{equation}
\label{890}
z_{\nu}=x^{(\nu-1)}-\sum_{i=p}^{\nu-1} \sum_{\alpha=1}^m  \Psi_{i\alpha}{u}_{\alpha}^{(\nu-1-i)}, \quad \nu=\overline{p+1,n}.
\end{equation}
Substituting \eqref{880} and \eqref{890} into \eqref{870}, we get
\begin{equation}
\label{900}
\begin{gathered}
 x^{(n)} +A_{1}x^{(n-1)}+A_{2}x^{(n-2)}+ \ldots +A_{n}x
\\ =
\sum_{i=p}^{n} \sum_{\alpha=1}^m  \Psi_{i\alpha}{u}_{\alpha}^{(n-i)}+
A_1 \sum_{i=p}^{n-1} \sum_{\alpha=1}^m  \Psi_{i\alpha}{u}_{\alpha}^{(n-1-i)}+\ldots+
A_{n-p} \sum_{i=p}^{p} \sum_{\alpha=1}^m  \Psi_{i\alpha}{u}_{\alpha}^{(p-i)} \\ =
\sum_{\varkappa=p}^n A_{n-\varkappa} \sum_{i=p}^{\varkappa} \sum_{\alpha=1}^m  \Psi_{i\alpha}{u}_{\alpha}^{(\varkappa-i)}.
\end{gathered}
\end{equation}
Denote the right-hand side of \eqref{900} by $\zeta$. We change the summation order in \eqref{900}, replacing $\sum\limits_{\varkappa=p}^{n}  \sum\limits_{i=p}^{\varkappa} $ by $\sum\limits_{i=p}^{n}  \sum\limits_{\varkappa=i}^{n}$. Then,
\begin{equation}
\label{910}
\zeta = \sum\limits_{i=p}^{n}  \sum\limits_{\varkappa=i}^{n}
\sum_{\alpha=1}^m    A_{n-\varkappa} \Psi_{i\alpha}{u}_{\alpha}^{(\varkappa-i)}.
\end{equation}
In \eqref{910}, we replace the summation index $\varkappa$ by $l=n-(\varkappa-i)$. Then, since $\varkappa$ ranges from $i$ to $n$, it follows that $l$ ranges  from $i$ to $n$ as well. Therefore,
\begin{equation}
\label{920}
\zeta = \sum\limits_{i=p}^{n}  \sum\limits_{l=i}^{n}
\sum_{\alpha=1}^m    A_{l-i} \Psi_{i\alpha}{u}_{\alpha}^{(n-l)}.
\end{equation}
Change the summation order in \eqref{920}, replacing $\sum\limits_{i=p}^{n}  \sum\limits_{l=i}^{n}$ by $\sum\limits_{l=p}^{n}  \sum\limits_{i=p}^{l}$. Then,
\begin{equation}
\label{930}
\zeta = 
\sum_{\alpha=1}^m   \sum\limits_{l=p}^{n}  \sum\limits_{i=p}^{l} A_{l-i} \Psi_{i\alpha}{u}_{\alpha}^{(n-l)}.
\end{equation}
From \eqref{810}, we have
\begin{equation}\label{935}
\mathcal{B}=P\Psi.
\end{equation}
Hence, taking
into account \eqref{835}, from \eqref{935}, we obtain that, for $l=\overline{p,n}$ and $\alpha=\overline{1,m}$,
\begin{equation}\label{937}
B_{l\alpha}=\sum\limits_{i=1}^{l} A_{l-i} \Psi_{i\alpha}=\sum\limits_{i=p}^{l} A_{l-i} \Psi_{i\alpha}.
\end{equation}
From \eqref{930} and \eqref{937}, we have
\begin{equation}
\label{940}
\zeta = 
\sum_{\alpha=1}^m   \sum\limits_{l=p}^{n}  B_{l \alpha}{u}_{\alpha}^{(n-l)}.
\end{equation}
We see that \eqref{940} coincides with the right-hand side of \eqref{210}. Thus, we obtain that equality \eqref{900} coincides with equation \eqref{210}.

Equation \eqref{830} is equivalent to \eqref{220} due to equalities \eqref{840}, \eqref{850}, and
$\Omega=\mathcal{C}^{\mathcal{T}}$. Q.E.D.
\hfill $\square$

From Theorem \ref{teo-7-1}, it follows that system \eqref{210}, \eqref{220} can be written in
the form of system \eqref{820}, \eqref{830} with matrices \eqref{810}, which, in its turn, is
a special case of system \eqref{080}, \eqref{082} with matrices \eqref{400}, \eqref{410}, \eqref{420}. 

Suppose that condition \eqref{790} holds for the system \eqref{210}, \eqref{220}.
In terms of system \eqref{820}, \eqref{830}, this means that the following equalities are fulfilled:
\begin{gather}\label{950}
\Gamma_i=A_i-{\rm SP}_s\,(\mathcal{J}^{i-1}P\Psi Q\Omega), \quad i=\overline{1,n}.
\end{gather}
Denote $D=\Psi Q \Omega$. Then, $D$ has the form \eqref{435} and $\Sigma=\mathcal{A}$ has the form \eqref{400}. Equalities \eqref{950} takes the form
\eqref{440}. By Lemma \ref{lemma-03}, equalities \eqref{950} are equivalent to equalities
\begin{equation*}%\label{960}
\Gamma_i=A_i-{\rm SP}_s\,(\mathcal{N}^{i-1}\Psi Q\Omega), \quad i=\overline{1,n},
\end{equation*}
where
$\mathcal{N}_0\!:=I\in M_{ns}$, $\mathcal{N}_{\nu}\!:=\mathcal{N}_{\nu}\cdot \Sigma - (I\otimes A_{\nu}) \in M_{ns}$. Condition \eqref{800} for system \eqref{210}, \eqref{220}  turns into condition \eqref{428} for system \eqref{820}, \eqref{830} with matrices \eqref{810}. Thus, Theorem \ref{teo-AMCA-LSOF-1} is a generalization of (sufficient conditions of) Theorem~3 in \cite{Zaitsev_2021}.

\section{Sufficient conditions of arbitrary scalar coefficient assignment for the characteristic polynomial by LSOF \label{sect-ASCA}}

\begin{definition}\label{def-6-01} 
We say that, for system \eqref{080}, \eqref{082}, \eqref{084}, the problem of {\it arbitrary scalar coefficient assignment (ASCA) for the characteristic polynomial (CP) by LSOF is resolvable} if for any $\delta_i\in \mathbb{K}$, $i=\overline{1,ns}$, there exists a gain matrix $Q\in M_{ms,ks}(\mathbb{K})$  such that the characteristic polynomial $\chi(F+GQH,\lambda)$ of the matrix of the closed-loop system \eqref{088} satisfies the equality
\begin{equation*}%\label{970}
\chi(F+GQH,\lambda)=\lambda^{ns}+\delta_1 \lambda^{ns-1} + \ldots +\delta_{ns-1}\lambda +\delta_{ns}.
\end{equation*}
\end{definition}

\begin{lemma}\label{lem-big-polynom}
\it
For any numbers $\delta_i\in\mathbb{K}$, $i=\overline{1,ns}$, there exist matrices $\Gamma_j\in M_s(\mathbb{K})$, $j=\overline{1,n},$ such that the  characteristic polynomial $\chi(\Phi,\lambda)$ of the matrix \eqref{064}  satisfies the equality
\begin{equation*}%\label{980}
\chi(\Phi,\lambda)=\lambda^{ns}+\delta_1 \lambda^{ns-1} + \ldots +\delta_{ns-1}\lambda +\delta_{ns}.
\end{equation*}
\end{lemma}

The proof of Lemma \ref{lem-big-polynom} is given in Theorem 5 of \cite{Zaitsev_2021}.

\begin{theorem} \label{teo-AMESA-implies-AESA}
\it
If, for system \eqref{080}, \eqref{082}, \eqref{084}, the problem of AMCA for CMP by LSOF is resolvable, then, for system \eqref{080}, \eqref{082}, \eqref{084}, the problem of ASCA for CP by LSOF is resolvable, 
\end{theorem}

Theorem \ref{teo-AMESA-implies-AESA} follows from Lemma \ref{lem-big-polynom}.

\begin{theorem}\label{teo-ACSA-LSOF-1}
\it
For system \eqref{080}, \eqref{082},  with coefficients \eqref{400}, \eqref{410}, \eqref{420}, the problem of ASCA for CP by LSOF is resolvable, if condition \eqref{428} holds.
\end{theorem}

Theorem \ref{teo-ACSA-LSOF-1} follows from Theorem \ref{teo-AMCA-LSOF-1} and Theorem \ref{teo-AMESA-implies-AESA}.

\begin{remark}
\textrm{Theorem \ref{teo-ACSA-LSOF-1} extends Theorem 7 of \cite{Zaitsev_2021} from systems \eqref{210}, \eqref{220} to systems
\eqref{080}, \eqref{082}  with coefficients \eqref{400}, \eqref{410}, \eqref{420}.}
\end{remark}

\begin{remark}\label{rem-suff-not-ness-Theta}
\textrm{Condition ${\rm rank}\,\Theta = ns^2$ in Theorem \ref{teo-ACSA-LSOF-1} is only sufficient but not necessary. 
This is confirmed by the example in Remark 7 of \cite{Zaitsev_2021} for system  \eqref{210}, \eqref{220}.}
\end{remark}

\begin{remark}
\textrm{The converse assertion to Theorem \ref{teo-AMESA-implies-AESA} is not true. This is confirmed by Example \ref{example3} given below in Section~\ref{sec-examples}.}
\end{remark}

\section{Special cases \label{sect-special-cases}}

Consider system \eqref{080}, \eqref{082} with matrices \eqref{400}, \eqref{410}, \eqref{420}.

\subsection{Blocks of $H$ are scalar matrices}

Suppose that the blocks of the matrix $H$ are scalar matrices, i.e.,
\begin{equation}\label{1000}
\begin{gathered}
H=\{H_{\beta i}\}\in M_{ks,ns}, \quad
H_{\beta i}=h_{\beta i}I, \quad h_{\beta i}\in\mathbb{K}, \quad I\in M_s, \quad \beta=\overline{1,k}, \quad i=\overline{1,n}, \\
h_{\beta i}=0, \quad \beta=\overline{1,k}, \quad i=\overline{p+1,n}.
\end{gathered}
\end{equation}
Then $H^{\mathcal{T}}=H^T$. Hence, 
\begin{equation}\label{1010}
(H^{\mathcal{T}})^T=H.
\end{equation}
By applying \cite[Lemma 2]{Zaitsev_2021} to $X=F^{i-1}GQ$, $Y=H$, we obtain
\begin{equation}\label{1020}
{\rm
SP}_s(F^{i-1}GQH)=
{\rm SP}_s\,(HF^{i-1}GQ).
\end{equation}
By \eqref{1020}, equalities \eqref{690} take the form
\begin{equation}\label{1025}
T_i={\rm SP}_s\,(HF^{i-1}GQ), \quad i=\overline{1,n}.
\end{equation}
Taking into account \eqref{1010}, equalities  \eqref{1025} are equivalent to 
\begin{equation}\label{1030}
T_i={\rm SP}_s\,\big((H^{\mathcal{T}})^T F^{i-1}GQ\big), \quad i=\overline{1,n}.
\end{equation}
By \cite[Eq. (33), Lemma 6]{Zaitsev_2021},
\begin{equation*}
{\rm SP}_s \big((H^{\mathcal{T}})^T F^{i-1}GQ\big)=
{\rm VecRR}_s \big((H^{\mathcal{T}})^T F^{i-1}G\big)\cdot {\rm VecCC}_s Q. 
\end{equation*}
Thus, one can rewrite system \eqref{1030} with respect to coefficients of $Q$ in the following form:
\begin{equation}\label{1040}
\Omega\cdot V= W.
\end{equation}
Here
\begin{gather}\label{1050}
\Omega:=\left[\begin{matrix}
{\rm VecRR}_{s}\big((H^{\mathcal{T}})^T G\big)\\
{\rm VecRR}_{s}\big((H^{\mathcal{T}})^T FG\big)\\
\ldots\ldots\ldots\ldots\ldots\ldots \ldots \ldots \\
{\rm VecRR}_{s}\big((H^{\mathcal{T}})^T F^{n-1}G\big)
\end{matrix}\right]\in M_{ns,kms},\\
\notag
 W:= P^{-1}(\widehat{\Gamma}-\widehat{A})\in M_{ns,s},\\
\notag
V:={\rm VecCC}_{s}\,Q \in M_{kms,s}.
\end{gather}
System \eqref{1040} is solvable with respect to $V$ for any $\Gamma_i\in M_s(\mathbb{K})$, $i=\overline{1,n}$, iff
\begin{equation}\label{1070}
{\rm rank}\, \Omega=ns.
\end{equation}
In particular,  system \eqref{1040} has the solution 
\begin{equation*}%\label{1080}
V=\Omega^T(\Omega\Omega^T)^{-1}W.
\end{equation*}
The required
matrix $Q$ can be found from the equality  
\begin{equation*}%\label{1090}
Q={\rm VecCC}_s^{-1} V.
\end{equation*}
 Thus, the following theorem holds.

\begin{theorem}\label{teo-H-scalar}
\it
Suppose that, for system \eqref{080}, \eqref{082}  with coefficients \eqref{400}, \eqref{410}, \eqref{420}, 
the blocks of the matrix $H$ are scalar matrices. If condition \eqref{1070} holds, then, for this system, the problem of AMCA for CMP by LSOF is resolvable.
\end{theorem}

\begin{corollary}\label{cor-H-scalar}
\it
Suppose that, for system \eqref{080}, \eqref{082}  with coefficients \eqref{400}, \eqref{410}, \eqref{420}, 
the blocks of the matrix $H$ are scalar matrices. If condition \eqref{1070} holds, then, for this system, the problem of ASCA for CP by LSOF is resolvable.
\end{corollary}

Corollary \ref{cor-H-scalar} follows from Theorem \ref{teo-H-scalar} and Theorem \ref{teo-AMESA-implies-AESA}.

\begin{remark}
\textrm{Theorem \ref{teo-H-scalar} and Corollary \ref{cor-H-scalar}  extend sufficient conditions of Theorem 8 and Corollary~1 of \cite{Zaitsev_2021}, respectively,  from systems \eqref{210}, \eqref{220} to systems
\eqref{080}, \eqref{082}  with coefficients \eqref{400}, \eqref{410}, \eqref{420}.}
\end{remark}

\begin{remark}
\textrm{Condition \eqref{1070} in Corollary \ref{cor-H-scalar} is only sufficient but not necessary. This is confirmed by the
example in Remark 9 of \cite{Zaitsev_2021}, for system \eqref{210}, \eqref{220}.
The question of the necessity of condition \eqref{1070} in Theorem \ref{teo-H-scalar} is still open.}
\end{remark}

Suppose that system \eqref{080}, \eqref{082}  with coefficients \eqref{400}, \eqref{410}, \eqref{420} is a system with LSSF. From this, it follows that $k=n$ and $H=I\in M_{ns}$. Hence, $p=n$. This means that, in the matrix $G$ of \eqref{410}, the first $n-1$ block rows are zero, and the last block row is $G_n=[G_{n1} \, \cdots \, G_{nm}]\in M_{s,ms}$. Blocks of $H$ are scalar matrices. Construct the matrix $\Omega$ of \eqref{1050}. Due to the form of the matrices $F$, $G$, and $H$, we obtain that $\Omega$ has the form
\begin{equation*}
\Omega=\left[
\begin{matrix}
0 & \dots & 0 & G_n \\
0 & \dots & G_n & * \\
\vdots & \iddots & \vdots & \vdots \\
G_n & \dots & * & *
\end{matrix}\right]\in M_{ns,nms}.
\end{equation*}
Clearly, condition \eqref{1070} is equaivalent to
\begin{equation}\label{1090}
{\rm rank}\, G_n=n.
\end{equation}
So, we have the following corollary from Theorem \ref{teo-H-scalar}.

\begin{corollary}\label{cor-H=I}
\it
Suppose that, for system \eqref{080}, \eqref{082}  with coefficients \eqref{400}, \eqref{410}, \eqref{420}, the following conditions hold: $k=n$, $H=I\in M_{ns}$, and $p=n$.
If condition \eqref{1090} holds, then, for this system, the problem of AMCA for CMP by LSSF is resolvable.
\end{corollary}

The question of the necessity in Corollary \ref{cor-H=I} requires additional research.

\subsection{Blocks of $F$ and $G$ are scalar matrices}

Suppose that the blocks of the matrices $F$ and $G$ are scalar matrices, i.e.,
\begin{gather}\label{1100}
A_i=a_i I, \quad a_i\in\mathbb{K},\quad i=\overline{1,n},\quad I\in M_s, \\
\label{1110}
\begin{gathered}
G=\{G_{j\alpha}\}\in M_{ns,ms}, \quad
G_{j\alpha}=g_{j\alpha}I, \quad g_{j\alpha}\in\mathbb{K}, \quad I\in M_s, \quad j=\overline{1,n}, \quad \alpha=\overline{1,m},\\
g_{j\alpha}=0, \quad j=\overline{1,p-1}, \quad \alpha=\overline{1,m}.
\end{gathered}
\end{gather}
Then, for any $i=\overline{1,n}$, the blocks of the matrices $F^{i-1}G$ are scalar matrix as well. Hence,
\begin{equation}\label{1120}
(F^{i-1}G)^{\mathcal{T}}=(F^{i-1}G)^{{T}}.
\end{equation}
By applying \cite[Lemma 2]{Zaitsev_2021} to $X=QH$, $Y=F^{i-1}G$, we obtain
\begin{equation}\label{1130}
{\rm
SP}_s(F^{i-1}GQH)=
{\rm SP}_s(QHF^{i-1}G).
\end{equation}
By \eqref{1130}, equalities \eqref{690} take the form
\begin{equation}\label{1135}
T_i={\rm SP}_s\,(QHF^{i-1}G), \quad i=\overline{1,n}.
\end{equation}
By \cite[Eq. (34), Lemma 6]{Zaitsev_2021},
\begin{equation*}
{\rm SP}_s(QHF^{i-1}G)=
{\rm VecCR}_s Q\cdot {\rm VecRC}_s (HF^{i-1}G). 
\end{equation*}
Thus, one can rewrite system \eqref{1135} with respect to coefficients of $Q$ in the following form:
\begin{equation}\label{1140}
X\cdot \Xi =Y.
\end{equation}
Here
\begin{gather*}
\Xi:=\left[{\rm VecRC}_{s}(HG), \ldots, {\rm VecRC}_{s}(H F^{n-1}G)
\right]\in M_{mks,ns},\\
Y:= \big( P^{-1}(\widehat{\Gamma}-\widehat{A})\big)^{\mathcal{T}}\in M_{s,ns},\\
X:={\rm VecCR}_{s}\,Q \in M_{s,mks}.
\end{gather*}

System \eqref{1140} is solvable with respect to $X$ for any $\Gamma_i\in M_s(\mathbb{K})$, $i=\overline{1,n}$, iff ${\rm rank}\, \Xi=ns$. In particular, system \eqref{1140} has the solution 
\begin{equation*}%\label{1170}
X=Y(\Xi^T\Xi)^{-1}\Xi^T.
\end{equation*}
  The required
matrix $Q$ can be found from the equality 
\begin{equation*} %\label{1180}
Q = {\rm VecCR}_s^{-1}X.
\end{equation*}

Let us rewrite system \eqref{1140} in the form
\begin{equation*}
\Xi^T\cdot X^T =Y^T.
\end{equation*}
Consider the matrix
\begin{equation}\label{1190}
\Xi^T=\left[
\begin{matrix}
\big[{\rm VecRC}_{s}(HG)\big]^T\\
\big[{\rm VecRC}_{s}(HFG)\big]^T\\
 \ldots \ldots \ldots \ldots \ldots \ldots \ldots \\
\big[{\rm VecRC}_{s}(HF^{n-1} G)\big]^T
\end{matrix}
\right]\in M_{ns,mks}.
\end{equation}
For any $i=\overline{1,n}$, by \eqref{336-5}, we have
\begin{equation}\label{1200}
\big[{\rm VecRC}_s (HF^{i-1}G)\big]^T = 
{\rm VecCR}_s \big( (HF^{i-1}G)^T \big).
\end{equation}
By \eqref{336-3},
\begin{equation}\label{1210}
{\rm VecCR}_s \big( (HF^{i-1}G)^T \big)=
{\rm VecRR}_s \Big(\big( (HF^{i-1}G)^T \big)^{\mathcal{T}}\Big).
\end{equation}
By \eqref{1120},
\begin{equation}\label{1220}
(HF^{i-1}G)^T=(F^{i-1}G)^T H^T = (F^{i-1}G)^{\mathcal{T}} H^T.
\end{equation}
By \cite[Eq. (30), Lemma 5]{Zaitsev_2021},
\begin{equation}\label{1230}
(F^{i-1}G)^{\mathcal{T}} H^T= \big( (H^T)^\mathcal{T} F^{i-1} G\big)^{\mathcal{T}}.
\end{equation}
It follows from \eqref{1220}, \eqref{1230}, and assertion 1 of \cite[Lemma 5]{Zaitsev_2021} that
\begin{equation}\label{1240}
\big( (HF^{i-1}G)^T \big)^{\mathcal{T}}=
\big( (F^{i-1}G)^{\mathcal{T}} H^T \big)^{\mathcal{T}}=
 \Big(\big( (H^T)^\mathcal{T} F^{i-1} G\big)^{\mathcal{T}}\Big)^{\mathcal{T}}=
 (H^T)^\mathcal{T} F^{i-1} G.
\end{equation}
It follows from \eqref{1200}, \eqref{1210}, and \eqref{1240} that the matrix \eqref{1190} has the form
\begin{equation}\label{1250}
\Xi^T=\left[
\begin{matrix}
{\rm VecRR}_{s}\big((H^T)^\mathcal{T} G \big)\\
{\rm VecRR}_{s}\big((H^T)^\mathcal{T}  F G \big)\\
 \ldots \ldots \ldots \ldots \ldots \ldots \ldots \ldots \\
{\rm VecRR}_{s}\big((H^T)^\mathcal{T} F^{n-1} G \big)
\end{matrix}
\right]\in M_{ns,mks}.
\end{equation}
From \eqref{1250}, it follows that the matrix $\Xi^T$  coincides with \eqref{1050}.  Thus, the following theorem holds.

\begin{theorem}\label{teo-FG-scalar}
\it
Suppose that, for system \eqref{080}, \eqref{082}  with coefficients \eqref{400}, \eqref{410}, \eqref{420}, 
the blocks of the matrices $F$ and $G$ are scalar matrices. If condition \eqref{1070} holds, then, for this system, the problem of AMCA for CMP by LSOF is resolvable.
\end{theorem}

\begin{corollary}\label{cor-FG-scalar}
\it
Suppose that, for system \eqref{080}, \eqref{082}  with coefficients \eqref{400}, \eqref{410}, \eqref{420}, 
the blocks of the matrices $F$ and $G$  are scalar matrices. If condition \eqref{1070} holds, then, for this system, the problem of ASCA for CP by LSOF is resolvable.
\end{corollary}

Corollary \ref{cor-FG-scalar} follows from Theorem \ref{teo-FG-scalar} and Theorem \ref{teo-AMESA-implies-AESA}.

\begin{remark}
\textrm{Theorem \ref{teo-FG-scalar} and Corollary \ref{cor-FG-scalar}  extend sufficient conditions of Theorem 9 and Corollary~2 of \cite{Zaitsev_2021}, respectively,  from systems \eqref{210}, \eqref{220} to systems
\eqref{080}, \eqref{082}  with coefficients \eqref{400}, \eqref{410}, \eqref{420}.}
\end{remark}

\begin{remark}
\textrm{Condition \eqref{1070} in Corollary \ref{cor-FG-scalar} is only sufficient but not necessary. This is confirmed by the
example in Remark 11 of \cite{Zaitsev_2021} for system \eqref{210}, \eqref{220}.
The question of the necessity of condition \eqref{1070} in Theorem \ref{teo-H-scalar} is still open.}
\end{remark}

\subsection{Blocks of $F$, $G$ and $H$ are scalar matrices}

Suppose that the blocks of the matrices $F$, $G$ and $H$ are scalar matrices, i.e., for the matrices \eqref{400}, \eqref{410}, \eqref{420}  the relations \eqref{1100}, \eqref{1110}, \eqref{1000} hold. Denote 
\begin{gather*}
F_0=J+e_n\psi, \quad \psi=[-a_n,\ldots,-a_1]\in M_{1,n},\\
G_0=\{g_{i\alpha}\}\in M_{n,m}, \quad H_0=\{h_{\beta i}\}\in M_{k,n}, 
\quad i=\overline{1,n}, \quad
\alpha=\overline{1,m}, \quad \beta=\overline{1,k}.
\end{gather*}
Then,
\begin{equation*}%\label{1260}
F=F_0\otimes I, \quad G=G_0\otimes I , \quad H=H_0\otimes I, \quad I\in M_s.
\end{equation*}

By using properties of the Kronecker product, for all $i=\overline{1,n}$, we have
\begin{gather*}
{\rm VecRR}_s \big((H^{\mathcal{T}})^T F^{i-1} G\big)=
{\rm VecRR}_s \big(H F^{i-1} G\big)=
\big((H_0\otimes I) (F_0^{i-1}\otimes I)(G_0\otimes I)\big)=\\
{\rm VecRR}_s \big((H_0 F_0^{i-1} G_0)\otimes I\big)=\Big({\rm vecr}\,(H_0 F_0^{i-1} G_0)\Big)\otimes I.
\end{gather*}
Thus, the matrix \eqref{1050} has the form
\begin{equation*}
\Omega =
\left[\begin{matrix}
\Big({\rm vecr}\,(H_0 G_0)\Big)\otimes I\\
\Big({\rm vecr}\,(H_0 F_0 G_0)\Big)\otimes I\\
\ldots\ldots\ldots\ldots\ldots\ldots\ldots \ldots \\
\Big({\rm vecr}\,(H_0 F_0^{n-1} G_0)\Big)\otimes I
\end{matrix}\right]=\Omega_0\otimes I,
\end{equation*}
where 
\begin{equation*}
\Omega_0 =
\left[\begin{matrix}
{\rm vecr}\,(H_0 G_0)\\
{\rm vecr}\,(H_0 F_0 G_0)\\
\ldots\ldots\ldots\ldots\ldots \\
{\rm vecr}\,(H_0 F_0^{n-1} G_0)
\end{matrix}\right]\in M_{n,km}.
\end{equation*}
Since ${\rm rank}\, (X\otimes Y)={\rm rank}\,X\cdot {\rm
rank}\,Y$, we obtain that ${\rm rank}\,\Omega=ns$ iff ${\rm rank}\,\Omega_0=n$. This condition is equivalent to linear independence of matrices
\begin{equation} \label{1270}
H_0 G_0, \quad H_0 F_0 G_0, \quad \ldots, \quad H_0 F_0^{n-1} G_0.
\end{equation}
Thus, from Theorem \ref{teo-H-scalar}, it follows the theorem.

\begin{theorem}\label{teo-FGH-scalar}
\it
Suppose that, for system \eqref{080}, \eqref{082}  with coefficients \eqref{400}, \eqref{410}, \eqref{420}, 
the blocks of the matrices $F$, $G$ and $H$ are scalar matrices. If the matrices \eqref{1270} are linearly independent, then, for this system, the problem of AMCA for CMP by LSOF is resolvable.
\end{theorem}

\begin{corollary}\label{cor-FGH-scalar}
\it
Suppose that, for system \eqref{080}, \eqref{082}  with coefficients \eqref{400}, \eqref{410}, \eqref{420}, 
the blocks of the matrices $F$ and $G$  are scalar matrices. If the matrices \eqref{1270} are linearly independent, then, for this system, the problem of ASCA for CP by LSOF is resolvable.
\end{corollary}

Corollary \ref{cor-FGH-scalar} follows from Theorem \ref{teo-FGH-scalar} and Theorem 
\ref{teo-AMESA-implies-AESA}.

\begin{remark}
\textrm{Theorem \ref{teo-FGH-scalar} and Corollary \ref{cor-FGH-scalar}  extend sufficient conditions of Theorem 10 and Corollary 3 of \cite{Zaitsev_2021}, respectively,  from systems \eqref{210}, \eqref{220} to systems
\eqref{080}, \eqref{082}  with coefficients \eqref{400}, \eqref{410}, \eqref{420}.}
\end{remark}

\begin{remark} \textrm{The converse assertion to Corollary \ref{cor-FGH-scalar} is not true, in general case. This is confirmed by the example in Remark 12 of \cite{Zaitsev_2021} for system \eqref{210}, \eqref{220} for the case $\mathbb{K}=\mathbb{C}$.  The converse assertion to Theorem \ref{cor-FGH-scalar} is true if $s=1$ (see Theorem \ref{theo-0-2}).  
We hypothesize that this statement is true for any $s\ge 1$ but it has not yet been proven.}
\end{remark}

\section{Sufficient conditions to solving the problem of AMCA for CMP by LSOF for systems with a lower block Hessenberg matrix \label{sect-Hessenberg}}

Consider the system of the form \eqref{080}, \eqref{082}:
\begin{gather}\label{1500}
\dot{\widetilde x} = \widetilde F \widetilde x + \widetilde G \widetilde u, \quad \widetilde x \in\mathbb{K}^{ns}, \quad
\widetilde u \in\mathbb{K}^{ms}, \\
\label{1510}
\widetilde y= \widetilde H \widetilde x, \quad \widetilde y \in\mathbb{K}^{ks}.
\end{gather}
Here we suppose that the coefficients of this system have the following special form: for some $p\in\{1,\ldots,n\}$, the first $p-1$ block rows of the matrix $\widetilde G$ are zero, the last $n-p$ block columns of the matrix $\widetilde H$ are zero, the matrix $\widetilde F$ is a lower block Hessenberg matrix, i.e.,
\begin{gather}\label{1520}
\widetilde F=
\left[\begin{matrix}
\widetilde{F}_{11} & \widetilde{F}_{12} & 0 & \ldots & 0 \\
\widetilde{F}_{21} & \widetilde{F}_{22} & \widetilde{F}_{23} & \ldots & 0 \\
\vdots & \vdots & \vdots & \ddots & \vdots  \\
\widetilde{F}_{n-1,1} & \widetilde{F}_{n-1,2} & \widetilde{F}_{n-1,3} & \ldots & \widetilde{F}_{n-1,n} \\
\widetilde{F}_{n1} & \widetilde{F}_{n2} & \widetilde{F}_{n3} & \ldots & \widetilde{F}_{nn}
\end{matrix} \right], \qquad 0,\widetilde{F}_{ij}\in M_s; \quad \det\widetilde{F}_{i,i+1}\ne 0, \quad i=\overline{1,n-1};\\
\label{1530}
\widetilde{G}=\begin{bmatrix} 0 & \ldots & 0\\
\vdots & & \vdots\\
0 & \ldots & 0\\
\widetilde{G}_{p1}& \ldots & \widetilde{G}_{pm}\\
\vdots & & \vdots\\
\widetilde{G}_{n1} & \ldots & \widetilde{G}_{nm}
\end{bmatrix}, \qquad 0,\widetilde{G}_{j\alpha}\in M_s, \quad j=\overline{p,n}, \quad \alpha=\overline{1,m}, \\
\label{1540}
\widetilde{H}=\begin{bmatrix} \widetilde{H}_{11} & \ldots & \widetilde{H}_{1p} & 0 & \ldots & 0\\
\vdots & & \vdots & \vdots & & \vdots\\
\widetilde{H}_{k1} & \ldots & \widetilde{H}_{kp} & 0 & \ldots & 0
\end{bmatrix}, \qquad 0,\widetilde{H}_{\beta i}\in M_s, \quad \beta=\overline{1,k}, \quad i=\overline{1,p}.
\end{gather}
By using Lemma \ref{lemma-00}, from the matrix $\widetilde F$, construct the non-degenerate lower block triangular matrix $\widetilde S$ such that $\widetilde S \widetilde F \widetilde S^{-1}$ is a lower block Frobenius matrix. Let us make the replacement of variables in system \eqref{1500}, \eqref{1510}:
\begin{equation*}% \label{1550}
x=\widetilde S\widetilde x, \quad u=\widetilde u, \quad y=\widetilde y.
\end{equation*}
Then system \eqref{1500}, \eqref{1510} is transformed into system  \eqref{080}, \eqref{082}, the matrices of which are related to the matrices of system \eqref{1500}, \eqref{1510}  by the equalities:
\begin{equation} \label{1560}
F=\widetilde S \widetilde F \widetilde S^{-1}, \quad G= \widetilde S \widetilde G, \quad H= \widetilde H \widetilde S^{-1}.
\end{equation}
By Lemma \ref{lemma-00}, the matrix $F$ has the form \eqref{400}. Since the matrix $\widetilde S$ is lower block triangular, the matrices $G$ and $H$ have the form \eqref{410} and \eqref{420} respectively. Construct the matrix
\begin{equation}\label{1570}
\widetilde\Theta=\left[\begin{matrix}
{\rm VecRR}_{s^2}\big(((\widetilde{H}\widetilde S^{-1})^{\mathcal{T}})^T\star (\widetilde S\widetilde{G})\big)\\
{\rm VecRR}_{s^2}\big(((\widetilde{H}\widetilde S^{-1})^{\mathcal{T}})^T\star (\widetilde S\widetilde{F}\widetilde{G})\big)\\
\ldots\ldots\ldots\ldots\ldots\ldots\ldots \ldots\ldots\ldots \ldots \ldots  \\
{\rm VecRR}_{s^2}\big(((\widetilde{H}\widetilde S^{-1})^{\mathcal{T}})^T\star (\widetilde S\widetilde{F}^{n-1}\widetilde{G})\big)
\end{matrix}\right]\in M_{ns^2,kms^2}.
\end{equation}
Then, due to \eqref{1560}, we obtain that $\widetilde\Theta=\Theta$.

\begin{theorem}\label{teo-Hessenberg}
\it
For system \eqref{1500}, \eqref{1510}, with coefficients \eqref{1520}, \eqref{1530}, \eqref{1540}, the problem of AMCA for CMP by LSOF is resolvable, if the following condition holds:
\begin{equation}
\label{1580}
{\rm rank}\,\widetilde\Theta = ns^2.
\end{equation}
\end{theorem}

\emph{Proof.}
System \eqref{1500}, \eqref{1510} closed-loop by LSOF $\widetilde u=Q\widetilde y$ has the form
\begin{equation}\label{1585}
\dot{\widetilde x}=\big(\widetilde F + \widetilde G  Q \widetilde H\big) \widetilde x.
\end{equation}
By \eqref{1560}, we have $F+GQH=\widetilde S \big(\widetilde F + \widetilde G Q \widetilde H\big)  \widetilde S^{-1}$, i.e. the matrices of systems \eqref{1585} and \eqref{088} are similar. Hence, the problem of AMCA for CMP by LSOF is resolvable for system \eqref{1500}, \eqref{1510} with coefficients \eqref{1520}, \eqref{1530}, \eqref{1540} iff this problem is resolvable for system \eqref{080}, \eqref{082},  with coefficients \eqref{400}, \eqref{410}, \eqref{420}. The latter takes place if \eqref{428} is fulfilled. Due to $\widetilde\Theta=\Theta$, we obtain the required. Q.E.D.
\hfill $\square$

Consider a partial case. Suppose that the blocks of the matrix $\widetilde S$ are scalar matrices. 
Then, the blocks of $\widetilde S^{-1}$ are also. Therefore, $(\widetilde S^{-1})^{\mathcal{T}}=(\widetilde S^{-1})^{{T}}$. Hence, by \cite[Lemma 5]{Zaitsev_2021},
\begin{equation} \label{1600}
((\widetilde{H}\widetilde S^{-1})^{\mathcal{T}})^T = 
\big((\widetilde S^{-1})^{\mathcal{T}}(\widetilde{H}^{\mathcal{T}})\big)^T = 
\big((\widetilde S^{-1})^{{T}}(\widetilde{H}^{\mathcal{T}})\big)^T = 
(\widetilde{H}^{\mathcal{T}})^T \widetilde S^{-1}.
\end{equation}
 Now, by applying \eqref{1600} and Lemma \ref{lem-F-star}, we obtain that, for any $i=\overline{1,n}$,
\begin{equation} \label{1610}
((\widetilde{H}\widetilde S^{-1})^{\mathcal{T}})^T\star (\widetilde S\widetilde{F}^{i-1}\widetilde{G}) = 
(\widetilde{H}^{\mathcal{T}})^T \widetilde S^{-1} \star (\widetilde S\widetilde{F}^{i-1}\widetilde{G}) = 
(\widetilde{H}^{\mathcal{T}})^T\star \widetilde{F}^{i-1}\widetilde{G}.
\end{equation}
From equalities \eqref{1610}, it follows that the matrix $\widetilde \Theta$ coincides with the matrix
\begin{equation}\label{1620}
\widehat\Theta=\left[\begin{matrix}
{\rm VecRR}_{s^2}\big((\widetilde{H}^{\mathcal{T}})^T\star \widetilde{G}\big)\\
{\rm VecRR}_{s^2}\big((\widetilde{H}^{\mathcal{T}})^T\star \widetilde{F}\widetilde{G}\big)\\
\ldots\ldots\ldots\ldots\ldots\ldots\ldots \ldots\ldots \\
{\rm VecRR}_{s^2}\big((\widetilde{H}^{\mathcal{T}})^T\star \widetilde{F}^{n-1}\widetilde{G}\big)
\end{matrix}\right]\in M_{ns^2,kms^2}.
\end{equation}
So, we have the following statement.

\begin{lemma} \label{lem-Hess-S}
\it
Suppose that the blocks of the matrix $\widetilde S$ are scalar matrices. Then, $\widetilde \Theta = \widehat\Theta$. In particular, condition \eqref{1580} is equivalent to the following condition:
\begin{equation}\label{1630}
{\rm rank}\,\widehat\Theta = ns^2.
\end{equation}
\end{lemma}

\begin{remark}\label{rem-scalar-S}
\textrm{A necessary and sufficient condition for the blocks of the matrix $\widetilde S$ to be scalar matrices is the condition that all blocks in the first $(n-1)$ block rows of the matrix $\widetilde F$ are scalar matrices. The sufficiency is evident. The proof of necessity is easy to obtain using induction.}
\end{remark}

So, from Remark \ref{rem-scalar-S}, Lemma \ref{lem-Hess-S}, and Theorem \ref{teo-Hessenberg}, we get the following theorem.

\begin{theorem}\label{teo-Hessenberg-2}
\it
Suppose that all blocks in the first $(n-1)$ block rows of the matrix \eqref{1520} are scalar matrices.
For system \eqref{1500}, \eqref{1510}, with coefficients \eqref{1520}, \eqref{1530}, \eqref{1540}, the problem of AMCA for CMP by LSOF is resolvable, if condition  \eqref{1630} holds.
\end{theorem}

\begin{remark}
\textrm{One can see that Theorem \ref{teo-AMCA-LSOF-1} is a partial case of Theorem \ref{teo-Hessenberg-2}.}
\end{remark}

\begin{remark}
\textrm{If not all blocks of the matrix $\widetilde S$ are scalar matrices, then the conclusion of Lemma~\ref{lem-Hess-S}
is not true, in general case; see Example \ref{example2} in Section \ref{sec-examples}. This emphasizes the difference between the cases $s=1$ and $s>1$.}
\end{remark}

Based on the proofs of Theorem \ref{teo-Hessenberg} and Theorem \ref{teo-AMCA-LSOF-1}, we present an algorithm for solving the problem of AMCA for CMP by LSOF for systems with a lower block Hessenberg matrix.

{\bf Algorithm 2.} Let the system
\eqref{1500}, \eqref{1510} with coefficients \eqref{1520}, \eqref{1530}, \eqref{1540} be given.

1. By using Lemma \ref{lemma-00}, construct the non-degenerate lower block triangular matrix $\widetilde S$ such that the replacement \eqref{1560} leads the system \eqref{1500}, \eqref{1510} with coefficients \eqref{1520}, \eqref{1530}, \eqref{1540} to a system \eqref{080}, \eqref{082} with coefficients \eqref{400}, \eqref{410}, \eqref{420}.

2. Construct the matrix $\widetilde\Theta$ of \eqref{1570}.

3. Check the condition \eqref{1580}. If this condition is satisfied, then the problem is solvable.

4. Let arbitrary matrices $\Gamma_1,\ldots,\Gamma_n\in M_s$ be given.

5. For the transformed system \eqref{080}, \eqref{082} with coefficients \eqref{400}, \eqref{410}, \eqref{420} with a lower block Frobenius matrix, apply Algorithm 1 (steps 4--10) and construct the gain matrix $Q$ and  the matrix $S$ such that $S(F+GQH)S^{-1}=\Phi$.

6. From the last equality and equalities \eqref{1560}, it follows that
$S\widetilde{S}(\widetilde{F}+\widetilde{G}Q\widetilde{H}){\widetilde{S}}^{-1}S^{-1}=\Phi$. Set $\mathcal{R}:=S\widetilde{S}$. Then,
$\mathcal{R}(\widetilde{F}+\widetilde{G}Q\widetilde{H}){\mathcal{R}}^{-1}=\Phi$.

Below we consider Example \ref{exam-Hess}, which illustrates the operation of this algorithm in relation to Theorem \ref{teo-Hessenberg-2}, which is a special case of Theorem \ref{teo-Hessenberg}. 
Under the conditions of Theorem \ref{teo-Hessenberg-2}, the equality $\widetilde \Theta = \widehat\Theta$ is satisfied. Therefore, in Example \ref{exam-Hess}, the matrix $\widehat\Theta$ is constructed instead of the matrix $\widetilde \Theta$,  and  condition \eqref{1630} is checked instead of condition \eqref{1580}.

\section{Corollaries for systems with a lower block Hessenberg matrix \label{Hessenberg-corollaries}}

Let us deduce corollaries from the results of Sections \ref{sect-Hessenberg} and \ref{sect-special-cases}.
Consider system \eqref{1500}, \eqref{1510}, with coefficients \eqref{1520}, \eqref{1530}, \eqref{1540}. By analogy with \eqref{1050}, we construct the matrix
\begin{equation*}%\label{1640}
\widetilde\Omega:=\left[\begin{matrix}
{\rm VecRR}_{s}\big((\widetilde H^{\mathcal{T}})^T \widetilde G\big)\\
{\rm VecRR}_{s}\big((\widetilde H^{\mathcal{T}})^T \widetilde{F} \widetilde{G}\big)\\
\ldots\ldots\ldots\ldots\ldots\ldots \ldots \ldots \\
{\rm VecRR}_{s}\big((\widetilde H^{\mathcal{T}})^T \widetilde{F}^{n-1}\widetilde{G}\big)
\end{matrix}\right]\in M_{ns,kms}.
\end{equation*}

\begin{corollary}\label{teo-Hessenberg-H-scalar}
\it
Suppose that, for system \eqref{1500}, \eqref{1510}, with coefficients \eqref{1520}, \eqref{1530}, \eqref{1540}, 
all blocks in the first $(n-1)$ block rows of the matrix \eqref{1520} are scalar matrices, and the blocks of the matrix \eqref{1540} are scalar matrices. Suppose that the following condition holds:
\begin{equation}\label{1650}
{\rm rank}\, \widetilde\Omega=ns.
\end{equation}
Then, for this system: (a) the problem of AMCA for CMP by LSOF is resolvable; (b) the problem of ASCA for CP by LSOF is resolvable.
\end{corollary}

Corollary \ref{teo-Hessenberg-H-scalar} follows from Theorem \ref{teo-Hessenberg-2} and Theorem \ref{teo-H-scalar} and Corollary \ref{cor-H-scalar}. In its turn, Theorem \ref{teo-H-scalar} and Corollary \ref{cor-H-scalar} are special cases of Corollary \ref{teo-Hessenberg-H-scalar} $(a)$ and $(b)$, respectively.

\begin{corollary}\label{teo-Hessenberg-FG-scalar}
\it
Suppose that, for system \eqref{1500}, \eqref{1510}, with coefficients \eqref{1520}, \eqref{1530}, \eqref{1540}, 
the blocks of the matrices \eqref{1520} and \eqref{1530} are scalar matrices.
If condition \eqref{1650} holds, then, for this system: (a) the problem of AMCA for CMP by LSOF is resolvable; (b) the problem of ASCA for CP by LSOF is resolvable.
\end{corollary}

Corollary \ref{teo-Hessenberg-FG-scalar} follows from Theorem \ref{teo-Hessenberg-2} and Theorem \ref{teo-FG-scalar} and Corollary \ref{cor-FG-scalar}. In its turn, Theorem \ref{teo-FG-scalar} and Corollary \ref{cor-FG-scalar} are special cases of Corollary \ref{teo-Hessenberg-FG-scalar} $(a)$ and $(b)$, respectively.

Now, let us suppose that the blocks of the matrices $\widetilde F$, $\widetilde G$, $\widetilde H$ are scalar matrices, i.e., matrices \eqref{1520}, \eqref{1530}, \eqref{1540} have the form
\begin{equation*}%\label{1660}
\widetilde F=\widetilde F_0\otimes I, \quad \widetilde G=\widetilde G_0\otimes I , \quad \widetilde H=\widetilde H_0\otimes I, \quad I\in M_s.
\end{equation*}
We construct the matrices
\begin{equation} \label{1670}
\widetilde H_0 \widetilde G_0, \quad \widetilde H_0 \widetilde F_0 \widetilde G_0, \quad \ldots, \quad \widetilde H_0 \widetilde F_0^{n-1} \widetilde G_0.
\end{equation}

\begin{corollary}\label{teo-Hessenberg-FGH-scalar}
\it
Suppose that, for system \eqref{1500}, \eqref{1510}, with coefficients \eqref{1520}, \eqref{1530}, \eqref{1540}, 
the blocks of the matrices $\widetilde F$, $\widetilde G$ and $\widetilde H$  are scalar matrices.
 If the matrices \eqref{1670} are linearly independent, then, for this system: (a) the problem of AMCA for CMP by LSOF is resolvable; (b) the problem of ASCA for CP by LSOF is resolvable.
\end{corollary}

Corollary \ref{teo-Hessenberg-FGH-scalar} follows from Theorem \ref{teo-Hessenberg-2} and Theorem \ref{teo-FGH-scalar} and Corollary \ref{cor-FGH-scalar}. In its turn, Theorem \ref{teo-FGH-scalar} and Corollary \ref{cor-FGH-scalar} are special cases of Corollary \ref{teo-Hessenberg-FGH-scalar} $(a)$ and $(b)$, respectively.

\begin{remark} \textrm{The converse assertion to Corollary \ref{teo-Hessenberg-FGH-scalar} $(b)$ is not true, in general case. This is confirmed by the example in Remark 12 of \cite{Zaitsev_2021} for system \eqref{210}, \eqref{220} for the case $\mathbb{K}=\mathbb{C}$.  The converse assertion to Corollary \ref{teo-Hessenberg-FGH-scalar} $(a)$ is true if $s=1$ (see Theorem \ref{theo-0-2}).  
We hypothesize that this statement is true for any $s\ge 1$ but it has not yet been proven.}
\end{remark}

\section{Block pole assignment by linear static output feedback \label{sect-assignment}}

In this section we develop the results of \cite{Leang_Shieh_1983} on block pole assignment from systems
\eqref{066}, \eqref{067} (with $m=1$) with LSSF \eqref{068}  to systems \eqref{080}, \eqref{082}, \eqref{084} with LSOF \eqref{086}. As mentioned in Introduction, for matrix polynomials, the problem of assigning arbitrary matrix coefficients is not equivalent to the problem of assigning a matrix spectrum. The questions of the relationship between the coefficients of a matrix polynomial and its roots (solvents), and factorization of $\lambda$-matrices have been studied in \cite{Dennis_1976} (see \cite{Shieh_1981} as well). Let an $n$th degree $s$th order monic  $\lambda$-matrix be given:
\begin{equation}\label{10-100}
\Psi(\lambda)=I \lambda^n+ \Gamma_1 \lambda^{n-1}+\ldots + \Gamma_{n-1} \lambda + \Gamma_n, \quad I\in M_s.
\end{equation}
The associated {\it left} matrix polynomial is given by
\begin{equation*}%\label{10-110}
\Psi_L(X)=X^n+ X^{n-1}\Gamma_1 +\ldots + X\Gamma_{n-1} + \Gamma_n, \quad X\in M_s(\mathbb{C}).
\end{equation*}
Let $\lambda_i$ be a complex number such that $\det \big(\Psi(\lambda_i)\big)=0$, then $\lambda_i$ is called a {\it latent root} of $\Psi(\lambda)$. The lambda-matrix $\Psi(\lambda)$ has $ns$ latent roots \cite[p.\,832]{Dennis_1976}.

If $L_j\in M_s(\mathbb{C})$ is such that
\begin{equation*}%\label{10-120}
\Psi_L(L_j)=L_j^n+ L_j^{n-1}\Gamma_1 +\ldots + L_j\Gamma_{n-1} + \Gamma_n=0\in M_s(\mathbb{C}),
\end{equation*}
then $L_j$ is referred to as a {\it left} solvent of the $\lambda$-matrix $\Psi(\lambda)$. The solvents play an important role in the spectral decomposition of $\lambda$-matrices. However, unlike the
roots of a scalar polynomial, the solvents of a $\lambda$-matrix may not exist (see the example in Introduction).

If the eigenvalues of left solvents $L_j$, $j=1,\ldots,n$, contain all the latent roots of $\Psi(\lambda)$, counting multiplicity, then the $L_j$ are referred to as a {\it complete set} of left solvents. If the elementary divisors of $\Psi(\lambda)$ are linear, then $\Psi(\lambda)$ has a complete set of left solvents \cite[Theorem 2.4]{Shieh_1981}, and $\Psi(\lambda)$ can be factored into
the product of $n$ monic linear $\lambda$-matrices (called a complete set of linear
spectral factors): $\Psi(\lambda)=(\lambda I-S_n)(\lambda I-S_{n-1})\cdots(\lambda I-S_1)$, $I\in M_s$ \cite[Theorem~2.5]{Shieh_1981}. Notice that the spectral factors are not unique. Transformations of left (right) solvents to spectral factors and  vice versa are presented in \cite{Shieh_1981}.

As a generalization of the problem of assigning eigenvalues spectrum from the scalar case to the block matrix case can be considered the problem of assigning left (or right) solvents of the $\lambda$-matrix $\Psi(\lambda)$. However, there exist sets containing $n$ matrices which are not a set of left solvents for any monic matrix polynomial of degree $n$ \cite[Corollary~6.1]{Dennis_1976}, that is, additional restrictions are required to impose on the set $\{{L}_j,\; j=1,\ldots,n\}$. They are as follows.

Define a block Vandermonde matrix as
$V(X_1,\ldots,X_n):=\{X^{i-1}_j\}_{i,j=1}^n\in M_{ns}(\mathbb{C})$, $X_j\in M_{s}(\mathbb{C})$, $j=\overline{1,n}$. Let $\{\lambda_i,\; i=1,\ldots,ns\}$ be the set of specified closed-loop eigenvalues of $\Psi(\lambda)$ and $\{\widehat{L}_j,\; j=1,\ldots,n\}$ be a set of $s\times s$-matrices such that \cite[(15)]{Leang_Shieh_1983}
\begin{equation}\label{10-125}
\bigcup_{j=1}^n \sigma(\widehat{L}_j)=\{\lambda_i,\; i=1,\ldots,ns\}
\end{equation}
and $V(\widehat{L}_1,\ldots,\widehat{L}_n)$ is nonsingular. Then $\{\widehat{L}_j,\; j=1,\ldots,n\}$ can be chosen as a complete set of left solvents of the closed-loop right characteristic $\lambda$-matrix $\Psi(\lambda)$ in \eqref{10-100}, and 
\begin{equation*}%\label{10-130}
\widehat L_j^n+ \widehat L_j^{n-1}\Gamma_1 +\ldots + \widehat L_j\Gamma_{n-1} + \Gamma_n=0\in M_s(\mathbb{C}), \quad
j=1,\ldots,n.
\end{equation*}
The matrix coefficients of $\Psi(\lambda)$ can be determined by
\begin{equation}\label{10-130}
\left[\begin{matrix}
\Gamma_n \\
\Gamma_{n-1} \\
\vdots \\
\Gamma_1
\end{matrix}
\right]
=- \left(V^{\mathcal{T}}(\widehat{L}_1,\ldots,\widehat{L}_n)\right)^{-1}
\left[\begin{matrix}
\widehat L_1^n \\
\widehat L_2^n \\
\vdots \\
\widehat L_n^n \\
\end{matrix}
\right].
\end{equation}
Thus, we have the following result.
\begin{theorem}\label{teo-10-001}
Suppose that, for system \eqref{080}, \eqref{082},  with coefficients \eqref{400}, \eqref{410}, \eqref{420}, condition \eqref{428} holds. Then, for an arbitrary set  $\{\lambda_i,\; i=1,\ldots,ns\}$ of numbers and arbitrary set  $\{\widehat{L}_j,\; j=1,\ldots,n\}$ of $s\times s$-matrices such that \eqref{10-125} holds
and $V(\widehat{L}_1,\ldots,\widehat{L}_n)$ is nonsingular, there exists LSOF \eqref{086} such that the matrix $F+GQH$ of the closed-loop system \eqref{088} is similar to the matrix \eqref{064} whose  characteristic matrix polynomial \eqref{10-100} has the prescribed left solvents $\{\widehat{L}_j,\; j=1,\ldots,n\}$.
\end{theorem}

Theorem \ref{teo-10-001} partly extends Theorem 1 of \cite{Leang_Shieh_1983} from systems with LSSF to systems with LSOF.

\section{Examples \label{sec-examples}}

\begin{example}
\label{exam-Hess}
Consider an example illustrating Theorem \ref{teo-Hessenberg-2}. Consider system \eqref{1500}, \eqref{1510} with $n=3$, $s=2$, $m=k=p=2$ with the following matrices:
\begin{equation}\label{1720}
\widetilde F=
\left[\begin{matrix}
1 & 0 & -1 & 0 & 0 & 0\\
0 & 1 & 0 & -1 & 0 & 0\\
0 & 0 & -1 & 0 & 1 & 0\\
0 & 0 & 0 & -1 & 0 & 1\\
1 & -1 & -2 & 1 & -1 & 0\\
0 & 2 & 1 & -2 & -1 & 1
\end{matrix} \right]\!, \;
\widetilde{G}=\begin{bmatrix} 0 & 0 & 0 & 0\\
0 & 0 & 0 & 0\\
-1 & 0 & 1 & 0 \\
1 & -1 & -1 & 1 \\
1 & 1 & -1 & 0 \\
0 & -1 & 0 & 1
\end{bmatrix}\!,  \;
\widetilde{H}=\begin{bmatrix} -1 & -1 & 1 & 0 & 0 & 0 \\
0 & -1 & -1 & 1 & 0 & 0 \\
1 & -1 & -1 & 1 & 0 & 0 \\
1 & -1 & 0 & 0 & 0 & 0
\end{bmatrix}\!.
\end{equation}
The matrices \eqref{1720} have the form \eqref{1520}, \eqref{1530}, \eqref{1540}. The first two block rows of the matrix 
$\widetilde F$ have the following blocks:
\begin{equation*}
\widetilde F_{11}=I, \quad
\widetilde F_{12}=-I, \quad
\widetilde F_{13}=0, \quad
\widetilde F_{21}=0, \quad
\widetilde F_{22}=-I, \quad
\widetilde F_{23}=I.
\end{equation*}
These matrices are scalar matrices. So, the first condition of Theorem \ref{teo-Hessenberg-2} is fulfilled. Next, from \eqref{1720}, we have
\begin{equation}\label{1730}
\widetilde F \widetilde{G} =
\left[\begin{matrix}
1 & 0 & -1 & 0 \\
-1 & 1 & 1 & -1 \\
2 & 1 & -2 & 0 \\
-1 & 0 & 1 & 0 \\
2 & -2 & -2 & 1 \\
-4 & 0 & 4 & -1
\end{matrix} \right], \qquad
\widetilde F^2 \widetilde{G}=\begin{bmatrix} 
-1 & -1 & 1 & 0 \\
0 & 1 & 0 & -1 \\
0 & -3 & 0 & 1 \\
-3 & 0 & 3 & -1 \\
-5 & -1 & 5 & 0 \\
-4 & 5 & 4 & -4
\end{bmatrix}.
\end{equation}
From \eqref{1720} and \eqref{1730}, we obtain:
\begin{gather*}
(\widetilde H^{\mathcal{T}})^T\star \widetilde{G}=\left[\begin{matrix}X_{11}&X_{12}\\X_{21}&X_{22}\end{matrix}\right],\\
X_{11}=\left[\begin{matrix}-1&0&1&0\\1&-1&-1&1\\0&0&-1&0\\0&0&1&-1\end{matrix}\right],\quad 
X_{12}=\left[\begin{matrix}1&0&-1&0\\-1&1&1&-1\\0&0&1&0\\0&0&-1&1\end{matrix}\right],\\
X_{21}=\left[\begin{matrix}1&0&0&0\\-1&1&0&0\\-1&0&0&0\\1&-1&0&0\end{matrix}\right],\quad 
X_{22}=\left[\begin{matrix}-1&0&0&0\\1&-1&0&0\\1&0&0&0\\-1&1&0&0\end{matrix}\right],
\end{gather*}
\begin{gather*}
(\widetilde H^{\mathcal{T}})^T\star \widetilde{F}\widetilde{G}=\left[\begin{matrix}Y_{11}&Y_{12}\\Y_{21}&Y_{22}\end{matrix}\right],\\
Y_{11}=\left[\begin{matrix}1&1&-2&-1\\0&-1&1&0\\-1&0&1&1\\1&-1&0&-1\end{matrix}\right],\quad 
Y_{12}=\left[\begin{matrix}-1&0&2&0\\0&1&-1&0\\1&0&-1&0\\-1&1&0&1\end{matrix}\right],\\
Y_{21}=\left[\begin{matrix}-1&-1&1&0\\0&1&-1&1\\1&1&-1&0\\0&-1&1&-1\end{matrix}\right],\quad 
Y_{22}=\left[\begin{matrix}1&0&-1&0\\0&-1&1&-1\\-1&0&1&0\\0&1&-1&1\end{matrix}\right],
\end{gather*}
\begin{gather*}
(\widetilde H^{\mathcal{T}})^T\star \widetilde{F}^2\widetilde{G}=\left[\begin{matrix}Z_{11}&Z_{12}\\Z_{21}&Z_{22}\end{matrix}\right],\\
Z_{11}=\left[\begin{matrix}1&-2&0&3\\-3&-1&3&0\\1&1&1&-2\\0&-1&-3&-1\end{matrix}\right],\quad 
Z_{12}=\left[\begin{matrix}-1&1&0&-1\\3&0&-3&1\\-1&0&-1&1\\0&1&3&0\end{matrix}\right],\\
Z_{21}=\left[\begin{matrix}-1&2&-1&-1\\3&1&0&1\\1&-2&1&1\\-3&-1&0&-1\end{matrix}\right],\quad 
Z_{22}=\left[\begin{matrix}1&-1&1&0\\-3&0&0&-1\\-1&1&-1&0\\3&0&0&1\end{matrix}\right].
\end{gather*}
Constructing matrix \eqref{1620}, we get
\begin{equation}\label{1740}
\widehat\Theta=\left[\begin{matrix}
X_{11}&X_{12}&X_{21}&X_{22}\\Y_{11}&Y_{12}&Y_{21}&Y_{22}\\Z_{11}&Z_{12}&Z_{21}&Z_{22}\end{matrix}\right].
\end{equation}
Calculating the rank of matrix \eqref{1740}, we obtain that ${\rm rank}\,\widehat\Theta=12$. So, conditions of Theorem \ref{teo-Hessenberg-2} are fulfilled. Hence, by Theorem \ref{teo-Hessenberg-2}, for system \eqref{1500}, \eqref{1510} with coefficients \eqref{1720}, the problem of AMCA for CMP by LSOF is resolvable. Let us construct this feedback control. Suppose, for example, that
\begin{equation}\label{1750}
\Gamma_1 =
\left[\begin{matrix}
-1 & 0 \\
0 & -1
\end{matrix} \right], \qquad
\Gamma_2 =
\left[\begin{matrix}
1 & 0 \\
0 & 0
\end{matrix} \right], \qquad
\Gamma_3 =
\left[\begin{matrix}
-2 & 1 \\
0 & 1
\end{matrix} \right].
\end{equation}
First, we reduce the system  \eqref{1500}, \eqref{1510}, \eqref{1720} with a lower block
Hessenberg matrix to a system \eqref{080}, \eqref{082} with coefficients \eqref{400}, \eqref{410}, \eqref{420}. By Lemma \ref{lemma-00}, we set 
\begin{equation*}%\label{1760}
\widetilde S_1:=
\left[\begin{matrix}
I & 0 & 0 \\
I & -I & 0 \\
0 & -I & I
\end{matrix} \right], \quad
\widetilde S_2:=
\left[\begin{matrix}
I & 0 & 0 \\
0 & I & 0 \\
0 & I & -I
\end{matrix} \right].
\end{equation*}
Construct $\widetilde S:=\widetilde S_2 \cdot \widetilde S_1$. Then,
\begin{equation*}%\label{1770}
\widetilde S=
\left[\begin{matrix}
I & 0 & 0 \\
I & -I & 0 \\
I & 0 & -I
\end{matrix} \right].
\end{equation*}
Construct matrices \eqref{1560}. We get 
\begin{gather}\label{1780}
F=\left[\begin{matrix}0&I&0\\0&0&I\\-A_3&-A_2&-A_1\end{matrix}\right],\quad 
A_1=\left[\begin{matrix}1&0\\1&-1\end{matrix}\right],\quad A_2=\left[\begin{matrix}1&-1\\-1&1\end{matrix}\right],\quad A_3=\left[\begin{matrix}-2&0\\0&1\end{matrix}\right],\\
\notag
G=\left[\begin{matrix}0&0\\G_{21}&G_{22}\\G_{31}&G_{32}\\\end{matrix}\right],\quad G_{21}=\left[\begin{matrix}1&0\\-1&1\end{matrix}\right],\quad G_{22}=\left[\begin{matrix}-1&0\\1&-1\end{matrix}\right], \quad
G_{31}=\left[\begin{matrix}-1&-1\\0&1\end{matrix}\right],\quad G_{32}=\left[\begin{matrix}1&0\\0&-1\end{matrix}\right],\\
\notag
H=\left[\begin{matrix}H_{11}&H_{12}&0\\H_{21}&H_{22}&0\end{matrix}\right], \;
H_{11}=\left[\begin{matrix}0&-1\\-1&0\end{matrix}\right],\; H_{12}=\left[\begin{matrix}-1&0\\1&-1\end{matrix}\right],\;
H_{21}=\left[\begin{matrix}0&0\\1&-1\end{matrix}\right],\; H_{22}=\left[\begin{matrix}1&-1\\0&0\end{matrix}\right].
\end{gather}
Construct the matrix \eqref{425}. By Lemma \ref{lem-Hess-S}, we obtain that $\Theta$ is equal to $\widehat\Theta$ of \eqref{1740}; therefore, ${\rm rank}\,\Theta=12$. By \eqref{430}, we construct matrix $P=\left[\begin{matrix}I&0&0\\A_1&I&0\\A_2&A_1&I\end{matrix}\right]$  where $A_1$, $A_2$ are from \eqref{1780}. Substituting $A_1$, $A_2$, $A_3$ from \eqref{1780} and $\Gamma_1$, $\Gamma_2$, $\Gamma_3$  from \eqref{1750} and $P$ into \eqref{680}, we obtain that
\begin{equation*}%\label{1810}
\widehat{T}=\left[\begin{matrix} T_1 \\ T_2 \\ T_3 \end{matrix}\right], \quad
T_1=\left[\begin{matrix}2&0\\1&0\end{matrix}\right],\quad T_2=\left[\begin{matrix}-2&-1\\-2&1\end{matrix}\right],\quad T_3=\left[\begin{matrix}1&0\\1&2\end{matrix}\right].
\end{equation*}
Let us construct $w$ by formula \eqref{710}. Then
$w={\rm col}\,(2,1,0,0,-2,-2,-1,1,1,1,0,2)$. Now, resolving the system $\Theta v =w$ of \eqref{720} by the formula \eqref{730}, we obtain that 
$$
v={\rm col}\,(-1,-1,1,0,1,2,-1,-7,-1,-1,1,0,1,6,-1,-2).
$$ From \eqref{740}, we get
\begin{equation}\label{1820}
Q=\left[\begin{matrix}-1&1&-1&1\\-1&0&-1&0\\1&-1&1&-1\\2&-7&6&-2\end{matrix}\right].
\end{equation}
We have
\begin{equation*}
F + G  Q  H=
\left[\begin{matrix}
0&0&1&0&0&0\\
0&0&0&1&0&0\\
0&0&2&0&1&0\\
-5&1&1&0&0&1\\
2&-1&-3&0&-1&0\\
-5&0&4&-1&-1&1
\end{matrix} \right].
\end{equation*}
Applying step 10 of Algorithm 1, we obtain that 
\begin{equation*}
S=\left[\begin{matrix} I  & 0  &  0\\  0 & I & 0  \\ S_{31} & S_{32} & I
\end{matrix}\right], \quad
S_{31}=\left[\begin{matrix} 0  & 0  \\  -5 &  1  
\end{matrix}\right], \quad
S_{32}=\left[\begin{matrix} 2  & 0  \\  1 &  0  
\end{matrix}\right].
\end{equation*}
Feedback control $\widetilde u =Q \widetilde y$ with \eqref{1820} leads the system \eqref{1500}, \eqref{1510}, \eqref{1720} to the closed-loop system \eqref{1585} with the matrix
\begin{equation}\label{1830}
\widetilde F + \widetilde G  Q \widetilde H=
\left[\begin{matrix}
1 & 0 & -1 & 0 & 0 & 0\\
0 & 1 & 0 & -1 & 0 & 0\\
-2 & 0 & 1 & 0 & 1 & 0\\
4 & -1 & 1 & -1 & 0 & 1\\
3 & 1 & -4 & 0 & -1 & 0\\
2 & 1 & 4 & -2 & -1 & 1
\end{matrix} \right].
\end{equation}
Construct $\mathcal{R}:=S\widetilde{S}$. Then,
\begin{equation*}
\mathcal{R}=
\left[\begin{matrix}
1 & 0 & 0 & 0 & 0 & 0\\
0 & 1 & 0 & 0 & 0 & 0\\
1 & 0 & -1 & 0 & 0 & 0\\
0 & 1 & 0 & -1 & 0 & 0\\
3 & 0 & -2 & 0 & -1 & 0\\
-4 & 2 & -1 & 0 & 0 & -1
\end{matrix} \right].
\end{equation*}
Then, the matrix \eqref{1830} is similar to the matrix $\Phi=\left[\begin{matrix}0&I&0\\0&0&I\\-\Gamma_3&-\Gamma_2&-\Gamma_1\end{matrix}\right]$ by means of the matrix $\mathcal{R}$: one can check that $\Phi =\mathcal{R}\big(\widetilde F + \widetilde G  Q \widetilde H){\mathcal{R}}^{-1}$.
\end{example}

\if0
. Let us check it. Matrix \eqref{1830} is a lower block Hessenberg matrix. By applying Lemma \ref{lemma-00} to this matrix, we construct
\begin{equation*}
\widetilde{\mathcal{S}}_1:=
\left[\begin{matrix}
1 & 0 & 0 & 0 & 0 & 0\\
0 & 1 & 0 & 0 & 0 & 0\\
1 & 0 & -1 & 0 & 0 & 0\\
0 & 1 & 0 & -1 & 0 & 0\\
-2 & 0 & 1 & 0 & 1 & 0\\
4 & -1 & 1 & -1 & 0 & 1
\end{matrix} \right], \quad
\widetilde{\mathcal{S}}_2:=
\left[\begin{matrix}
I & 0 & 0 \\
0 & I & 0 \\
0 & I & -I
\end{matrix} \right],
\end{equation*}
and $\widetilde{\mathcal{S}}:=\widetilde{\mathcal{S}}_2\widetilde{\mathcal{S}}_1$.
Then
\begin{equation*}%\label{1850}
\widetilde{\mathcal{S}}=
\left[\begin{matrix}
1 & 0 & 0 & 0 & 0 & 0\\
0 & 1 & 0 & 0 & 0 & 0\\
1 & 0 & -1 & 0 & 0 & 0\\
0 & 1 & 0 & -1 & 0 & 0\\
3 & 0 & -2 & 0 & -1 & 0\\
-4 & 2 & -1 & 0 & 0 & -1
\end{matrix} \right].
\end{equation*}
The matrices \eqref{1830} and $\Phi$ are similar by means of $\widetilde{\mathcal{S}}$: 
one can check that $\Phi =\widetilde{\mathcal{S}}\big(\widetilde F + \widetilde G  Q \widetilde H)\widetilde{\mathcal{S}}^{-1}$.
\fi

\begin{example}\label{example2}
Let us give an example illustrating that if not all blocks of the matrix $\widetilde S$ are scalar matrices, then the conclusion of Lemma \ref{lem-Hess-S} is not true, in general case. Let $n=4$, $s=2$, $m=2$, $k=2$, $p=2$,
\begin{equation*}
\widetilde F=
\left[\begin{matrix}
0 & A & 0 & 0\\
 0 & 0 & B & 0 \\
 0 & 0 & 0 & I \\
 0 & 0 & B & 0 
\end{matrix} \right], \;
\widetilde G=
\left[\begin{matrix}
0 &  0\\
 I & 0  \\
 0 & I \\
 I & 0
\end{matrix} \right], \;
\widetilde H=
\left[\begin{matrix}
I &  0 & 0 & 0 \\
0 & I & 0 & 0 
\end{matrix} \right],
 \;
A=
\left[\begin{matrix}
2 &  0 \\
0 & 1 
\end{matrix} \right],
\;
B=
\left[\begin{matrix}
0 &  1 \\
1 & 0 
\end{matrix} \right],\; I,0\in M_2.
\end{equation*}
Then, $\widetilde{S}_1= \mathrm{diag} \{I,A,B,I\}$, 
$\widetilde{S}_2= \mathrm{diag} \{I,I,A,B\}$, 
$\widetilde{S}_3= \mathrm{diag} \{I,I,I,A\}$. Hence, 
$\widetilde{S}=\widetilde{S}_3 \widetilde{S}_2 \widetilde{S}_1= \mathrm{diag} \{I,A,AB,AB\}$.
By constructing the matrix $\widehat\Theta$ of \eqref{1620}, we get that ${\rm rank}\,\widehat\Theta = 16$.
By constructing the matrix $\widetilde\Theta$ of \eqref{1570}, we get that ${\rm rank}\,\widetilde\Theta = 12$.
So, condition \eqref{1580} is not equivalent to condition \eqref{1630}, and $\widetilde\Theta \ne \widehat\Theta$.
\end{example}

\begin{example}\label{example3}
Let us show that the converse assertion to Theorem \ref{teo-AMESA-implies-AESA} is not true. Let $n=2$, $s=2$, $m=1$, $k=2$, $p=2$; $F=
\left[\begin{matrix}
0 & I \\
 0 & 0 
\end{matrix} \right]\in M_4$, $I,0\in M_2$; $G=\left[\begin{matrix}
0  \\
 I
\end{matrix} \right]\in M_{4,2}$, $I,0\in M_2$; $H=
\left[\begin{matrix}
H_{11} & H_{12} \\
 H_{21} & H_{22} 
\end{matrix} \right]\in M_4$, $H_{11}=I\in M_2$, $H_{12}=H_{21}=0\in M_2$, $H_{22}=\left[\begin{matrix}
1 & 0 \\
 0 & 0
\end{matrix} \right] \in M_2$.
Denote $Q=\{q_{ij}\}\in M_{2,4}$, $q_{ij}\in\mathbb{K}$, $i=1,2$, $j=\overline{1,4}$. Then,
\begin{equation}\label{9-300}
F+GQH=\left[\begin{matrix}
0 & 0 & 1 & 0\\
0 & 0 & 0 & 1\\
q_{11} & q_{12} & q_{13} & 0 \\
q_{21} & q_{22} & q_{23} & 0 
\end{matrix} \right].
\end{equation}
For any $\delta_i\in\mathbb{K}$, $i=\overline{1,4}$, let us take 
\begin{equation}\label{9-310}
q_{11}:=-\delta_2, \quad q_{12}:=1, \quad q_{13}:=-\delta_1, \quad q_{21}:=-\delta_4, \quad q_{22}:=0, \quad q_{23}:=-\delta_3.
\end{equation}
Calculating the characteristic polynomial of the matrix \eqref{9-300} with coefficients \eqref{9-310}, we obtain that $\chi(F+GQH, \lambda)= \lambda^4+\delta_1 \lambda^3+\delta_2\lambda^2+\delta_3\lambda+\delta_4$. Thus,
for system \eqref{080}, \eqref{082}, \eqref{084}, the problem of ASCA for CP by LSOF is resolvable.

Let us show that not for any matrices $\Gamma_1, \Gamma_2 \in M_2$ there exists $Q$ such that \eqref{9-300} is similar to the matrix 
\begin{equation}\label{9-320}
\Phi=\left[\begin{matrix}
0 & I\\
-\Gamma_2 & -\Gamma_1 
\end{matrix} \right].
\end{equation}
Set
\begin{equation}\label{9-330}
\Gamma_1=\left[\begin{matrix}
-2a & 0\\
0 & -2a
\end{matrix} \right], \quad
\Gamma_2=\left[\begin{matrix}
a^2 & 0\\
0 & a^2
\end{matrix} \right], \quad \text { where } \; a\ne 0.
\end{equation}
Then
\begin{equation}\label{9-340}
\chi(\Phi,\lambda)=\lambda^4-4a\lambda^3+6a^2\lambda^2-4a^3\lambda + a^4.
\end{equation}
Let us assume the contrary: there exists $Q$ such that \eqref{9-300} is similar to \eqref{9-320} with \eqref{9-330}. Then, the following conditions are necessarily satisfied:
\begin{gather}\label{9-350}
\chi(F+GQH,\lambda)=\chi(\Phi,\lambda), \\
\label{9-360}
{\rm rank}\,(F+GQH- aI)={\rm rank}\,(\Phi- aI).
\end{gather}
From \eqref{9-300}, we have
\begin{equation}\label{9-370}
\chi(F+GQH,\lambda)=\lambda^4-q_{13} \lambda^3 - (q_{11}+q_{22})\lambda^2 + (q_{13}q_{22}-q_{12}q_{23})\lambda+(q_{11}q_{22}-q_{12}q_{21}).
\end{equation}
From \eqref{9-340}, \eqref{9-370} and equality \eqref{9-350}, it follows that, in particular,
\begin{gather}\label{9-380}
q_{13}=4a, \\
\label{9-390}
q_{11}+q_{22}=-6a^2.
\end{gather}
Next, we have ${\rm rank}\,(\Phi- aI)=2$. Substituting \eqref{9-380} into \eqref{9-300} and constructing $F+GQH- aI$, we obtain
\begin{equation}\label{9-400}
F+GQH - a I=\left[\begin{matrix}
-a & 0 & 1 & 0\\
0 & -a & 0 & 1\\
q_{11} & q_{12} & 3 a & 0 \\
q_{21} & q_{22} & q_{23} & -a  
\end{matrix} \right].
\end{equation}
By \eqref{9-360}, the rank of the matrix \eqref{9-400} should be equal to 2. 
Consequently, any third-order minor of matrix \eqref{9-400} must necessarily be equal to zero.
Taking the minor of matrix \eqref{9-400} consisting of 1, 2, 3 rows and 1, 2, 3 columns and equating to zero, we obtain that $q_{11}=-3a^2$. 
Taking the minor of matrix \eqref{9-400} consisting of 1, 2, 4 rows and 2, 3, 4 columns and equating to zero, we obtain that $q_{22}=a^2$.  Then, $q_{11}+q_{22}=-2a^2$. This contradicts to \eqref{9-390} because $a\ne 0$. Q.E.D.
\end{example}

\begin{example}
\label{ex-10}
Let us illustrate  Theorem \ref{teo-10-001} using the example of the system \eqref{080}, \eqref{082} with matrices \eqref{ex1-1720}. Let the set of eigenvalues be given: $\{-1,-1,-2,-2,-3,-3\}$. Take
\begin{equation}\label{15-010}
\widehat{L}_1=\left[\begin{matrix}
-1 & 0 \\
0 & -1 
\end{matrix} \right], \quad
\widehat{L}_2=\left[\begin{matrix}
-2 & 0 \\
0 & -2 
\end{matrix} \right], \quad
\widehat{L}_1=\left[\begin{matrix}
-3 & 0 \\
0 & -3
\end{matrix} \right].
\end{equation}
Then \eqref{10-125} holds. Construct the block Vandermonde matrix:
\begin{equation*}%\label{15-020}
V(\widehat{L}_1,\widehat{L}_2,\widehat{L}_3)=\left[\begin{matrix}
I & I & I \\
\widehat{L}_1 & \widehat{L}_2 & \widehat{L}_3 \\
\widehat{L}_1^2 & \widehat{L}_2^2 & \widehat{L}_3^2
\end{matrix} \right] = \left[\begin{matrix}
I & I & I \\
-I & -2 I & - 3 I \\
I & 4 I & 9 I
\end{matrix} \right].
\end{equation*}
It is clear that the matrix $V(\widehat{L}_1,\widehat{L}_2,\widehat{L}_3)$ is non-singular.
Thus, one can assign the prescribed left solvents \eqref{15-010}. To do this, calculate the coefficients of the characteristic matrix using formula \eqref{10-130}, we obtain that 
\begin{equation}\label{15-210}
\Gamma_1 =
\left[\begin{matrix}
6 & 0 \\
0 & 6
\end{matrix} \right], \qquad
\Gamma_2 =
\left[\begin{matrix}
11 & 0 \\
0 & 11
\end{matrix} \right], \qquad
\Gamma_3 =
\left[\begin{matrix}
6 & 0 \\
0 & 6
\end{matrix} \right].
\end{equation}
In Example \ref{ex0}, the matrix \eqref{ex1-1820}  was constructed that provides the similarity of the matrix $F+GQH$ (with coefficients \eqref{ex1-1720}) to matrix \eqref{064} with coefficients \eqref{15-210} and hence, the prescribed left solvents \eqref{15-010}. 
\end{example}

\section{Conclusion \label{sect-conclusion}}

In this work,  we have introduced the formulation of the problem of AMCA  for CMP by LSSF and LSOF. This problem is a generalization of the problem of assigning the scalar spectrum. Such problems have applications in problems of assigning an eigenstructure, in stabilization problems and other problems (see Sect.~\ref{sect1}). The following main results have been obtained.

(1) Sufficient conditions have been obtained for resolving the problem of AMCA for CMP by LSOF
when the state matrix is a lower block Frobenius matrix, and the input and output block matrix coefficients contain some zero blocks.

(2) It is proved that these sufficient conditions generalize the sufficient conditions obtained 
earlier in Theorems 2 and 3 of \cite {Zaitsev_2021} for system \eqref{210}, \eqref{220}.

(3) Corollaries are derived from Theorem \ref{teo-AMCA-LSOF-1}, when some matrix coefficients of system \eqref{080}, \eqref{082}, \eqref{084} have a simpler form, in particular, their blocks are scalar matrices.

(4) The sufficient conditions of Theorem \ref{teo-AMCA-LSOF-1} are extended to the case when the state matrix is a lower block Hessenberg matrix.

(5) Sufficient conditions are obtained for the solvability of the problem of assigning the matrix spectrum of solvents.

All results are of a constructive nature. Based on the proofs, algorithms for solving the problems have been constructed. These algorithms have been demonstrated on a number of computational examples. The examples have been performed on a computer using Maple program, LinearAlgebra package.

Previously, other works \cite{Vitoria_1982,Victoria_1982-2,Martins_2014} also considered block systems and, in particular, issues of stability of such systems. In these works, a restriction was imposed: matrix blocks must commute pairwise. The advantage of our work is that this restriction is not imposed in our work.

Further, the problem under study exhibits complications compared to the case $s=1$. For example, the necessity in almost all theorems remains unproven yet. One of the technical obstacles to this is the following fact. For $s=1$, if two lower Frobenius matrices $F_1$ and $F_2$ are similar, then they necessarily coincide. For $s>1$, if two lower block Frobenius matrix $\mathcal{F}_1$ and $\mathcal{F}_2$ are similar, then they do not necessarily coincide (see \cite[Sect. 1]{Zaitsev_2024-2}).

There are other open questions. For example, one of the important open question is the analogs of the statements (St1) and (St2) (see Remark \ref{rem-dop-2}) for the case $m>1$. Another question is whether it is possible to replace arbitrary matrices in Definitions \ref{def-01} and \ref{def-02} with triangular (or diagonal) ones. Partial answers to these questions have been given in \cite{Zaitsev_2024-2}.
We plan to study open questions in future work.

%\section*{Acknowledgments}

%The research was supported by the grant from the Russian Science Foundation No.~24-21-00311, https://rscf.ru/project/24-21-00311/ .

% You may incorporate your references as follows in your main tex file.
% Using BibTex is not recommended but can be handled.

%\bibliographystyle{AIMS}
%\bibliography{ref}

\end{document}